\documentclass[fontsize=12pt,a4paper,headings=normal,
twoside=false,leqno,parskip=half-,abstract=true]{scrartcl}
\usepackage[english]{babel}
\usepackage[utf8]{inputenc}
\usepackage[nottoc,notlot,notlof]{tocbibind} %makes References appear in toc
\usepackage{csquotes}
\usepackage{hyperref}
\hypersetup{
 pdftitle={eternal is not entire},
 pdfauthor={Bernold Fiedler},
 colorlinks=true,
 linkcolor=blue,
 citecolor=blue,
 filecolor=blue,
 urlcolor=blue}

 \usepackage{afterpage}

\usepackage{graphicx}
\usepackage{wrapfig} %to position figures left, right, etc of text; https://www.overleaf.com/learn/latex/Positioning_images_and_tables
\usepackage[format=plain,labelfont=bf,font=small]{caption}
\usepackage{xcolor}
\usepackage[arrow, matrix, curve]{xy}
\usepackage{float}

\usepackage{caption}
\usepackage{tabulary}
\usepackage{array}
\newcolumntype{N}[1]{>{\centering\arraybackslash}m{#1}}

\usepackage{amsmath,amsthm}
\swapnumbers % Switch number/label style
\usepackage{amssymb, eurosym} 
\makeatletter
\newcommand{\tpitchfork}{%
  \vbox{
    \baselineskip\z@skip
    \lineskip-.52ex
    \lineskiplimit\maxdimen
    \m@th
    \ialign{##\crcr\hidewidth\smash{$-$}\hidewidth\crcr$\pitchfork$\crcr}
  }%
}
\makeatother
\usepackage{latexsym}
\usepackage{enumerate}% http://ctan.org/pkg/enumerate

\usepackage[notref,notcite,color, final %make final instead of draft, by % %
]{showkeys}
\definecolor{refkey}{rgb}{1,0,0}
\definecolor{labelkey}{rgb}{1,0,0}

\usepackage{cancel}

\usepackage{tikz}

  \mathchardef\ordinarycolon\mathcode`\:
  \mathcode`\:=\string"8000
  \begingroup \catcode`\:=\active
    \gdef:{\mathrel{\mathop\ordinarycolon}}
  \endgroup

\theoremstyle{plain}
\newtheorem{thm}{Theorem}[section]
\newtheorem{lem}[thm]{Lemma}

\newtheorem{cor}[thm]{Corollary}

\newcommand\eps{\varepsilon}
\newcommand\mi{\mathrm{i}}
\renewcommand\theta{\vartheta}
\renewcommand\rho{\varrho}
\renewcommand\phi{\varphi}
\renewcommand\Re{\mathrm{Re}\,}
\renewcommand\Im{\mathrm{Im}\,}

\begin{document}

\title{\LARGE{Real eternal PDE solutions\\
are not complex entire:\\
a quadratic parabolic example}}
\vspace{1cm}
{\subtitle{}
	\vspace{1ex}
	{}}\vspace{1ex}

\author{
 \\
\emph{-- Dedicated to the memory of Marek Fila  --}\\
{~}\\
Bernold Fiedler* and Hannes Stuke*\\
\vspace{2cm}}

\date{\small{version of \today}}
\maketitle
\thispagestyle{empty}

\vfill

*\\
Institut für Mathematik\\
Freie Universität Berlin\\
Arnimallee 3\\ 
14195 Berlin, Germany

%%%%%%%%%%%%%%%%%%%%%%%%%%%%%%%%%%%%%%%%%%%%%%%%%%%%%%%%%%%

\newpage
\pagestyle{plain}
\pagenumbering{roman}
\setcounter{page}{1}

\begin{abstract}
\noindent
In parabolic or hyperbolic PDEs, solutions which remain uniformly bounded for all real times $t=r\in\mathbb{R}$ are often called \emph{PDE entire} or \emph{eternal}.
For a nonlinear example, consider the quadratic parabolic PDE
\begin{equation*}
\label{*}
w_t=w_{xx}+6w^2-\lambda\,,  \tag{*}
\end{equation*}
for $0<x<\tfrac{1}{2}$, under Neumann boundary conditions. 
By its gradient-like structure, all \emph{real eternal} non-equilibrium orbits $\Gamma(r)$ of \eqref{*}  are heteroclinic among equilibria $w=W_n(x)$.
For parameters $\lambda>0$, the trivial homogeneous equilibria are locally asymptotically stable $W_0=-\sqrt{\lambda/6}$, and $W_\infty=+\sqrt{\lambda/6}$ of unstable dimension (Morse index) $i(W_\infty)=1,2,3,\ldots$\,, depending on $\lambda$.
All nontrivial real $W_n$ are rescaled and properly translated real-valued Weierstrass elliptic functions with Morse index $i(W_n)=n$\,.

\smallskip\noindent
We show that the complex time extensions $\Gamma(r+\mathrm{i}s)$, of analytic real heteroclinic orbits $\Gamma(r)$ towards $W_0$\,, are \emph{not complex entire}.
For example, consider the time-reversible complex-valued solution $\psi(s)=\Gamma(r_0-\mi s)$ of the nonlinear and nonconservative quadratic Schrödinger equation
\begin{equation*}
\label{**}
\mathrm{i}\psi_s=\psi_{xx}+6\psi^2-\lambda  %\tag{**}
\end{equation*}
with real initial condition $\psi_0=\Gamma(r_0)$.
Then there exist real $r_0$ such that $\psi(s)$ blows up at some finite real times $\pm s^*\neq 0$\,. 

\smallskip\noindent
Abstractly, our results are formulated in the setting of analytic semigroups. 
They are based on Poincaré non-resonance of unstable eigenvalues at equilibria $W_n$\,, near pitchfork bifurcation.
Technically, we have to except discrete sets of parameters $\lambda$, and are currently limited to unstable dimensions $i(W_n)\leq22$, or to fast unstable manifolds of dimensions $d<1+\tfrac{1}{\sqrt{2}}i(W_n)$\,.

\end{abstract}

\newpage
%\vspace{2cm}
\tableofcontents

%%%%%%%%%%%%%%%%%%%%%%%%%%%%%%%%%%%%%%%%%%%%%%%%%%%%%%%%%%%

\newpage
\pagenumbering{arabic}
\setcounter{page}{1}

\section{Introduction and main results} \label{Int}

\numberwithin{equation}{section}
\numberwithin{figure}{section}
\numberwithin{table}{section}

\subsection{A paradigm}\label{Para}
A major part of the work by Marek Fila and his co-authors has been dedicated to parabolic blow-up in real time.
See for example the many references to his work in the standard monograph \cite{Quittner}.
Interestingly, he devoted several papers to real continuation after (incomplete) real-time blow-up \cite{FilaPolacik, FilaMatano, FilaMatanoPolacik, FilaMizoguchi}.
First attempts to circumnavigate complete real-time parabolic blow-up, by detours through complex time, go back to Kyûya Masuda \cite{Masuda1,Masuda2}.
In contrast, we explore blow-up in complex time, for real solutions which remain uniformly bounded, eternally, for all forward and backward real time.
See section \ref{Sing} for further discussion of some PDE literature.

Our guiding paradigm will be:
\begin{center}
\emph{Real eternal solutions are not complex entire.}
\end{center}
In PDE context, solutions $\Gamma(r)$ which exist for all real times $r$ are often called \emph{eternal} or \emph{(PDE) entire}.
Throughout our present paper, we restrict this notion further, by the additional requirement that 
\textbf{eternal solutions $\Gamma(r)$ remain uniformly bounded, for all real $t$}.
For example, this more restrictive notion characterizes global attractors for dissipative compact semigroups.
See \cite{chvi02, ha88, haetal02, la91, te88} and others, for general background on global attractors.

Here and below, the terms \emph{``analytic''} and \emph{``analyticity''} emphasize local expansions by convergent power series.
\emph{``Holomorphic''} emphasizes complex differentiability and, therefore, Cauchy-Riemann equations.
For continuously real differentiable functions, e.g., the two notions coincide.
In complex analysis, unfortunately, \emph{entire} functions are defined to be globally analytic (alias holomorphic or complex differentiable) on the whole complex plane, or on $\mathbb{C}^N$.
To avoid confusion, we talk of \emph{complex entire} functions, and we avoid the term ``PDE entire'', entirely.

Of course, there are obvious counterexamples to our paradigm: equilibria, i.e.~time-independent solutions, or real sines and cosines come to mind, among others.
Specifically, the present paper therefore explores how real-time heteroclinic solutions $\Gamma(t)$, between equilibria, of analytic semigroups give rise to blow-up behavior, in complex time $t=r+\mi s$.
Heteroclinic solutions $\Gamma:W_-\leadsto W_+$\,, of course, are examples of real eternal solutions.
In that terminology, our results on heteroclinicity will support the paradigm that real eternal solutions are not complex entire.

The study of finite time blow-up, for complex-valued PDEs and in complex time, is fraught with technical obstacles.
Not the least among them are the absence of, both, variational structure and effective comparison principles.
We will use spectral non-resonance and Poincaré linearization, to overcome some of these difficulties.

\subsection{Semigroup setting} \label{Sgr}

Our abstract setting are analytic semigroups.
To be specific, let $\Phi^t(w_0):=w(t)$ denote the local analytic semigroup defined by solutions $w=w(t)$ of semilinear equations
\begin{equation}
\label{sgr}
\dot w = Aw+f(w)
\end{equation}
with initial condition $w(0)=w_0$\,. 
See \cite{Henry, Pazy} for a general background.
Here the complex linear sectorial operator $A$ is the infinitesimal generator of a linear analytic semigroup $\exp (At)$ on a complex Banach space $X$.
The fractional power space $X^\alpha$ is equipped with the graph norm $\|\cdot\|_\alpha$ of the fractional power $(-A)^\alpha$. The nonlinearity $f:X^\alpha\rightarrow X$ is assumed to be \emph{complex entire}, i.e. complex Fréchet differentiable on all $X^\alpha$, for some fixed $0\leq\alpha<1$.
We assume locally uniform bounds
\begin{equation}
\label{f'}
\|f'(w)\|\leq C(\|w\|_\alpha)
\end{equation}
on the Fréchet derivatives $f'(w):X^\alpha\rightarrow X$. 
For given initial conditions $w=w_0$ at time $t=0$, let $\Phi^t(w_0):=w(t)$ denote the local solution of \eqref{sgr}.
Then the semigroup $\Phi^t:X^\alpha\rightarrow X^\alpha$ is well-defined and strongly continuous, for small $|t|<\rho(w_0)$ in a closed complex positive sector 
\begin{equation}
\label{sector}
\mathfrak{S}_\Theta:=\{t\in\mathbb{C}\,|\ \lvert\arg t\rvert\leq\Theta<\pi\}.
\end{equation}
Complex differentiability of $\Phi^t$ holds in the interior of $\mathfrak{S}_\Theta$\,.
Indeed this follows from complex differentiability of $f$.
For $t=r+\mi s$ in the sector, just invoke real $C^1$ differentiability with respect to $r,s$, following \cite{Henry}.
Then note that solutions satisfy the Cauchy-Riemann equations; see also \cite{ChowHale}.
Alternatively, one may directly invoke a complex version of the implicit function theorem.
For an explicit exposition based on the majorant method, see for example \cite{Berger}.

By construction and Cauchy's theorem, the local semigroup satisfies
\begin{equation}
\label{flow}
\Phi^{t_2}\circ\Phi^{t_1}= \Phi^{t_1+t_2}, \qquad \Phi^0=\mathrm{id},
\end{equation}
for any initial condition $w_0\in X^\alpha$ and, locally, for all arguments  $t_1,t_2 \in \mathfrak{S}_\Theta$ such that the (small) closed complex parallelogram spanned by $t_1,t_2$ is contained in the domain of existence of the local semiflow $\Phi^t(w_0)$.
More generally, commuting local flows $\Phi_\iota^t$ are generated by vector fields $\dot{w}=F_\iota(w),\ \iota=1,2,$ with vanishing commutator
\begin{equation}
\label{[]}
[F_1,F_2] := F_1'\cdot F_2-F_2'\cdot F_1=0\,.
\end{equation}
In our complex case, the commutator vanishes by complex linearity of $A,f'$.

For real times $t=r\in\mathbb{R}>0$, let $r\in[0,r^*(w_0))$ denote the maximal interval of existence for the solution $w(r)$ of \eqref{sgr}.
Note $0<r^*\leq\infty$.
By the locally uniform Lipschitz bound \eqref{f'}, the case $0<r^*<\infty$ of \emph{finite time blow-up} is then characterized by solutions $w(r)$ escaping to infinity:
\begin{equation}
\label{defblow-up}
\|w(r)\|_\alpha \rightarrow\infty \quad\textrm{for}\quad r\nearrow r^*<\infty\,.
\end{equation}
An analogous definition applies for complex times $t\rightarrow t^*=t^*(w_0)$, with $t$ approaching $t^*$ from any backward sector $t^*-\mathfrak{S}_\Theta$\,.
In reverse time $r\searrow r^*>-\infty$, where applicable, we sometimes speak of \emph{blow-down}.

\emph{Equilibria} $W\in X^\alpha$, i.e.~time-independent solutions of \eqref{sgr}, satisfy 
\begin{equation}
\label{W}
0=AW+f(W)\,.
\end{equation}
We call $W$ \emph{hyperbolic} if the spectrum of the linearization $L=A+f'(W)$ at $W$ does not intersect the imaginary axis.
By \cite{Henry}, hyperbolic equilibria $W$ possess unique locally invariant local stable and unstable manifolds $W^s$ and $W^u$ in $X^\alpha$ which are graphs over their tangent spaces $E^s, E^u$ at W, by the spectral decomposition $X^\alpha=E^s\oplus E^u$.
They consist of all initial conditions $w_0$ such that, for all $\pm t\in \mathfrak{S}_\Theta$\,, the solution $w(t)=\Phi^t(w_0)$ exists, remains close to $W$, and converges to $W$, respectively, for $\pm\Re t \rightarrow +\infty$.
Under our differentiability assumptions, the local invariant manifolds are complex differentiable, and hence analytic.
See also \cite{ChowHale}.

\emph{Heteroclinic} solutions $\Gamma(r) \in X^\alpha$ in real time $r\in\mathbb{R}$ connect time-independent equilibria $w=W_\pm$\,, i.e.
\begin{equation}
\label{heteroclinic}
\Gamma(r) \rightarrow W_\pm \ \textrm{in}\ X^\alpha , \quad\textrm{for}\quad r\rightarrow\pm\infty.
\end{equation}
Depending on context, such $\Gamma$ are also called \emph{connecting orbits}, \emph{traveling fronts} (or \emph{backs}), or \emph{solitons}.
By definition, heteroclinic orbits between hyperbolic equilibria enter the local invariant manifolds $W_+^s$ and $W_-^u$ at $W_\pm$\,, respectively, in real time $r\rightarrow\pm\infty$.
We abbreviate heteroclinicity as $\Gamma:W_-\leadsto W_+$\,.
Unless specified otherwise, explicitly, we subsume the \emph{homoclinic} case $W_+=W_-$ of non-constant $\Gamma(r)$ under the heteroclinic label.

Complex ODEs have most expertly been presented in \cite{Ilyashenko}.
The above properties, in the setting of analytic semigroups,  closely mimic the ODE case.
However, we have to exercise caution when using real versus imaginary time.
The sectorial resolvent property defines local analytic semigroups along complex rays $0<r\mapsto r\exp(\mi\theta)$, in a sector $|\theta|\leq\Theta<\pi/2$.
Even for Laplacian generators $A$, however, the purely imaginary case $\theta=\pi/2$ only defines a strongly continuous semigroup -- as the Schrödinger case demonstrates \cite{Pazy}.
Instead, we use complex continuation $t\mapsto\Gamma(t)$ of real-time heteroclinic orbits, only, to venture into purely imaginary times $t=\mi s$.

Along these lines, we do not favor the language of complex foliations, which so elegantly blurs the distinction between real and imaginary time.
In fact, our PDE motivation prods us to insist on the real axis, as a distinguished direction of time.

\subsection{Abstract result} \label{AbsRes}

Our main abstract result, in theorem \ref{thmmain} below, indicates how real analytic heteroclinic solutions $\Gamma(t):W_-\leadsto W_+$, for real times $t=r\in\mathbb{R}$, cannot be complex entire, for all complex times $t=r+\mi s \in\mathbb{C}$.
This supports the initial paradigm of section \ref{Para}.

\emph{Spectral non-resonance} of the unstable part of the linearization at the source equilibrium $W_-$ will take center stage in our analysis.
In a finite-dimensional complex analytic ODE setting $X=\mathbb{C}^N$, suppose an equilibrium $W$ possesses diagonalizable linearization $L=A+f'(W)$ and spectrum $\mathrm{spec}\,L=\{\mu_0,\ldots,\mu_k\}$.
We call $W$ \emph{spectrally resonant}, if there exists a nonnegative integer multi-index $\mathbf{m}=(m_0,\ldots,m_k)$ of order $\lvert\mathbf{m}\rvert := m_0+\ldots+m_k\geq 2$, such that
\begin{equation}
\label{defresonance}
\mu_j=\mathbf{m}\cdot\mu:= m_0\mu_0+\ldots+m_k\mu_k
\end{equation}
holds, for some $\mu_j$\,.
Else, we call $W$ \emph{non-resonant}. 
For analytic semigroups we apply the same terminology to finite unstable spectrum.

Our main result excludes complex entire heteroclinic orbits $\Gamma(t):W_-\leadsto W_+$\,.

\begin{thm}\label{thmmain}
In the above complex analytic setting, assume 
\begin{enumerate}[(i)]
  \item $\Gamma: W_-\leadsto W_+$ is heteroclinic between hyperbolic equilibria $W_\pm$ in real time $t=r\in\mathbb{R}$;
  \item the target equilibrium $W_+$ is linearly asymptotically stable in real time $t=r\geq 0$;
  \item the local unstable manifold $W_-^u$ of $W_-$ is finite-dimensional; 
  \item the linearization at $W_-$ within $W_-^u$ is diagonalizable, with real non-resonant unstable spectrum.
\end{enumerate}
Then $\Gamma(t)$ cannot possess any complex entire time extension to all $t=r+\mi s\in\mathbb{C}$.
\end{thm}

To our knowledge, such a link has first been established in a PDE context in the dissertation thesis \cite{Stukediss,Stukearxiv}.
See section \ref{Sing} for further discussion.
%For a more detailed survey of the ODE case, we refer to our recent survey in \cite{fiedlerClaudia}, and to the discussion sections \ref{1000}, \ref{EntirePer} below.

The significance of non-resonance arises from \emph{Poincaré's theorem on analytic linearization}; see \cite{Ilyashenko} for a complete proof.
We postpone further comments to section \ref{Ete}, where we prove theorem \ref{thmmain}.

Unstable targets $W_+$ are admissible, if $\Gamma$ emanates from the one-dimensional fastest unstable manifold of $W_-$\,. 
See corollary \ref{cormain} in section \ref{Fastest} below.
A much less explicit foreboding, in the different and very restrictive setting of complex entire diffeomorphisms in $\mathbb{C}^2$, can be attributed to Ushiki \cite{Ushiki-2}.
See also \cite{Ushiki-N} for the special heteroclinic case $i(W_-)=1,\ i(W_+)=N-1$ in higher finite dimensions $\mathbb{C}^N$.
Remarkably early results by Rellich and Wittich excluded the existence of \emph{any} non-constant complex entire solutions for large classes of second order ODEs; see \cite{Rellich,Wittich1,Wittich2}.
Sections \ref{ODE2} and \ref{ODEd} will deepen our discussion of ODE aspects.

In the ODE case $X=\mathbb{C}^N$, heteroclinic orbits $\Gamma:W_-\leadsto W_+$ to unstable hyperbolic targets $W_+$ are admissible.
We still assume real non-resonant spectrum in the unstable or fast unstable manifolds at $W_-$\,.
At the target $W_+$\,, we also have to assume that the stable spectrum is real and, among itself, non-resonant.
See our recent study \cite{fiedlerClaudia} for details.
In the PDE case, infinite-dimensional non-resonance at $W^+$ cannot be brought to bear, because Poincaré linearization remains elusive in infinite dimensions.

\subsection{The quadratic heat equation} \label{Qua}

Our main specific example is the quadratic parabolic heat equation
\begin{equation}
\label{heat}
w_t=w_{xx}+6w^2-\lambda\,.
\end{equation}
Indices indicate partial derivatives.       
We consider $0<x<\tfrac{1}{2}$ under Neumann boundary conditions, and assume real parameter values $\lambda>0$.
Note how any real or complex solution $w=w(t,x)$, at parameter $\lambda$, extends periodically, by iterated reflection through the boundary, and therefore rescales to solutions
\begin{equation}
\label{nscale}
w_n=n^2w(n^2t,nx), \qquad \lambda_n=n^4 \lambda,
\end{equation}
for integer $n=1,2,3,\ldots$, at rescaled parameters $\lambda_n$\,.
                
\begin{figure}[t]
\centering \includegraphics[width=0.9\textwidth]{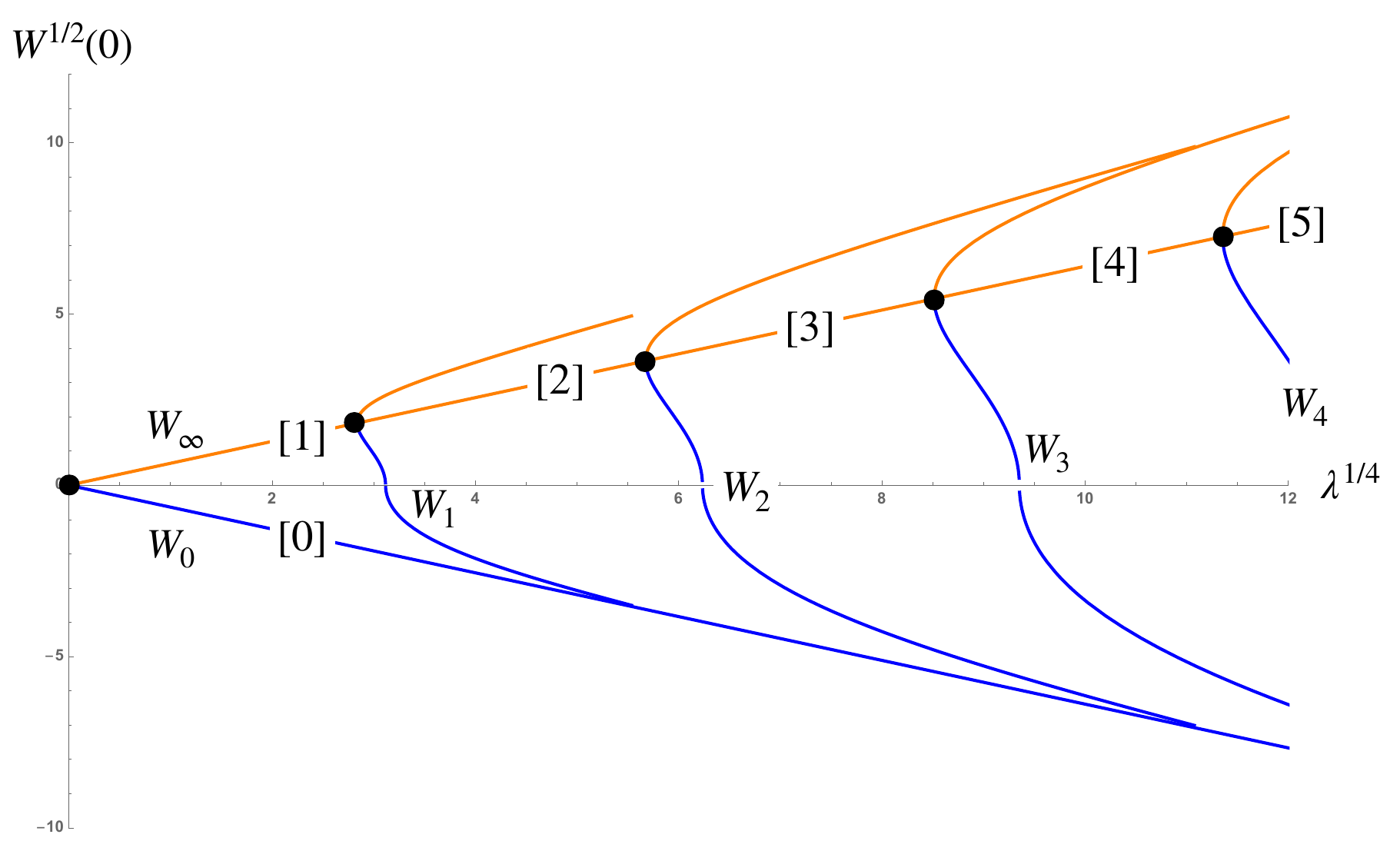}
\caption{\emph{
Bifurcation diagram of equilibria $W_n$ for the quadratic heat equation \eqref{heat}, \eqref{heat0} and parameters $\lambda>0$.
The horizontal axis $\lambda^{1/4}$ and the vertical axis of signed roots $W^{1/2}(0)$ render the equilibrium scaling \eqref{nscale} as stretching by a factor $n$.
Morse indices $i(W)$ are indicated in brackets $[\cdot]$, for $W=W_0\,, W_\infty$\,.
Note how $W_0$ (blue) is the unique asymptotically stable equilibrium.
Black dots mark bifurcations.
The saddle node bifurcation of the trivial spatially homogeneous branches $W_0\,, W_\infty$  is rendered as a wedge, at $\lambda=W=0$.
The remaining bifurcations are supercritical pitchforks of branches $W_n$ with Morse index $i=n$, emanating from $W_\infty$ (orange) at $\lambda=\tfrac{2}{3}(n\pi)^4$ with successively higher unstable dimensions; see \eqref{biflambda}.
The upper branches $W_n(0)=-n^2 e_2$ (orange) and the lower branches $W_n(0)=-n^2 e_3$ (blue) are related by a half-period shift in $x$.
The vertical tangents at $W_n(0)=0$  (blue) are an artifact of axis scaling, and do not indicate linear degeneracy.
}}
\label{fig1}
\end{figure}

Consider real solutions $w$ in real time $t=r$.
Scalar parabolic PDEs in one space dimension possess a gradient-like structure, in great generality.
See \cite{LappicyAttractors, LappicyBlowup}, most recently, and the many earlier references there. 
By the gradient-like structure, all non-equilibrium eternal solutions $\Gamma(r)$ of \eqref{heat}  are then heteroclinic among real equilibria $w=W_n(x)$.
The two trivial homogeneous equilibria are $W_\infty=+\sqrt{\lambda/6}$ and $W_0=-\sqrt{\lambda/6}$.
All nontrivial equilibria $W=W(x)$ are spatially non-constant solutions of the equilibrium ODE
\begin{equation}
\label{heat0}
0=W_{xx}+6W^2-\lambda\,.
\end{equation}
Therefore, they are real-valued \emph{Weierstrass elliptic functions} $-\wp(z)$, properly translated, and scaled according to \eqref{nscale}.
See \cite{Akhiezer, Fricke, Lang, Lamotke} for some background, including standard notation.
The functions $\wp(z)$ are meromorphic and doubly periodic in $z\in\mathbb{C}$.
Their period lattice is generated by $1, \tau$ with $\Im\tau>0$.
In our case, $W_n(x)=-n^2 \wp(nx+\tfrac{1}{2}\tau)$ or $-n^2 \wp(nx+\tfrac{1}{2}+\tfrac{1}{2}\tau)$, at $\lambda=\tfrac{1}{2}n^4 g_2(\tau)$, are real, of Morse indices (alias unstable dimensions) $i(W_n)=n$. 
See figure \ref{fig1} for a bifurcation diagram, and section \ref{wp} for details on this standard notation.

The analytic semigroup setting for \eqref{heat} is standard:  $A:=\partial_{xx}\,,\ X=L^2,\ \alpha=\tfrac{1}{2}$, and $X^{1/2}=H^1=W^{1,2}$ with the $H^1$-norm $\|w\|_{H^1}^2 := \|w\|_{L^2}^2 + \|\partial_xw\|_{L^2}^2$\,.
Any complex sector $t\in \mathfrak{S}_\Theta$ with opening angle $0<\Theta<\pi/2$ is admissible for the semigroup.
The analytic nonlinearity $f:H^1\rightarrow L^2$ is induced by the polynomial $f(w):= 6w^2-\lambda$.
As a technical consequence of this setting, theorem \ref{thmmain} addresses blow-up in $H^1$.
For $L^\infty$ blow-up, see our comments in section \ref{Sup}.

Note that the resolvent of $A$ is compact.
Hence unstable manifolds are automatically finite-dimensional.
By standard Sturm-Liouville theory \cite{Hartman}, the selfadjoint linearization $L=A+f'(W)$ at any real equilibrium $W$ possesses simple real spectrum
 \begin{equation}
\label{mu}
\mu_0>\mu_1>\ldots\searrow-\infty
\end{equation} 
with $L^2$-complete orthonormal eigenfunctions $\phi_0, \phi_1,\ldots$\,; see also \eqref{linmu}, \eqref{linbc} below.
Hyperbolic $W=W_n$ of Morse index $n=i(W)$ are characterized by $\mu_{n-1}>0>\mu_{n}$\,.

Real heteroclinic orbits have been studied in much detail and generality, for decades, starting with Chafee and Infante \cite{chin74}, and many others.
For some results and surveys on general dissipative nonlinearities $f(w)(x)=f(x,w,w_x)$ see \cite{raugel, firoSFB, firoFusco} and the many references there.
In addition to the gradient-like structure mentioned above, intersection and nodal properties like the zero number play an essential role here.
The results apply, in particular, to cubic dissipative regularizations $f^\eps(w):= 6w^2-\lambda-\eps w^3$ of \eqref{heat}.
For $\eps> 0$, dissipativeness implies that the global attractors $\mathcal{A}^\eps$ consist of precisely all eternal solutions.
By the gradient-like structure, the global attractors decompose into equilibria, and heteroclinic orbits among them.

The global attractors $\mathcal{A}^\eps$ of the dissipative regularizations are all $C^0$ orbit equivalent to the standard and celebrated Chafee-Infante attractors.
In the limit $\eps\searrow 0$, the largest positive homogeneous equilibrium escapes to $+\infty$.
The remaining set of eternal solutions, at $\eps=0$, is bounded.
\,In fact, for any two hyperbolic equilibria $W_\pm$ with Morse indices $i(W_\pm)$, we encounter heteroclinic orbits
\begin{equation}
\label{heathet}
\Gamma:W_-\leadsto W_+ \quad\Leftrightarrow\quad i(W_-)>i(W_+)\,.
\end{equation}
This is quite analogous to the dissipative Chafee-Infante case, except for those orbits in the unstable manifolds $W_-^u$ which blow up in finite real time.
For a result on prescribed sign-changing profiles at real-time blow-up, which involves nodal properties, see \cite{FiedlerMatano}.
We omit further details, because these constructions and results are not a main aspect of our present paper.

Instead, we focus on blow-up in complex time.
Based on theorem \ref{thmmain}, we obtain the following result.

\begin{thm}\label{thmheat}
In the bifurcation diagram of figure \ref{fig1}, consider heteroclinic orbits \eqref{heathet} emanating from hyperbolic source equilibria $W_-$ along branches $W_n$\,, or along $W_\infty$\,. 
Then unstable spectral non-resonance holds at most parameters $\lambda$.\\
More precisely, assume unstable dimensions $1\leq n\leq 22$.
Then, at most, unstable spectral resonance \eqref{defresonance} along $W_n$ occurs for a discrete set of resonant parameters $\lambda\geq 0$. \\
For spatially homogeneous $W_-=W_\infty$\,, no restriction is imposed on the unstable dimension.
However, the discrete set of spectrally resonant parameters does in fact accumulate, from above, at the pitchfork bifurcation parameters 
\begin{equation}
\label{biflambda}
\lambda_{n0}:= \tfrac{2}{3}(n\pi)^4\,,\quad n=1,2,3,\ldots \,.
\end{equation}
In either case, $\Gamma:W_-\leadsto W_0$ cannot possess any entire extension to all complex times $t=r+\mi s\in\mathbb{C}$, at non-resonant parameters $\lambda$.
\end{thm}

\subsection{Fast unstable manifolds}\label{Fast}

A variant of theorem \ref{thmheat} addresses local \emph{fast unstable manifolds} $W_-^{uu}$ at any equilibrium $W_-=W_n$\,.
These locally invariant submanifolds of dimension $d<n$ are uniquely characterized by asymptotic decay to equilibrium at exponential  rates faster than $\mu_d>0$, in real time $t=r\searrow -\infty$.
They are therefore tangent to the span of the associated eigenfunctions $\phi_0,\ldots,\phi_{d-1}$\,.
Like the local unstable manifold $W_-^u$ itself, the fast unstable manifolds $W_-^{uu}$ are complex analytic in a neighborhood of $W_-$\,.

\begin{thm}\label{thmheatd}
In the setting of theorem \ref{thmheat} and figure \ref{fig1}, assume the heteroclinic orbit $\Gamma: W_-\leadsto W_0$ emanates from $W_-$ via the local fast unstable manifold $W_-^{uu}$ of dimension $d<i(W_-)$.
Assume
\begin{equation}
\label{dn}
d<1+\tfrac{1}{\sqrt{2}}\,i(W_-)\,.
\end{equation}
Then unstable spectral non-resonance within $W_-^{uu}$ holds at most parameters $\lambda$\,. 
I.e., unstable spectral resonance \eqref{defresonance} only occurs at a discrete set of resonant parameters $\lambda\geq 0$.
This holds for all\, $W_-=W_n$ and includes the case $W_-=W_\infty$\,, without any exceptional accumulations at the bifurcation values \eqref{biflambda}.\\
In particular, $\Gamma:W_-\leadsto W_0$ cannot possess any complex entire extension to all complex times $t=r+\mi s\in\mathbb{C}$, at non-resonant parameters $\lambda$.
\end{thm}

At bifurcation, the non-resonance bound \eqref{dn} is optimal, asymptotically for large $i(W_-)$; see \eqref{pythagoras} in the proof section \ref{Wuu} below.

\subsection{Heat versus Schrödinger equation} \label{Schro}

Given heteroclinic orbits $\Gamma(t)$ as in the theorems above, which are real analytic but not complex entire, let $\psi(s):=\Gamma(r_0-\mathrm{i}s)$ track complex time extensions $\Gamma(t)$ in the imaginary time direction, for fixed real $r_0$.
Then $\psi(s)\in H^1$ solves the nonlinear and nonconservative Schrödinger equation
\begin{equation}
\label{psi}
\mi\psi_s=\psi_{xx}+6\psi^2-\lambda\,.  
\end{equation}
The solutions of \eqref{psi} define a strongly continuous (but not analytic) local semigroup $\psi(s)=S^s\psi_0$ in $H^1$ \cite{Pazy, Tanabe}.
We can still define \emph{blow-up} as
\begin{equation}
\label{blow-up}
\|\psi(s)\|_{H^1}\rightarrow\infty \quad\textrm{for}\quad s\nearrow s^*\ \textrm{and}\ s\searrow -s^*\,,
\end{equation}
at some finite positive $s^*$.

Note \emph{time reversibility} of \eqref{psi} here, under complex conjugation $\psi\mapsto\overline{\psi}$.
Indeed, $\psi(s)$ is a solution, if and only if $\overline{\psi}(-s)$ is.
In other words, the strongly continuous local solution semigroup $S^s$ satisfies
\begin{equation}
\label{reversibility}
S^s(\overline{\psi_0})=\overline{S^{-s}(\psi_0)}.
\end{equation}
In particular, we have a strongly continuous local solution \emph{group} $S^s$, for the Schrödinger equation \eqref{psi}.

\begin{thm}\label{thmpsi}
Under the assumptions of theorems \ref{thmheat} or \ref{thmheatd}, there exists a real initial condition $\psi_0:=\Gamma(r_0)$ for the Schrödinger equation \eqref{psi}, such that the solution $\psi(s)$ blows up at some finite real time $s^*>0$ and, correspondingly, blows down from $-s^*<0$.
\end{thm}

\subsection{The fastest unstable manifold} \label{Fastest}

The case of the  \emph{fastest unstable manifold} $W_-^{uu}$ of dimension $d=1$  is particularly intriguing. Under Dirichlet boundary conditions, it has already been addressed in \cite{Stukediss,Stukearxiv}. 
We comment on the spatially homogeneous Neumann case $W_-=W_\infty$ in the next section \ref{ODE2}.

For one-dimensional $W_-^{uu}$, we can drop any stability assumption on the target $W_+$\,, in the semigroup setting of theorem \ref{thmmain}. 
By \cite{brfi86}, a heteroclinic orbit $\Gamma(r)$ emanates from the fastest unstable manifold $W_-^{uu}$ of dimension $d=1$, at $W_-$\,, if and only if $\Gamma(r)$ becomes tangent to the eigenvector $\pm\phi_0$ of the (unique) fastest unstable eigenvalue $\mu_0$\,, for $r\searrow-\infty$.

\begin{cor}\label{cormain}
In the analytic semigroup setting of section \ref{Sgr}, assume 
\begin{enumerate}[(i)]
  \item $\Gamma: W_-\leadsto W_+$ is heteroclinic between hyperbolic equilibria $W_\pm$ in real time $t=r\in\mathbb{R}$;
  \item the unstable eigenvalue $\mu_0$ of largest real part, at  $W_-$\,,  is unique, real, and algebraically simple; 
  \item $\Gamma(r)$ is tangent to $\pm\phi_0$\,, at $W_-$\,.
\end{enumerate}

Then $\Gamma(t)$ cannot possess any complex entire time extension to all $t=r+\mi s\in\mathbb{C}$.

In particular, theorems \ref{thmheatd}, \ref{thmpsi} remain valid under these modified assumptions.
\end{cor}

\begin{proof}
We will eventually suppose $t\mapsto\Gamma(t)$ is complex entire, and obtain a contradiction.
For the moment, we proceed without this indirect extra assumption.

The simple eigenvalue $\mu_0$ alone cannot be resonant; see \eqref{defresonance}.
Local Poincaré linearization in $W_-^{uu}$ at $W_-$ therefore implies periodicity
\begin{equation}
\label{Gammaperiod}
\Gamma(t+\mi p)=\Gamma(t)
\end{equation}
for all $t=r+\mi s\in\mathbb{C}$ with, say, $r\leq 0$.
For details see \eqref{P} below.
Here $p=2\pi/\mu_0>0$ is the minimal period of $\Gamma$ in the imaginary time direction.
At fixed $s$, let $-\infty<r<r^*(s)\in(0, +\infty]$ denote the maximal interval of existence for $\Gamma(r+\mi s)$, in real time $r$.
Define the analyticity domain
\begin{equation}
\label{D}
\mathcal{D}:=\{r+\mi s\,|\,-\infty<r<r^*(s)\}\,
\end{equation}
by the semigroup solution of \eqref{sgr} in real time.

We now invoke our indirect assumption that $\Gamma$ is complex entire, i.e. $\mathcal{D}=\mathbb{C}$. We claim
\begin{equation}
\label{rper}
\Gamma(r+\mi s)\rightarrow W_\pm \quad \textrm{for}\quad r\rightarrow \pm\infty\,,
\end{equation}
and uniformly for $0\leq s\leq p$.
For $r\searrow -\infty$, convergence follows from linearization \eqref{P}.
For $r\nearrow +\infty$, the locally $s$-uniform limit $\Gamma(r+\mi s)\rightarrow W_+=W_0$ follows from strong continuity and commutativity \eqref{flow} in sectors.

In particular, \eqref{rper} implies that $\Gamma$ is uniformly bounded, hence constant by Liouville's theorem.
This contradiction proves the corollary, in the semigroup setting.
\end{proof}

For the quadratic heat and Schrödinger equations, \eqref{heat} and \eqref{psi}, we claim that the target $W_+$ necessarily coincides with the only asymptotically stable equilibrium $W_0$\,.
Therefore theorems \ref{thmheatd} and \ref{thmpsi} apply directly, as stated. 

We fill in some details.
Originating from spatially non-homogeneous $W_-=W_n$\,, the fastest real heteroclinic orbits $\Gamma$ are spatially non-homogeneous..
In fact, they emanate from $W_-$\,, tangentially along the eigenfunction $-\phi_0<0$ of the largest Sturm-Liouville eigenvalue $\mu_0>0$.
By \cite{brfi86}, the real heteroclinic orbit $\Gamma:W_n\leadsto W_+$ must decrease monotonically to its target equilibrium $W_+<W_-$\,.
In fact, $\Gamma_r<0$ for all real times $r$.
Monotone convergence to $W_+$ implies convergence tangent to $\phi_{+0}$\,, the positive eigenfunction of the largest eigenvalue $\mu_{+0}$ at $W_+$\,.
Therefore $\mu_{+0}< 0$, and the hyperbolic target $W_+$ is an asymptotically stable equilibrium.
This observation applies to general nonlinearities $f(w)(x)=f(x,w(x),w_x(x))$ in \eqref{sgr}, \eqref{heat}. 
In our quadratic case \eqref{heat}, it identifies the target $W_+=W_0$ as the lower homogeneous equilibrium.

For $f(w)(x)=f(w(x))$, we can alternatively invoke phase plane analysis of $(W,W_x)\in\mathbb{R}^2$, as in \cite{brfi88,brfi89, FiedlerRochaWolfrumHam, FiedlerRochaHam}.
Spatially non-homogeneous real Neumann equilibria $W_\pm$ of \eqref{sgr} are then parts of non-constant, $x$-periodic orbits of the Hamiltonian pendulum $W_{xx}+f(W)=0$.
Such orbits are excluded as target options $W_+$\,, because they would have to be nested with $W_-$\,, in the real phase plane, rather than located strictly below $W_-$ in $\mathbb{R}$.
Therefore, the target $W_+$ coincides with the largest spatially homogeneous equilibrium $f(W)=0$ below $W_-$\,.
In our quadratic case, this equilibrium is $W_+=W_0$ again.

On the fastest unstable manifold $W_-^{uu}$, we can also identify some blow-up, more specifically.
Together with reversibility \eqref{reversibility}, periodicity \eqref{Gammaperiod} implies that $\Gamma(r+\mi s)$ is real (for all $x$), if and only if \begin{equation}
\label{Gammareal}
s\equiv 0  \quad \textrm{or}\quad s\equiv p/2 \mod p\,.
\end{equation}
 
The real orbit $w(r)=\Gamma(r+\mi p/2)$ emanates from $W_-=W_n$\,, tangentially to the eigenfunction $+\phi_0>0$. 
It is therefore monotonically increasing with $r$, rather than decreasing.
Let $r\in (-\infty,r^*(s))$ denote the maximal interval of existence of $r\mapsto \Gamma(r+\mi s)$.
By absence of any equilibrium $W>W_n$\,, we have
\begin{equation}
\label{Gammablowup}
\|w(r)\|_{H^1} \rightarrow\infty \quad\textrm{for}\quad r\nearrow r^*(p/2)\leq\infty\,.
\end{equation}
For \eqref{heat}, the real comparison principle can then be invoked to show finite time blow-up $r^*(p/2)<\infty$, in the sense of \eqref{defblow-up}, rather than grow-up $r^*(p/2)=\infty$ as in \cite{Juliana}.

Next consider Schrödinger blow-up in \eqref{psi} for solutions $\psi(s)=\Gamma(r_0+\mi s)$ starting from some suitable $\psi_0=\Gamma(r_0)$ on that fastest real heteroclinic orbit $\Gamma:W_n\leadsto\Gamma_0$\,; see theorem \ref{thmpsi}.
By periodicity \eqref{Gammaperiod}, reversibility \eqref{reversibility}, and because $\Gamma(\mathbb{R})$ is real, the boundary $s\mapsto r^*(s)\leq+\infty$ of the analyticity domain $\mathcal{D}$ in \eqref{D} is $p$-periodic and symmetric:
\begin{equation}
\label{r*symmetry}
r^*(p+s)=r^*(s)=r^*(p-s),
\end{equation}
for all $s\in\mathbb{R}$.
Since $\mathcal{D}\subsetneq\mathbb{C}$, the  lower semicontinuous boundary $r^*(s)\not\equiv +\infty$ attains its finite minimum $r_0$ at some minimal $s_0>0$.
Note $0<s_0\leq p/2$, by symmetry \eqref{r*symmetry} and minimality of $s_0$\,.
Moreover $r^*(s_0)=r_0$\,, and $s^*(r_0)=s_0\leq p/2$ is the positive Schrödinger blow-up time of $\psi(s)$, starting from $\psi_0=\Gamma(r_0)$.
For details on this last step, see our proof of theorem \ref{thmpsi} in section \ref{HeaSchro} and, in particular, the discussion of the rectangle $\mathcal{R}$ with upper right corner $r_0+\mi s_0$ in \eqref{rectangle}.

As a caveat, we add that we have \emph{not} excluded Schrödinger blow-up times $s^*(r)\in(p/2,+\infty)$, for other $\psi_0=\Gamma(r)$ and $r>r_0$.
However, reversibility then implies violation $\psi(p/2)\neq\psi(-p/2)$ of $p$-periodicity, at nonvanishing imaginary part.
This would be reminiscent of Masuda's discrepancy \cite{Masuda1, Masuda2} between upper and lower complex detours around real-time blow-up.
See section \ref{Sing} for further discussion.

As a peculiarity of dimension $d=1$, and an advantage of our complex viewpoint, we now understand how fastest real heteroclinicity $\Gamma: W_n\leadsto W_0$, at fixed imaginary level $s=0$, and real blow-up, at fixed imaginary level $s_0\leq p/2=\pi/\mu_0$\,, are just two aspects of one and the same underlying complex fastest unstable manifold $W_-^{uu}$.
Indeed, the two phenomena are contingently related, half a period apart, in the foliation by fixed imaginary times $\mi s$, for $s\,\mathrm{mod}\, p$. 
In $s$, the bounded fastest heteroclinic orbit $\Gamma$ also provides certain real initial conditions $\psi_0=\Gamma(r_0)$ with complex Schrödinger blow-up of $\psi(s)$ at times $\pm s^*$ not exceeding the half-period $p/2$. We will return to a more global viewpoint of the fastest unstable manifold, and continuation beyond the analyticity domain $\mathcal{D}$, in discussion sections \ref{ODEd} on spatially homogenous ODE solutions, and \ref{Boundary} on the closure and boundary of unstable manifolds.
See also figures \ref{fig2}, \ref{fig3}, and the next section.

\subsection{The quadratic ODE} \label{ODE2}

\begin{figure}[t]
\centering \includegraphics[width=0.86\textwidth]{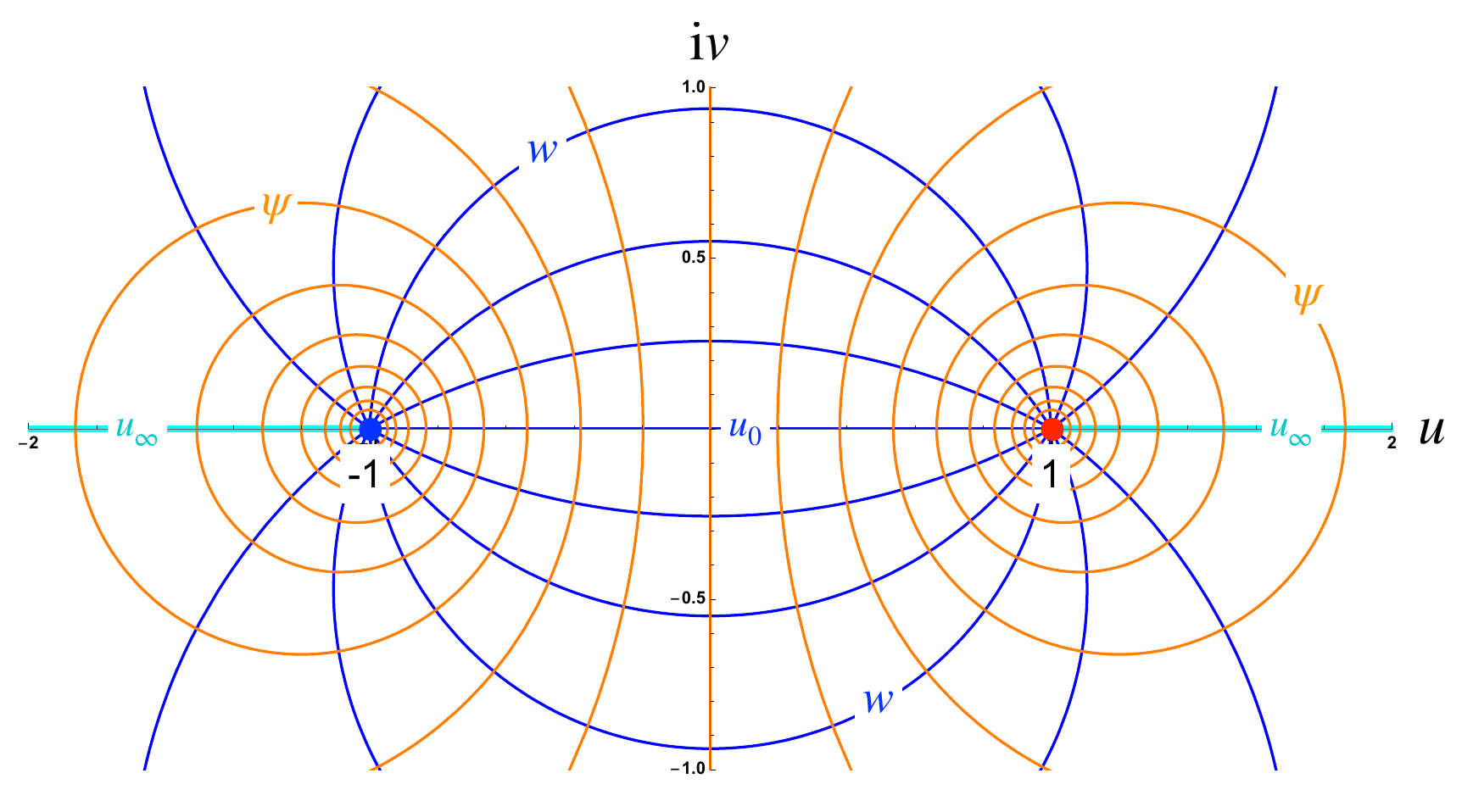}
\caption{\emph{
Phase portrait of complex ODEs \eqref{ODEwr} -- \eqref{ODEuv} for the spatially homogeneous heteroclinic orbit $\Gamma:W_-\leadsto W_+$\,.
Orbits $r\mapsto\Gamma(r+\mi s)=u+\mi v$ of \eqref{ODEwr}, \eqref{ODEuv}, in real time $r$, are circle arcs (blue). 
They are heteroclinic from source  $W_-=W_\infty=+1$ (red dot) to sink $W_+=W_0=-1$ (blue dot).
The $\pi$-periodic orbits $s\mapsto\psi(s)=\Gamma(r-\mi s)=u+\mi v$ of \eqref{ODEpsi}, in contrast, are full circles (orange), each surrounding one of the two equilibria $W_\pm=\mp 1$.
Because the flows in real and imaginary time commute, the orange circles also serve as isochrones which globally synchronize the blue heteroclinics.
Conversely, the blue circle arcs are isochrones which globally synchronize the orange periodic orbits.
For any fixed time $t$, the holomorphic flow map $w_0\mapsto \Phi^t(w_0)$ is conformal, i.e. angle-preserving. 
In particular, the blue and orange circle families are mutually orthogonal.\\
The imaginary $v$-axis (orange) blows up in finite time $s^*=\pm\pi/2$; see \eqref{s*}.
On the Riemann sphere $\widehat{\mathbb{C}}$, this is just a longitude circle of period $\pi$ through the North and South Poles at $w=\infty$ and $w=0$.
The real $u$-axis features the real heteroclinic orbit $u_0: +1 \leadsto -1$, in blue, as well as the two cyan parts of the blow-up heteroclinic $u_\infty:+1\leadsto\infty$ and the blow-down heteroclinic $u_\infty:\infty\leadsto -1$.
On the Riemann sphere $\widehat{\mathbb{C}}$, these three line segments combine to another longitude circle, perpendicular to the first at the intersections 0 and $\infty$.
In conclusion, the blue bounded heteroclinic orbit $u_0$ on the real $u$-axis, in real time, gives rise to finite time blow-up and blow-down on the orange imaginary $v$-axis, in imaginary time, when started at $u=v=0$.
The complex viewpoint also reveals how the blue real heteroclinicity $u_0$ and the cyan real blow-up/blow-down segments $u_\infty$\,, in real time, become three aspects of one and the same underlying trajectory $\Gamma$, in complex time $t=r+\mi s$.\\
}}
\label{fig2}
\end{figure}

Even in the most elementary and explicit ODE case $\dot{w}=6w^2-\lambda$
of spatially homogeneous solutions, it is instructive to see the above results at work; see figure \ref{fig2}.
In complex notation and real time $t=r$, the associated spatially homogeneous equations \eqref{heat} and \eqref{psi} then read 
\begin{align}
    \label{ODEwr}
    w_r &=w^2-1, \\
    \label{ODEpsi}
    \mi\psi_s &= \psi^2-1\,.
\end{align}
Here we have first rescaled $w$ by a factor $\sqrt{\lambda/6}$\,, and then rescaled real time by a factor $\lambda>0$.
%(Complex quadratic vector fields can be rescaled, similarly,  a complex rescaling of time and the purely quadratic case $\lambda=0$.)
For $w=u+\mi v$ in \eqref{ODEwr} we obtain the equivalent real system
\begin{equation}
\label{ODEuv}
\begin{aligned}
    \dot{u} &= u^2-v^2-1\,,  \\
    \dot{v} &= \quad 2uv\,.
\end{aligned}
\end{equation}
The equilibria are the attractor $W_0=-1$ (blue dot) and the repellor $W_\infty=+1$ (red dot).
All nonstationary orbits (blue) are heteroclinic $\Gamma: W_\infty\leadsto W_0$.
The real $u$-axis $\{v=0\}$ is invariant and contains the monotonically decreasing heteroclinic orbit $u_0:1\leadsto -1$.
The two unbounded real orbits $u_\infty$ on the $u$-axis (cyan) look like an exception, at first: they blow up or blow down at $u=\pm\infty$ in finite positive or negative time.

Unlike the general PDE case of $W_-=W_n$\,, the ODE \eqref{ODEwr} can be regularized on the Riemann sphere $\widehat{\mathbb{C}}:=\mathbb{C}\cup\{\infty\}$\,.
Just note equivariance of $\dot{w}=w^2-1$ under the involution $w\mapsto 1/w$.
In other words, $w(t)$ is a solution, if and only if $1/w(t)$ is.
In particular, the two blow-up pieces of $u_\infty$ can be joined to indicate yet another heteroclinic orbit $u_\infty : 1\leadsto -1$ on the Riemann sphere.
Indeed, cyan $u_\infty$ is just the involutive copy of $u_0$\,, on $\widehat{\mathbb{C}}$\,.

Because $\dim _\mathbb{C}w = 1$, Poincaré linearization near $W_\infty=1$ shows that the complex heteroclinic orbit $\Gamma$ is periodic of minimal period $p=\pi$ in imaginary time.
In other words, the meromorphic heteroclinic orbit $t\mapsto\Gamma(t)$ provides a bi-holomorphic map between complex cylinders,
\begin{equation}
\label{ODEhol}
\Gamma: \mathbb{C}/\pi\mi\mathbb{Z} \ \rightarrow\  \widehat{\mathbb{C}}\setminus\{\pm 1\}\,.
\end{equation}

By corollary \ref{thmpsi}, blow-up to $\psi(s)=\infty$ must occur, for some $r=r_0$ and at some $s=s^*$.
Bi-holomorphy \eqref{ODEhol} implies that $r_0$ and $s^*$ are unique.
Here real $\Gamma(r):=u_0(r)\in (-1,1)$ tracks the bounded decreasing real heteroclinic orbit. 
By time reversibility \eqref{reversibility} of \eqref{psi} in $s$, blow-up must therefore occur at 
\begin{equation}
\label{s*}
s^*\equiv\pi/2 \mod\pi,
\end{equation}
for some $\psi_0=\Gamma(r_0)$.
Also, $\Gamma(r\pm\pi\mi/2)\in\mathbb{R}\cup\{\infty\}$.

Without loss of generality, let us fix $\Gamma(0)=0$.
Explicitly, $v=0$ in \eqref{ODEuv} then identifies the biholomorphic cylinder map \eqref{ODEhol} as
\begin{equation}
\label{Gammaexplicit}
\Gamma(t)=-\tanh t,
\end{equation}
first for real $t=r\in\mathbb{R}$, i.e. for $u_0:1\leadsto -1$, and then for all $t\in\mathbb{C}$, by analytic continuation.
This identifies $r_0=0$ and the blow-up solution
\begin{equation}
\label{psiexplicit}
\psi(s)=-\mi\tan s\,.
\end{equation}
Complex time shifts of $\Gamma$, in \eqref{ODEhol}, induce the fractional linear, biholomorphic Möbius transformations of the Riemann sphere $\widehat{\mathbb{C}}$ which fix $\pm 1$. 
In \eqref{Gammaexplicit}, the flow property \eqref{flow} then amounts to the elementary addition theorems for the hyperbolic tangent.

In passing, we note that the real and imaginary parts of the inverse function $r+\mi s=-\mathrm{arctanh}\,w$ of the solution \eqref{Gammaexplicit} are first integrals of \eqref{ODEpsi} and \eqref{ODEwr}, respectively.
In figure \ref{fig2} their level curves are colored orange and blue/cyan.
Specifically, the blue/cyan heteroclinic orbits of \eqref{ODEuv}, \eqref{ODEwr}, in real time $\mathbb{R}\ni r\mapsto t=r+\mi s$ and for fixed $s$, foliate the cylinder $\widehat{\mathbb{C}}\setminus\{\pm 1\}$ into arcs of invariant Euclidean circles through $w=\pm1$, with centers on the imaginary $v$-axis.
The orange $\pi$-periodic orbits of \eqref{ODEuv}, in the Schrödinger variant \eqref{ODEpsi} of imaginary time $t=\mi s$, define a perpendicular foliation into circles with centers $w=a\in\mathbb{R},\ 1<|a|\leq\infty$ on the real $u$-axis.
In other words, their radii $0<\sqrt{a^2-1}\leq\infty$ are defined by the tangents from $a\in\mathbb{R}$ to the complex unit circle.
The straight axes themselves are "circles with centers at infinity", as already Cusanus has conceptualized \cite{Cusanus}.

In electrostatics, the blue/cyan and orange circle families of figure \ref{fig2} are known as field lines and potential levels of an electrical dipole with charges $\pm1$ located at $W=\pm1$.
Similarly they illustrate flow lines and potential of a static, incompressible, irrotational planar fluid flow with source and sink at $W=\pm1$.
See for example \cite{Marsden}, ch 5.3.

Only in the spatially homogeneous cases \eqref{ODEwr}, \eqref{ODEpsi}, of course, we can also trivialize the cylinder map \eqref{ODEhol} by the fractional linear, biholomorphic Möbius transformation $w\mapsto z=(w+1)/(w-1)$ of the Riemann sphere $\widehat{\mathbb{C}}$\,.  
This maps the equilibria $W=\pm 1$ to $z=\infty$ and $z=0$, respectively. 
It also globally linearizes the quadratic complex ODE $\dot{w}=w^2-1$, alias \eqref{ODEuv}, to become $\dot{z}=-2z$.
In $z$, the heteroclinic orbits of \eqref{ODEwr} are inward radial, and the periodic orbits of \eqref{ODEpsi} are circles of constant period $\pi$ and radius $|z|$.
Circles and lines in figure \ref{fig2} illustrate how Möbius transformations, just as circle inversions, map collections of circles and lines to circles and lines, preserving angles of intersection.
For polynomial scalar ODEs $\dot w=f(w)$ we continue this discussion in section \ref{ODEd}.

\subsection{Outline} \label{Out}

We proceed as follows.
In section \ref{Ete} we prove theorem \ref{thmmain}.
Section \ref{wp} adapts the classical approach via Weierstrass elliptic functions to derive the bifurcation diagram \ref{fig1} of all real equilibria.
In particular, lemma \ref{Wn} summarizes known Fourier expansions which hold globally along the bifurcating branches $W_n$\,.
In section \ref{NonRes}, explicit eigenvalue expansions to second order, at the pitchfork bifurcations, then allow us to establish unstable spectral non-resonance for all but a discrete set of parameters $\lambda$.
This proves theorem \ref{thmheat}, for the quadratic heat equation \eqref{heat}.
The proof of theorem \ref{thmheatd} on fast unstable manifolds is deferred to section \ref{Wuu}.
Blow-up in the quadratic Schrödinger equation \eqref{psi}, collateral to heteroclinicity in the heat equation \eqref{heat}, is established in section \ref{HeaSchro}.
We conclude with a rather detailed discussion.
See section \ref{Over} for a summary.
Readers interested in just some pertinent PDE literature may rush to section \ref{Sing}, directly and at their own risk.

\subsection{Acknowledgment}\label{Ack} 
This paper is dedicated to the memory of \emph{Marek Fila}, a source of inspiration and generous friendship for decades. 
He was the first to draw our attention to the topic of PDE blow-up, always most encouragingly and never losing faith in us, until we belatedly ventured to present this paper.
We will miss his cheerful comments forever.

We are also indebted to \emph{Vassili Gelfreich}, for very patient explanations of his profound work on exponential splitting asymptotics, to \emph{Anatoly Neishtadt} for lucid conversations on adiabatic elimination, to \emph{Jonathan Jaquette} and \emph{Jean-Philippe Lessard} for generously sharing and discussing numerical details of their work, and to \emph{Karsten Matthies, Carlos Rocha, Jürgen Scheurle} and late \emph{Claudia Wulff}, for their lasting interest.

\section{Real eternal contradicts complex entire} \label{Ete}

In this section we prove our main abstract result, theorem \ref{thmmain}.
We work in the setting of analytic semigroups $\Phi^t$ in sectors $t\in \mathfrak{S}_\Theta$ on a complex Banach space $X$, with $Y:=X^\alpha,\ 0\leq\alpha<1$, and norm $\|\cdot\|$ on $Y$; see section \ref{Sgr}.
Our proof is indirect: we will assume that the heteroclinic orbit $\Gamma:W_-\leadsto W_+$ between hyperbolic equilibria $W_\pm$\,, in real time $r$, extends globally to a complex entire solution $\Gamma(t)$ of \eqref{sgr}, for all $t=r+\mi s\in\mathbb{C}$.
%See assumption \emph{(i)}.

Theorem \ref{thmmain}, as well as its proof, is quite asymmetric with respect to the two equilibria $W_\pm$\,.
One difficulty arises in case parts of the heteroclinic orbits $r \mapsto \Gamma(r+\mi s)$, in real time $r$ and for fixed $s\neq 0$, traverse local neighborhoods of $W_+$ multiple times, before settling into the limiting equilibrium $W_+$\,.
For large $|s|$, potentially, this could happen in smaller and smaller neighborhoods.
The local asymptotic stability assumption \emph{(ii)} at $W_+$ eliminates this difficulty.
%At $W_-$\,, of course, i.e.~in reverse time $r\rightarrow -\infty$, reverse ``asymptotic stability'' is not a semigroup option.

A second difficulty arises for estimates in the imaginary time direction $\mi s$.
Note the absence, in general, of a strongly continuous local semigroup
\begin{equation}
\label{imflow}
s\mapsto \Phi^{\mi s}
\end{equation}
on all of $Y$, in imaginary time $t=\mi s$.
In particular, nonexistent $\Phi^{\mi s}$ cannot commute with sectorial $\Phi^t,\,t\in \mathfrak{S}_\Theta$\,, in the sense of \eqref{flow}.
Restricted to the (supposedly) complex entire heteroclinic orbit $\Gamma(r+\mi s)$, however, the semigroup acts as a flow on $\Gamma$, globally in forward, backward, and sideways time directions $t\in\mathbb{C}$, by time shift
\begin{equation}
\label{shift}
\Phi^{t}\Gamma(t_0):= \Gamma(t+t_0).
\end{equation}
Commutativity, continuity, and even analyticity of $\Phi^t\Gamma(t_0)$ now hold trivially, for all $t,t_0\in\mathbb{C}$, under this restriction to $\Gamma$ itself.
So far, however, we are still lacking boundedness estimates and recurrence properties, in imaginary time.

Spectral non-resonance assumption \emph{(iv)} will come to our rescue, for real times $r\searrow-\infty$. 
The significance lies in \emph{Poincaré's theorem on analytic linearization}; see \cite{Ilyashenko} for a complete proof.
In absence of spectral resonances, the theorem provides a complex analytic, local diffeomorphism $\mathfrak{P}: (E_-^u,0) \rightarrow (W_-^u,W_-)$ which linearizes the ODE flow within the local unstable manifold $W_-^u$, near the equilibrium $W_-$\,. 
Note the identical embedding $\mathfrak{P}'(0): E_-^u\hookrightarrow Y$ of the finite-dimensional unstable eigenspace $E_-^u=T_{W_-}W_-^u$\,, in our case.

Restricted to the local unstable manifold $W_-^u$\,, all the way, let us now consider the linearization $L=A+f'(W_-)$, restricted to the unstable eigenspace $E_-^u$\,. 
By assumption \emph{(iv)} of theorem \ref{thmmain}, here $L$ is diagonalizable, with 
real and totally unstable spectrum, i.e. $\mathrm{spec}\,L >0$.
Imaginary time $t=\mi s$ rotates that real spectrum to become purely imaginary.
Poincaré linearization therefore keeps oscillations small, i.e. within the region of validity of the linearization, uniformly for all imaginary parts $s\in\mathbb{R}$.
Specifically, in the spectral notation of \eqref{defresonance}, we observe
\begin{equation}
\label{P}
\Gamma(t) = \mathfrak{P}\Big(\,\sum_{0\leq j< i(W_-)}\,a_j\,\exp(\mu_j t) \,\phi_j\, \Big)
\end{equation}
for all $r=\Re t \leq0$ and suitable coefficients $a_j\in\mathbb{C}$, provided $\|\Gamma(0)-W_-\|<\eps$ is chosen small enough.
Indeed, the coefficients $a_j\exp(\mu_j t)$ of any eigenvector $\phi_j$ remain unchanged, in absolute value, under the analytically equivalent linearized flow in imaginary time $t=\mi s$.

Moreover, any solution $w(\mi s)$ starting at $w_0$ near $W_-$ is quasi-periodic in imaginary time $s$.
See \cite{Bohr, Besicovitch, Corduneanu,Haraux, fiedlerClaudia} for details on almost- and quasi-periodicity.
In particular, $w_0$ is \emph{Poisson stable} in, both, forward and backward imaginary time $t=\mi s$, i.e.
\begin{equation}
\label{poisson}
w_0\in\boldsymbol{\alpha}_{\mi\mathbb{R}}(w_0)\cap\boldsymbol{\omega}_{\mi\mathbb{R}}(w_0).
\end{equation}
Here the sets of accumulation points of global solutions $w(\pm t_n)$ are called the $\boldsymbol{\omega}$\emph{-} and the $\boldsymbol{\alpha}$\emph{-limit set} $\boldsymbol{\omega}_\mathbb{R}(w_0)$ and $\boldsymbol{\alpha}_\mathbb{R}(w_0)$, respectively, for real times $t_n=r_n\nearrow+\infty$.
Real-time heteroclinic orbits $\Gamma:W_-\leadsto W_+$\,, for example, are characterized by $\boldsymbol{\alpha}_{\mathbb{R}}(\Gamma_0)=\{W_-\}$ and $\boldsymbol{\omega}_{\mathbb{R}}(\Gamma_0)=\{W_+\}$\,.
For imaginary times $t_n=\mi s_n$ and $s_n\nearrow+\infty$, we denote the corresponding limit sets by $\boldsymbol{\omega}_{\mi\mathbb{R}}(w_0)$ and $\boldsymbol{\alpha}_{\mi\mathbb{R}}(w_0)$.

Non-resonant Poincaré linearization, imaginary quasi-periodicity, and imaginary Poisson stability apply, in particular, within the complex analytic local unstable manifold $w_0:=\Gamma_0\in W_-^u$ of the hyperbolic heteroclinic source $W_-$\,.

\textbf{Proof of theorem \ref{thmmain}.}\\
\nopagebreak
We have assumed, indirectly, that real-time $\Gamma:W_-\leadsto W_+$ possesses a complex entire extension to all $t=r+\mi s\in\mathbb{C}$.
Assumptions \emph{(i), (iii),} and \emph{(iv)} allow us to fix
\begin{equation}
\label{t=0}
w_0=\Gamma_0=\Gamma(0)\in W_-^u
\end{equation} 
close enough to $W_-$ in the local unstable manifold $W_-^u$\,, such that Poincaré linearization \eqref{P} and   Poisson stability \eqref{poisson} hold, for all $r\leq 0$.
Let $C:=\sup \|\Gamma(\mathbb{R})\|$. 
Inductively, we will establish rectangle boundaries $0\leq r_n\,,s_n\nearrow+\infty$ on which the heteroclinic orbit $\Gamma$ is bounded by $C+1$, i.e.
\begin{align}
\label{right}
|r|=r_n\,,\  |s|\leq s_n    \quad&\Rightarrow\quad \|\Gamma(r+\mi s)\|\leq C+1 \,, \\
\label{top}
|r|\leq r_n\,,\  |s|= s_n    \quad&\Rightarrow\quad \|\Gamma(r+\mi s)\|\leq C+1 \,. 
\end{align}
Since $\Gamma(t)\in\mathbb{C}$ is assumed to be complex entire, it is complex differentiable in $t=r+\mi s\in\mathbb{C}$.
Therefore, the harmonic maximum principle establishes the same uniform bound $C+1$ in the interior of any rectangle, and hence for all $\Gamma(\mathbb{C})$.
By Liouville's theorem, $\Gamma$ must then be constant.
This contradicts our definition of heteroclinicity, and will complete the indirect proof.

To establish sequences $r_n\,,s_n\nearrow\infty$ of rectangle boundaries, where the bounds $\|\Gamma\|\leq C+1$ of \eqref{right}, \eqref{top} hold, we only consider imaginary parts $s\geq 0$; the case $s\leq 0$ is analogous.
However, the case $0\leq r \leq r_n\nearrow +\infty$ of forward convergence to the locally asymptotically stable equilibrium $W_+$\,, and the case $0\geq r \geq -r_n \searrow -\infty$ of convergence to the non-resonantly unstable equilibrium $W_-$ in backward time, are quite different.

\textbf{Step 1:} $0\geq r \geq -r_n \searrow -\infty$.

Our unstable spectral non-resonance assumption \emph{(iv)} and Poincaré linearization \eqref{P} have allowed us to pick $\Gamma_0=\Gamma(0)$ within the domain of Poincaré linearization \eqref{P}, in the \emph{local} unstable manifold $W_-^u$ of the unstable source equilibrium $W_-$ of $\Gamma: W_-\leadsto W_+$\,.
In particular, we may assume
\begin{equation}
\label{eps1-}
\|\Gamma(r+\mi s)-W_-\| < 1
\end{equation}
for all $r\leq 0$ and all $s\in\mathbb{R}$, in the first place.

This settles estimates \eqref{right}, \eqref{top} in the case $0\geq r \geq -r_n \searrow -\infty$ of backward time. 

\textbf{Step 2:} $0\leq r \leq r_n\nearrow +\infty$.

We start from a few elementary observations concerning asymptotic stability of equilibria $W$.
By standard definition of local asymptotic stability in real time, there exists $0<\eps<1$ such that the following two statements on stability and convergence hold for $w(t)=\Phi^t(w_0)$ and all $0\leq t=r\in\mathbb{R}$:
\begin{align}
\label{eps1}
\|w_0-W\|<\eps  \quad&\Rightarrow\quad\, \|w(r)-W\|<1 \,,  \\
\label{epslim}
\|w_0-W\|<\eps  \quad&\Rightarrow\qquad\ \,  \lim_{r\rightarrow\infty} w(r)=W \,.
\end{align}
Since linear asymptotic stability implies nonlinear local asymptotic stability, assumption \emph{(ii)} of theorem \ref{thmmain} implies that these statements hold true for $W:=W_+$\,, in our heteroclinic setting.
Decreasing the sectorial angle $\Theta>0$, if necessary, we reach the same conclusions for the complex sector $t\in \mathfrak{S}_\Theta$\,.

Starting an induction over $n$ at $r_1=s_1=0$, there is nothing to prove.
We assume $r_{n-1}\,, s_{n-1}$ have been chosen, and we construct $r_n\geq r_{n-1}+1,\ s_n\geq s_{n-1}+1$ next.

For the complex entire heteroclinic orbit $\Gamma:W_-\leadsto W_+$, we first observe
\begin{equation}
\label{lim}
\lim_{r\rightarrow \infty} \Gamma(t_0+r) = W_+\,
\end{equation}
for real $r$ and any fixed $t_0\in \mathbb{C}$\,.
This follows from \eqref{epslim}, strong continuity, and because the semigroups commute in the sector $\mathfrak{S}_\Theta$\,; see \eqref{flow}.

Consider any $w_0=\Gamma(t_0)\in\Gamma(\mathbb{C})$.
Then the following statements hold for all $\tilde{w}_0\in Y$:
\begin{align}
\label{r0W+}
	\forall\, w_0\ \exists\, r_0>0: \qquad\qquad\qquad\qquad\quad\  r\geq r_0 \ &\Rightarrow \ \| w(r)-W_+\| \phantom{_0} <\eps\,;\\
\label{scont}
	\forall\, w_0,\, r_0\,,\,\eps'>0\  \exists\, \delta\phantom{_0}>0: \ \| \tilde{w}_0-w_0\|<\delta,\ 0\leq r\leq r_0 \ &\Rightarrow \ \| \tilde{w}(r)-w(r)\| <\eps';\\
\label{r0eps}
	\forall\,w_0\ \exists\ r_0, \delta>0: \ \, \| \tilde{w}_0-w_0\|<\delta \qquad\qquad\quad\ \,&\Rightarrow \ \| \tilde{w}(r_0)-W_+\| <\eps\,;\\
\label{r0}
	\forall\, w_0\ \exists\ r_0, \delta>0:\,  \ \| \tilde{w}_0-w_0\|<\delta, \ \,\phantom{0\geq}r\geq r_0 \ &\Rightarrow \  \| \tilde{w}(r)-W_+\| \phantom{_0}<1\,.
\end{align}
Indeed, convergence \eqref{lim} implies \eqref{r0W+}.
Strong continuity on compact time intervals implies \eqref{scont}.
Together, \eqref{r0W+} with $r:=r_0$\,, and \eqref{scont} with $\eps':=\eps-\| w(r_0)-W_+\|>0$, imply \eqref{r0eps}.

Combining \eqref{r0eps} with \eqref{eps1} proves \eqref{r0}.
Only this last step uses local asymptotic stability property \eqref{eps1} of $W_+$\,.
Specifically, $w_0$ in \eqref{eps1} is chosen as $\tilde{w}(r_0)$ from \eqref{r0eps}.
\footnote{For unstable targets $W_+$\,, this step \eqref{r0} fails: although $ \tilde{w}(r_0)$ would be an element of the \emph{global} stable manifold $\mathbf{W}_+^s$ of the hyperbolic unstable equilibrium $W_+$\,, near $W_+$\,, it might fail to lie in the \emph{local} stable manifold $W_+^s$\,.
The local stability estimate \eqref{eps1}, however, only applies in the \emph{local} stable manifold $W_+^s$\, of  $W_+$\,.
In step 1, i.e. near $W_-$\,, unstable spectral non-resonance and Poincaré linearization have circumvented this difficulty.}

To establish the required estimate \eqref{top} at a suitable top boundary $s=s_n$\,, we now recall our choice of $\Gamma_0=\Gamma(0)$, in step 1, such that Poisson stability \eqref{poisson} holds for $w_0=\Gamma_0$\,, in imaginary time.
Pick $r_0$ and $\delta$ according to \eqref{r0}, for $w_0:=\Gamma_0$\,.
Reduce $\delta$, if necessary, such that \eqref{scont} holds, in addition, for $w_0:=\Gamma_0$ and $\eps':=1$.
Then
\begin{equation}
\label{delta}
\exists\,\delta>0\,: \quad r\geq 0,\ \| \tilde{w}_0-\Gamma_0\|<\delta \quad \Rightarrow \quad  \|\tilde{w}(r)\| \leq C+1\,.
\end{equation}

Since $w_0=\Gamma_0\in\boldsymbol{\omega}_{\mi\mathbb{R}}(\Gamma_0)$ by \eqref{poisson}, we may choose $s_n \geq s_{n-1}+1$ such that $\tilde{w}_0:=\Gamma(\mi s_n)$ satisfies $\| \tilde{w}_0-\Gamma_0\|<\delta$, in \eqref{delta}.
This proves claim \eqref{top}, at the top boundary, independently of any choice for the right boundary $r_n$\,.

To construct $r_n\geq r_{n-1}+1$, we can now cover the compact set $\Gamma(\mi [0,s_n])$ by a finite collection of $\delta =\delta_j$ neighborhoods of $w_0=w_j$\,, with associated $r=r_{0j}$\,, such that  \eqref{r0} holds for all $j$.
Let $r_n$ also exceed the maximal $r_{0j}$\,. 
Then \eqref{delta} establishes the claimed right boundary estimate \eqref{right}.
This settles the case $0\leq r \leq r_n\nearrow +\infty$ of forward time, and completes the proof of theorem \ref{thmmain}.
\hfill$\bowtie$

We recall how our proof of corollary \ref{cormain} has circumvented the difficulty concerning \eqref{r0}, mentioned in the footnote.
In the fastest unstable manifold of $W_-$\,, of dimension $d=1$, any entire $\Gamma(r+\mi s)$ has to be periodic in $s$ of minimal period $p=2\pi/\mu_0$\,; see \eqref{Gammaperiod} in section \ref{Fastest}.
Therefore we were able to replace the local stability argument \eqref{eps1} by strong continuity and uniform convergence \eqref{rper}, in the bounded period interval $0\leq s\leq p$.
This mended the gap and proved theorem \ref{thmmain}, as before.

In the ODE case $X=\mathbb{C}^N$ of \cite{fiedlerClaudia}, Poincaré linearization and quasi-periodicity at, both, $W_-^u$ and $W_+^s$, came to rescue.
In fact, the proof there was based on Cauchy's theorem, to compare the conflicting Fourier-coefficients with respect to $s$, at fixed $r$ near $\pm\infty$.
This was possible even though the supposedly entire heteroclinic orbit $\Gamma(r+\mi s)$ might have lost quasi-periodicity in $s$, at intermediate $r$.

\section{The quadratic heat equation: Weierstrass equilibria}\label{wp}

As announced in section \ref{Qua}, we now embark on our study of the quadratic heat equation \eqref{heat}.
Specifically, we address spatially non-homogeneous real equilibria $W$, i.e. non-constant solutions $W=W(x)$ of the Neumann ODE boundary value problem
\begin{align}
\label{ODEW}
W_{xx}+6W^2-\lambda &=0,\quad \textrm{for}\ 0<x<\tfrac{1}{2}\,,\\
\label{ODEbc}
W_x &=0,\quad \textrm{for}\ x=0,\tfrac{1}{2}\,;
\end{align}
see \eqref{heat0}. 
Neumann boundary conditions \eqref{ODEbc}, and our preoccupation with real heteroclinic PDE orbits $\Gamma:W_-\leadsto W_+$ between real equilibria $W=W_\pm$\,, require a few adaptations of the standard complex results.
We harmonize the bifurcation viewpoint of figure \ref{fig1} with classical Fourier-expansions of the Weierstrass function $\wp$, for modular lattice parameters $\tau\rightarrow \mi\infty$.
Throughout we refer to \cite{Akhiezer, Fricke, Lang, Lamotke} for details in complex analysis language.

For $m,n\in\mathbb{Z}$ and $\Im \tau>0$, let $\omega=m+n\tau$ denote the elements of the integer lattice $\Lambda:=\langle1,\tau \rangle_\mathbb{Z}$\,, with $\Lambda':=\Lambda\setminus\{0\}$.
Then the associated Weierstrass function $\wp(z)=\wp(\tau;z)$ is defined by the infinite partial fraction
\begin{equation}
\label{defwp}
\wp(\tau;z) := z^{-2}+\sum_{\omega\in\Lambda'}\,\big((z-\omega)^{-2}-\omega^{-2}\big)\,.
\end{equation}
Note how $\wp$ is meromorphic of degree $-2$, and doubly periodic in $z$ with lattice periods $1, \tau$.
Essentially by Liouville's theorem, $\wp=\wp(z)$ satisfies the ODE
\begin{align}
\label{ODE1wp}
\wp_{z}^2 &\phantom{:}= 4\wp^3-g_2\wp-g_3\,,\quad \textrm{for}\ z\in\mathbb{C}\setminus \Lambda,\ \textrm{with coefficients}\\
\label{g2fraction}
g_2 &:=60\phantom{1}\sum_{\omega\in\Lambda'}\,\omega^{-4}   \,,\\
\label{g3fraction}
g_3 &:=140\sum_{\omega\in\Lambda'}\,\omega^{-6}   \,.
\end{align}
Standard calculus of residues on $\wp_z$ locates the three zeros $\omega_j$ of $\wp_z$\,, and hence the zeros $\xi$=$e_j$=$\wp(\omega_j)$ of the cubic polynomial $4\xi^3-g_2\xi-g_3$\,, at the half-periods $\omega_j$ of $z\mapsto\wp(z)$:
\begin{equation}
\label{ej}
\omega_1=\tfrac{1}{2},\quad \omega_2=\tfrac{\tau}{2},\quad \omega_3=\tfrac{1}{2}+\tfrac{\tau}{2}\,.
\end{equation}

\begin{lem}\label{Wn}
Any spatially non-homogeneous real solution $W$\,of the Neumann ODE boundary value problem \eqref{ODEW}, \eqref{ODEbc} takes the form $W=W_n$ at $\lambda=\lambda_n$, where
\begin{align}
\label{lambdan}
\lambda_n &= \tfrac{1}{2} n^4 g_2\,, \qquad\qquad\textrm{and either}\\
\label{Wntop}
W_n(x)&=-n^2\wp(\omega_2 +n x)\,, \qquad\textrm{or}\\
\label{Wnbot}
W_n(x)&=-n^2\wp(\omega_3+n x)\,,
\end{align}
for some purely imaginary modular parameter $\tau=\mi\theta,\ \theta>0$ and $n=1,2,3,\ldots$\,.

Explicit Fourier expansion in $x$ and $\tau$, i.e. in terms of real $h=\pm\exp(-\pi\theta)$, reads
\begin{align}
\label{lambdanh}
\lambda_n &= \tfrac{2}{3}(n\pi)^4 \big(1+240 \sum_{k\geq 1}\, \sigma_3(k) h^{2k}\big)\,, \\
\label{Wnh}
W_n(x)&=(n\pi)^2\big(\eta+8\sum_{k\geq 1}\,k\,\frac{h^k}{1-h^{2k}}\,c_{nk}  \big)\,,\qquad\textrm{where}\\
\label{eta}
\eta&= \tfrac{1}{3}-8\sum_{k\geq 1}\,k\,\frac{h^{2k}}{1-h^{2k}} \,.
\end{align}
Here $\sigma_3(k)$ in \eqref{lambdanh} abbreviates the sum of the third powers of all divisors of k, including 1 and k itself.
In \eqref{Wnh}, the spatial Fourier terms $\cos(2\pi nkx)$ are abbreviated by $c_{nk}$\,.

As in figure \ref{fig1}, note the supercritical pitchfork bifurcations from the spatially homogeneous solution $h=0$, alias $\tau=\mi\infty,\ \theta=+\infty$, at
\begin{equation}
\label{bif}
\lambda_{n0}= \tfrac{2}{3}(n\pi)^4\,,\quad W_{n0}(x)\equiv \tfrac{1}{3}(n\pi)^2 = \sqrt{\lambda_{n0}/6}=W_\infty\,.
\end{equation}
\end{lem}

\begin{proof}
We consider complex solutions $W=W(x),\,x\in\mathbb{R}$ of the Neumann boundary value problem \eqref{ODEW}, \eqref{ODEbc} first.
The two spatially homogeneous solutions are obtained, trivially, as the square roots of $\lambda/6$.
Suppose therefore that $W$ is non-constant in $x$.
By integration of the Hamiltonian pendulum \eqref{ODEW} with respect to $x$, we obtain the conservation
\begin{equation}
\label{ODE1W}
W_{x}^2 = -4W^3+2\lambda W+2E\,,
\end{equation}
of some complex ``energy'' parameter $E\in\mathbb{C}$.
Reversibility of the second order equation with respect to $x\mapsto-x$, and reflection through the Neumann boundary conditions \eqref{ODEbc}, show that $W=W(x)$ must possess (not necessarily minimal) period 1 in $x$.
Let $1/n$ denote the minimal period of $x\mapsto W(x)$, and define
\begin{equation}
\label{defWn}
\tilde{\wp}(x):= -n^{-2}W(x/n)\,.
\end{equation}
Then $\tilde{\wp}(x)$ possesses minimal period 1 in $x$, and solves \eqref{ODE1wp} with
\begin{equation}
\label{g2g3}
g_2=2n^{-4}\lambda\quad\textrm{and}\quad g_3=-2n^{-6}E\,.
\end{equation}
Since $W(x)$ is neither constant, nor homoclinic, the discriminant $\Delta:=g_2^3-27g_3^2$ of the cubic polynomial $4\xi^3-g_2\xi-g_3$ is nonzero.
In particular,  \begin{equation}
\label{tildewp=wp}
\tilde{\wp}(x)=\wp(z_0+x)
\end{equation} coincides with the standard Weierstrass function, for a suitable integer lattice $\Lambda$ and up to a time-shift $z_0\in\mathbb{C}$.
The modular parameter $\tau$ is determined uniquely from $g_2\,,g_3$ in \eqref{g2g3}, up to an $\mathrm{SL}_2(\mathbb{Z})/\{\pm\textrm{id}\}$ fractional linear transformation of $\tau$, via Klein's modular invariant $J=g_2^3/\Delta$.

The Neumann boundary conditions \eqref{ODEbc} and reversibility in $x$, finally, determine boundary values $\tilde{\wp}(x)\in\{e_1,e_2,e_3\}$ at $x=0,1/2$; see \eqref{ej}.
Only from here on, we use that the parameter $\lambda>0$ and the solutions $W$ are assumed to be real.
We will discuss complex solutions in section \ref{Sing}; in particular see \eqref{sigma}.
In view of \eqref{ej}, real $W$ imply $z_0=\omega_2$ or $z_0=\omega_3$.
Indeed, the choice $z_0=\omega_1$ would lead to singular solutions $\tilde{\wp}(x)$ and $W_n(x)$.
Non-horizontal evaluations of $\wp$ would lead to non-real complex solutions; see \eqref{sigma}.
This proves the alternative claims of \eqref{Wntop} and \eqref{Wnbot}, respectively, and \eqref{lambdan}, in either case.

We now address the case \eqref{Wntop}, where $\tilde{\wp}(x)=\wp(\omega_2+x)$ tracks the 1-periodically rescaled version $-n^{-2}W(x/n)$ between $e_2$ and $e_3$\,, for $0<x<1/2$; see \eqref{defWn} and \eqref{ej}.
Standard Fourier expansions of $\wp$ and $g_2$ in terms of 
\begin{equation}
\label{htau}
h=\exp(\pi\mi\tau)
\end{equation} 
then prove claims \eqref{Wnh}, \eqref{eta}, and \eqref{lambdanh}, via the scaling \eqref{defWn}.
Note absolute convergence of the sums: $\Im\tau>0$ implies $0<|h|<1$.
Similar expansions have been derived in \cite{Fricke}, section I.4.9, (14)--(16), and, skipping $g_2$\,, in \cite{Akhiezer}, Table X, p. 204.
In different notation, see \cite{Lang}, section 4.2, and for $g_2$\,, but not $\wp$, also \cite{Apostol}, Theorem 1.18 in section 1.14.
The notationally confusing, but elementary, arithmetic details to derive the specific versions \eqref{lambdanh}--\eqref{eta} from any of these references are left to the diligent reader.

We consider real solutions $W_n$\,, from now on.
We will address case \eqref{Wnbot} after first identifying the modular parameter $\tau$. 
The Fourier expansion \eqref{Wnh} in terms of $c_{nk}$ has to be real now, for all $0<x<1/2$. 
In particular all coefficients must be real.
For $k=1$ and nonzero complex $h$, this implies that there exists $a\in\mathbb{R}$ such that
\begin{equation}
\label{hreal}
h/(1-h^2)=1/a\,.
\end{equation}
The solutions $h=\tfrac{1}{2}(-a\pm\sqrt{a^2+4}$) of the resulting quadratic equation for $h$ are therefore real, too. 
By \eqref{htau} and $0<|h|<1$, this implies the alternative
\begin{align}
\label{htop}
0<h=+\exp(-\pi\theta)<1,\quad &\tau=\mi\theta,\qquad\textrm{or} \\
\label{hbot}
-1<h=-\exp(-\pi\theta)<0,\quad &\tau=\mi\theta+1,
\end{align}
for modular parameters $\tau$ in the upper half plane $\Im\tau>0$, and for  $\Re\tau$ taken $\mathrm{mod}$ 2.
Here $0<\theta<\infty$. 
Comparing the cases $\tau=\mi\theta$ and $\tau=\mi\theta+1$ we see how the roles of $\omega_2$ and $\omega_3$ in \eqref{ej} interchange.
Since \eqref{Wnh} follows from \eqref{Wntop}, for the first case \eqref{htop} and $0<h<1$,
this also proves the same Fourier expansion \eqref{Wnh} by \eqref{Wnbot}, in the second case \eqref{hbot} and for $-1<h<0$.

Note how the two cases $\pm h>0$ correspond to the upper (orange) and lower (blue) branches of the pitchfork bifurcations in figure \ref{fig1}.
Also $\lambda_n\geq\lambda_{n0}$ increases strictly with $h$, globally from the bifurcation point \eqref{bif} at $h=0$, alias $\tau=\mi\infty,\ \theta=+\infty$ in \eqref{lambdanh}--\eqref{eta}.

This proves the lemma.
\end{proof}

For real solutions $W_n$\,, the parameter $E\in\mathbb{R}$ in \eqref{ODE1W} is the classical pendulum energy.
Note positive discriminant $\Delta>0$, since $W_n$ possesses two real extrema: $e_2$ and $e_3$.
In particular $g_2^3=\Delta+27 g_3^2>0$, and hence $\lambda>0$ by \eqref{g2g3}.
More precisely, $\lambda>\lambda_{n0}=\tfrac{2}{3}(n\pi)^4$, by \eqref{lambdanh}, \eqref{bif}.

The appearance of the modular parameter $\tau\ \mathrm{mod}\, 2$ in \eqref{htop}, \eqref{hbot} might surprise some experts (and non-experts, like the authors) who would properly have expected $\tau\ \mathrm{mod}\, 1$, from biholomorphic equivalence of complex tori $\mathbb{C}/\Lambda$.
The discrepancy is caused by our selection of Neumann boundary conditions, which distinguish between boundary values $-n^2e_2$ and $-n^2e_3$ of $W_n(x)$ appearing on the left boundary $x=0$, and thus distinguish between upper and lower branches $h>0$ and $h<0$ in figure \ref{fig1}, for the doubly periodic Weierstrass function $\wp$.

The same procedure which has identified all real solutions of the ODE boundary value problem \eqref{ODEW}, \eqref{ODEbc}, in lemma \ref{Wn}, can easily be modified to also identify all remaining \emph{complex solution branches}.
Indeed, to encounter an interval of real branch parameters $\lambda$, rather than just a totally disconnected set of singletons,
we have to work with modular parameters $h$ such that $\lambda$ remains real, on at least an interval of $h$.
Therefore $h^2=\exp(2\pi\mi\tau)$ must remain real, in the complex Fourier expansion \eqref{lambdanh} of $\lambda=\tfrac{1}{2}n^4 g_2$\,; see \eqref{g2g3}.
In other words, $\tau=\tfrac{1}{2}+\mi\theta$.
This parametrizes the two non-real branches which bifurcate subcritically at the same bifurcation points \eqref{bif}, for purely imaginary $h$.

Our present PDE techniques do not allow us to determine global heteroclinic orbits among complex branches, in real or complex time.
See the discussion section \ref{Sing} for further comments and examples, including the purely quadratic ``vertical'' case of $\lambda=0$.
%There, our branch argument fails due to complex rescaling \eqref{sigma}.

\section{The quadratic heat equation: non-resonance at bifurcation}\label{NonRes}

We continue our study of the quadratic heat equation \eqref{heat} with Fourier expansions of eigenfunctions and eigenvalues, up to order $h^2$ in the Fourier term $h=\exp(\pi\mi\tau)$ of the modular parameter $\tau\rightarrow \mi\infty$.
We prove theorem \ref{thmheat}.
Our starting point is lemma \ref{Wn} on the real equilibrium branches $(\lambda_n\,,W_n)$, parametrized over real $h$.
In view of theorem \ref{thmmain}, proved in section  \ref{Ete}, it only remains to meet assumption \emph{(iv)} of theorem \ref{thmmain}, by exclusion of unstable spectral resonances \eqref{defresonance}.

\textbf{Proof of theorem \ref{thmheat}.}\\
The spectrum of the linearization $L=A+f'(W)$ at $W=W_n$ or $W=W_\infty$ is given by the Sturm-Liouville eigenvalues $\mu=\mu_k$ and eigenfunctions $\phi=\phi_k$ of the Neumann boundary value problem
\begin{align}
\label{linmu}
\mu\,\phi&=\phi_{xx}+12W\phi\,,\quad \textrm{for}\ 0<x<1/2\,,\\
\label{linbc}
0&=\phi_x\,, \qquad\qquad\quad\  \textrm{for}\ x=0,1/2\,.
\end{align}
See \eqref{mu}.
We suppress the subscripts $k$ and $n$, for the moment.

\textbf{Case 1:} $W_-=W_\infty$\\
\nopagebreak
The spatially homogeneous case $W=W_\infty=\sqrt{\lambda/6}$ is explicit.
The unstable spectrum at any $\lambda>0$ is given by
\begin{equation}
\label{muinfty}
0<\mu=\mu_{\infty k} = 2\sqrt{6\lambda}-(2\pi k)^2, \quad \textrm{for}\ 0\leq k^2 < \pi^{-2}\sqrt{\tfrac{3}{2}\lambda}\,,
\end{equation}
with eigenfunctions $\phi=c_k(x):=\cos(2\pi kx)$. The resonance condition \eqref{defresonance} reads
\begin{equation}
\label{inftyres}
(\vert \mathbf{m}\vert-1)\pi^{-2}\sqrt{\tfrac{3}{2}\lambda} = -j^2+\sum_k m_k\, k^2 \leq \vert\mathbf{m}\vert\,\bar{k}^2\,.
\end{equation}
Here $m_k\geq 0,\ \vert m\vert\geq 2$, and $\bar{k}$ denotes the maximal $k$ such that $m_k>0$.
As long as $\bar{k}^2<\pi^{-2}\sqrt{\tfrac{3}{2}\lambda}$ stays away from its upper bound,
we obtain upper bounds on $\vert \mathbf{m}\vert$, and hence discrete sets of resonant $\lambda$.
The limiting cases $\pi^{-2}\sqrt{\tfrac{3}{2}\lambda}\searrow k^2=n^2 $ of vanishing eigenvalue $\mu_{\infty n}=0$, on the other hand, precisely identify the points $\lambda=\lambda_{n0}$ of pitchfork bifurcations; see \eqref{bif} and \eqref{biflambda}.

This establishes the result of theorem \ref{thmheat}, for the case $W_-=W_\infty$.

In fact, it is easy to construct examples of resonant parameters accumulating to the bifurcation points.
Integer $m\nearrow\infty$, for example, provide a sequence of 1:$m$ resonances $\mu_{\infty 0}=m\mu_{\infty n}$ in the sense of \eqref{defresonance}.
Resonance occurs at parameters $\lambda=\lambda_m'\searrow\lambda_{n0}$ defined by the second equality in
\begin{equation}
\label{inftyresex}
\mu_{\infty 0} = 2\sqrt{6\lambda_m'} = m\, \big(2\sqrt{6\lambda_m'}-(2\pi n)^2\big) = m\, \mu_{\infty n}\,.
\end{equation}

\textbf{Case 2:} $W_-=W_n\,,\ 1\leq n\leq 22$\\
\nopagebreak
Since all eigenvalues are algebraically simple and the equilibrium $W=W_n$ depends ana\-lytically on $h$, so do the eigenvalues $\mu$ and properly normalized eigenfunctions $\phi$.
Suppressing standing subscripts $n,k$ of $W,\mu,\phi$, all the way, 
we can therefore expand at the bifurcation point $h=0$ and write
\begin{align}
\label{Wh2}
W&= (n\pi)^2( W_0+W_1h+W_2h^2+\ldots) \,,            \\
\label{muh2}
\mu\,&=\ \, 4\pi^2\ (\,\mu_0\,+\,\mu_1h\,+\,\,\mu_2h^2\,+\ldots) \,,          \\
\label{phih1}
\phi\,&= \phantom{(n\pi)^2(}\, \phi _0\,+\,\phi_1h\,+\ldots\,,
\end{align}
with new indices $0,1,2$ indicating coefficients of the $h$-expansions.
In view of lemma \ref{Wn}, and \eqref{Wnh} in particular, we cover both bifurcating pitchfork branches \eqref{Wntop} and \eqref{Wnbot} by $0<\pm h<1$.

Truncating the expansions \eqref{Wnh}, \eqref{eta}, we obtain
\begin{align}
\label{W0}
    W_0 &= \ \tfrac{1}{3}\,,   \\ 
\label{W1}
    W_1 &=  \qquad\quad  8\, c_n\,,     \\
\label{W2}
    W_2 &=-8 + 16\,c_{2n}\,.
\end{align}

We now insert the $h$-expansions \eqref{Wh2}--\eqref{phih1} into the spectral linearization \eqref{linmu}, \eqref{linbc} and compare coefficients.
At level $h^0$ we get
\begin{equation}
\label{muphi0}
\mu_0=n^2-k^2,\qquad \phi_0=c_k\,.
\end{equation}
Note how our amplitude choice $\phi_0=c_k$ scales the whole eigenfunction expansion.
I.e., the Fourier term $c_k\in \ker(\mu_0-\partial_{xx}-12\,W_\infty)$ can be assumed absent in $\phi_1\,,\ldots$\ .
At level $h^1$, consequently, we have to distinguish the special case $k=n/2$, which only arises for even $n$:
\begin{align}
\label{muphi1}
\mu_1 &= 0\,, & \phi_1 &= 12\, n\big(\tfrac{1}{n-2k}c_{n-k}+\tfrac{1}{n+2k}c_{n+k}\big)\,, &  \textrm{for}\  k\neq n/2, \\
\label{muphi1/2}
\mu_1 &= 12\, n^2\,, & \phi_1&= 6\, c_{3n/2}\,, &  \textrm{for}\ k= n/2.
\end{align}
At level $h^2$, we again have to distinguish the special case $k=n/2$. 
Since our main objective here is the eigenvalue coefficient $\mu_2\, \phi_0$\,, scalar $L^2$-multiplication of \eqref{linmu} by $\phi_0$ eliminates the eigenfunction term $\phi_2$ as irrelevant for $\mu_2$\,. 
We obtain
\begin{align}
\label{muphi2}
\mu_2 &= 24\,n^2\cdot\tfrac{11n^2+4k^2}{n^2-4k^2}\,,   \quad \textrm{for}\  k\neq n/2\,, \\
\label{muphi2/2}
\mu_2 &= 48\, n^2\,,\phantom{\cdot\tfrac{11n^2+4k^2}{n^2-4k^2}} \quad\ \,  \textrm{for}\ k= n/2\,.
\end{align}

Any linear resonance conditions \eqref{defresonance}, as well as the eigenvalues $\mu_j\,,\mu_k$ themselves, are analytic in $h$.
Therefore, any unstable resonance either holds identically in $h$, or at most for a discrete set of $h$, and hence of parameters $\lambda$ by expansion \eqref{lambdanh}.
Since the pitchfork bifurcations are supercritical, we can extend this argument to include the bifurcation points $\lambda=\lambda_{n0}$ themselves, this time.
Indeed, the stable eigenvalue $\mu_n\nearrow 0$, for $\lambda\searrow\lambda_{n0}$\,, does not contribute to unstable resonance, this time.

We have to exclude identical resonance.
In other words, it is sufficient to exclude resonances \eqref{defresonance}, separately, for each set of expansion coefficients $\mu_0\,, \mu_1\,, \mu_2$\,, and with the same integer vector $\mathbf{m}$.
Computer assisted checking of the resulting three explicit linear Diophantine conditions for $\mathbf{m}$ excluded identical resonance, for $n\leq 22$.
(Disclaimer: we do not claim failure of our crude second order check at $n=23$.)

This completes the proof of theorem \ref{thmheat}. \hfill $\bowtie$

\section{The quadratic heat equation: fast unstable manifolds}\label{Wuu}

We continue our study of the quadratic heat equation \eqref{ODEW} with Neumann boundary \eqref{ODEbc}.
We address heteroclinic orbits $\Gamma:W_-\leadsto W_+$ which originate in the local $d$-dimensional fast unstable manifold $W_-^{uu}$ of $W_-$\,.
See section \ref{Fast}.
The hyperbolic equilibrium $W=W_-$ can be homogeneous, i.e. $W=W_\infty$\,\,, or non-homogeneous, i.e. $W=W_n$\,.
In either case, we consider dimensions $1\leq d<i(W)$.
In particular, we prove theorem \ref{thmheatd} on unstable spectral non-resonance, for $d<1+\tfrac{1}{\sqrt{2}}\,i(W)$; see \eqref{dn}.

Because the Sturm-Liouville spectrum at $W$ consists of algebraically simple eigenvalues, the local $d$-dimensional fast unstable manifold $W_-^{uu}$ always exists. It is characterized by exponential decay rates, in backward time, which exceed the eigenvalue $\mu_d>0$; see \eqref{mu} and \eqref{linmu}, \eqref{linbc}.
The standard Perron proof for $W^u$ in \cite{Henry}, adapted for decay rates by exponentially weighted spaces, establishes local analyticity of $W_-^{uu}$.

\textbf{Proof of theorem \ref{thmheatd}.}\\
\nopagebreak
Our proof follows the lines of section \ref{NonRes}, replacing the bound on $n:=i(W)$ by a bound on $d/n<1$.
Our treatment is identical for the cases $W=W_n$ and $W=W_\infty$\,, because our assumption $d<n=i(W)$ eliminates the small positive eigenvalue $\mu_{\infty n}$ generated, supercritically, along the primary branch $W_\infty$ at pitchfork bifurcation $\lambda\searrow\lambda_{n0}$\,.
Compare section \ref{NonRes}, case 1.

In fact we just show absence of identical resonances at order $h^0$, in \eqref{Wh2}--\eqref{phih1}. 
By \eqref{muphi0}, that resonance condition can be estimated by
\begin{equation}
\label{mu0resonance}
n^2\geq n^2-j^2=\ \sum_{k=0}^{d-1}\,m_k\,  (n^2-k^2) \geq 2\, (n^2-(d-1)^2)\,,
\end{equation}
for any choice of nonnegative integers $m_k$ of order $\vert\mathbf{m}\vert\geq 2$. This implies
\begin{equation}
\label{dresonance}
d\geq1+\tfrac{1}{\sqrt{2}}\,n\,,
\end{equation}
i.e. the opposite of assumption \eqref{dn}.
By contraposition, \eqref{dn} prohibits identical unstable resonances and proves theorem \ref{thmheatd}. \hfill $\bowtie$

As announced in section \ref{Fast}, we now indicate asymptotic optimality of the bound \eqref{dn} for non-resonance at level $h^0$.
Choosing $j=0$ and $\vert\mathbf{m}\vert=(0,\ldots,1,1)\in\mathbb{N}_0^d$\,, we just have to solve resonance condition \eqref{defresonance} for Pythagorean triples
\begin{equation}
\label{pythagoras}
(d-2)^2+(d-1)^2=n^2.
\end{equation}
Explicit resonances are generated by the successive co-prime convergents $a/b$ of the 1-periodic continued fraction 
\begin{equation}
\label{contfrac}
1+\sqrt{2}=[2;2,2,2,\ldots]\,.
\end{equation}
The integers $a,b$ are related by a linear recursion of Fibonacci type.
Indeed we may choose $n=a^2+b^2$ and $\{d-2,d-1\}=\{a^2-b^2,\,2ab\}$ to achieve $\lim\, n/d = \sqrt{2}$.
Explicit examples of worst case resonances \eqref{dresonance} are $a/b=5/2,\ n=29,\ d=22$ and $a/b=12/5,\ n=169,\ d=121$. And so on.

\section{The quadratic Schrödinger equation: blow-up and blow-down} \label{HeaSchro}

Following up on section \ref{Schro}, we prove theorem \ref{thmpsi} on blow-up $0<s\nearrow s^*$ (and symmetric blow-down $0>s\searrow -s^*$) in the complex Schrödinger equation \eqref{psi}.
Blow-up will be shown to occur for certain initial conditions $\psi_0:=\Gamma(r_0)$ on real eternal heteroclinic solutions $\Gamma:W_-\leadsto W_0$ of the quadratic heat equation \eqref{heat}.
We assume the equilibrium $W_-$ is non-resonant, and $\Gamma$ emanates from the unstable or strong unstable manifold of $W_-$ under the assumptions of theorems \ref{thmheat} or \ref{thmheatd}, respectively.

Our construction of Schrödinger blow-up for $s\mapsto\psi(s)=\Gamma(r_0+\mi s)$, at imaginary time $0<s=s^*<\infty$, will be based on blow-up of $r\mapsto w(r):=\Gamma(r+\mi s_0)$ for the heat equation, at real time $r\nearrow r^*(s_0)<\infty$.
Here $r$ carries sectorial analyticity, whereas $s$ only knows strong continuity.
We then seek Schrödinger blow-up at $s^*=s_0$\,, for $r_0:=r^*(s_0)$.
In other words, we aim for a setting where $r_0=r^*(s_0)$ implies $s_0=s^*(r_0)$.

Starting from just any $s_0>0$ such that $r^*(s_0)<\infty$, however, may run into trouble.
In fact, the real heat semiflow $\Phi^r$ and the imaginary Schrödinger semiflow $S^s=\Phi^{\mi s}$ may fail to commute, up to $r=r_0\,,\ s=s_0$\,, when other singularities of $\Gamma$ to the lower left eclipse full access from $t=0$ to $r_0+\mi s_0$\,.
See the parallelogram condition for commutativity \eqref{flow}.

In \eqref{commute}, our proof of theorem \ref{thmpsi} will avoid this obstacle by the construction \eqref{rectangle} of a left-infinite rectangle $\mathcal{R}\subset\mathbb{C}$ which fully supports complex analyticity of $\Gamma$, except at the upper right corner $r_0+\mi s_0$ where heat \emph{and} Schrödinger blow-up will occur in unison.

\textbf{Proof of theorem \ref{thmpsi}.}\\
By theorems \ref{thmheat} and \ref{thmheatd}, $\Gamma$ cannot be complex entire.
Therefore, real-time blow-up \eqref{defblow-up} of $r\mapsto \Phi^{r+\mi s}(\Gamma_0)=\Gamma(r+\mi s)$ has to occur, for some fixed $s=\hat{s}>0$ and at some finite real time $r\nearrow r^*(\hat{s})<\infty$.
Here we have used reversibility \eqref{reversibility} to restrict to $\hat{s}>0$ and blow-up, 
%rather than $\hat{s}<0$ and blow-down, 
without loss of generality.
We recall here, and will peruse, that due to non-resonance at $W_-$ the heteroclinic solutions $\Gamma(r+\mi s)$ emanating from $W_-$ via $W_-^u$ or $W_-^{uu}$ are complex analytic and remain local, for all $r\leq 0$; see \eqref{eps1-}.

As in \eqref{D}, let $-\infty<r<r^*(s)\leq\infty$ denote the maximal interval of existence of the solution $r\mapsto\Gamma(r+\mi s)=\Phi^r(\Gamma(\mi s))$, in real time $r$.
Note complex analyticity of $\Gamma(t)$ in the domain $t\in\mathcal{D}:=\{r+\mi s\,|\,-\infty<r<r^*(s)\}$.
By strong continuity of the sectorial semigroup $\Phi^t$, the domain $\mathcal{D}$ is open and, equivalently,  $s\mapsto r^*(s)$ is lower semicontinuous.
To reach a contradiction, as announced above, we will construct a left-infinite rectangle
\begin{equation}
\label{rectangle}
\mathcal{R}:=\{r+\mi s\,|\,-\infty<r\leq r_0<\infty,\ 0\leq s\leq s_0<\infty\} \subset \mathcal{D}\cup\{r_0+\mi s_0\},
\end{equation}
such that only the upper right corner $r_0+\mi s_0$ is located on the boundary $\partial \mathcal{D}$.
For the special case of fastest heteroclinic orbits $\Gamma$, see also section \ref{Fastest}.

The significance of the rectangle $\mathcal{R}$ is that the real heat semiflow $r\mapsto\Phi^r$ and the imaginary Schrödinger flow $s\mapsto S^s=\Phi^{\mi s}$ commute on $\Gamma$ for $r+\mi s \in \mathcal{R}\cap \mathcal{D}$\,:
\begin{equation}
\label{commute}
\Phi^r(\Gamma(\mi s)) =\Phi^rS^s(\Gamma_0)=\Gamma(r+\mi s)=S^s\Phi^r(\Gamma_0)= S^s(\Gamma(r)),
\end{equation}
for $r+\mi s\in\mathcal{R}\setminus\{r_0+\mi s_0\}$\,.
See also \eqref{flow}.
We defer the construction of $\mathcal{R}$, for a moment.

Let $0\leq s < s^*(r)\leq\infty$ denote the maximal forward interval of existence of the strongly continuous local Schrödinger flow $S^s(\Gamma(r))$.
For $r=r_0$\,, we claim Schrödinger blow-up at $s^*(r_0)=s_0<\infty$\, i.e. at the upper right corner of the rectangle $\mathcal{R}$.
Indeed $s^*(r_0)\geq s_0$\,, by our construction of $\mathcal{R}$.
Suppose $s^*(r_0)> s_0$\,.
Then $\|\Gamma(r+\mi s)\|$ remains uniformly bounded in some neighborhood of the upper right corner $r_0+\mi s_0$\,, by lower semicontinuity of $r\mapsto s^*(r)$.
Blow-up \eqref{defblow-up} of  $r\mapsto\|\Gamma(r+\mi s)\|$ at $r_0=r^*(s_0)$ under $\Phi^r$, on the other hand, shows unboundedness of the same $H^1$-norm $\|\Gamma(r+\mi s_0)\|$, for $r\nearrow r_0=r^*(s_0)$. 
This contradiction shows Schrödinger blow-up of $\psi(s):=\Gamma(r_0+\mi s)$ at $s=s^*(r_0)=s_0$\,, which will complete the proof.
In short: $r_0=r^*(s_0)$ implies $s_0=s^*(r_0)$, for our rectangle $\mathcal{R}$.

It only remains to construct $\mathcal{R}$ of \eqref{rectangle}, with only the upper right corner $r_0+\mi s_0$ in $\partial \mathcal{D}$.
Clearly we will have to choose $r_0=r^*(s_0)$.
To construct $s_0$\,, we may restrict attention to any compact interval $s\in I:=[0,\hat{s}]$, where the imaginary level $s=\hat{s}$ is already known for blow-up of $r\mapsto\Gamma(r+\mi s)$ in finite real time $r^*(s)$.
The lower semicontinuous function $s\mapsto r^*(s)$ must attain its finite minimum, on some nonempty compact subset of $0<\underline{s}_0\in I$.
Indeed, $\underline{s}_0=0$ does not qualify, because $r^*(0)=\infty$ on the real heteroclinic orbit $\Gamma(r)$.
Define $s_0:=\min \underline{s}_0$\, and $r_0:=r^*(s_0)$.
Then the rectangle $\mathcal{R}$ with upper right corner $r_0+\mi s_0$ satisfies \eqref{rectangle}.

This completes the proof of theorem \ref{thmpsi}. \hfill $\bowtie$

\section{Discussion}\label{Dis}

\subsection{Overview}\label{Over}

We conclude with a cursory collection of additional aspects, and shortcomings, of our results.
Above, we have defined blow-up in the topology of fractional power spaces $X^\alpha$ and $H^1$ Sobolev norms; see \eqref{defblow-up} and \eqref{blow-up}.
In section \ref{Sup} we show how this implies blow-up \eqref{supblow-up} in $L^\infty$ sup-norm. Section \ref{ODEd} returns to the spatially homogeneous, polynomial ODE case $\dot w = f(w)$ in some detail, generalizing the introductory remarks on quadratic $f$ of section \ref{ODE2}.
Based on extensive work in \cite{WolfrumCorr}, for the real case, section \ref{Boundary} explores the boundary of the fastest unstable manifold in complex time.
Blow-up and quasi-periodicity are discussed.
Section \ref{Sing} highlights a few results from the literature on PDE blow-up in the parabolic and Schrödinger settings \eqref{heat} and \eqref{psi}, mostly for the purely quadratic case $\lambda=0$.
After due celebration of the ground-breaking and pioneering results \cite{Masuda1,Masuda2} by Kyûya Masuda, we comment on exciting recent results by Jonathan Jaquette and co-workers, in some detail.
As a fundamentally new tool, they have brought rigorous numerics and computer-assisted proofs to the field.
Under periodic boundary conditions in $x$, section \ref{Trav} provides a glimpse at the construction of traveling wave solutions $w(s,x)=w(\zeta)$ in complex time $\zeta$.
Mixing real and imaginary traveling wave time, e.g.~as $\zeta=cr+\mi x $, we also construct finite time blow-up with a traveling wave flavor.
In section \ref{1000}, we follow \cite{fiedlerClaudia} to describe the potential relevance of complex entire homoclinic orbits, even in ODEs, for the elusive phenomenon of ultra-exponentially small homoclinic splittings and ultra-invisible chaos.
The existence of such homoclinic orbits remains a 1,000 \euro\ open question.
We also recall a reversible scalar ODE example of second order.
The example features a non-entire homoclinic orbit, and a coexisting complex entire (in fact a cosine) periodic orbit.
As a consequence, ultra-exponential homoclinic splitting may fail, but ultra-exponentially  sharp Arnold tongues do occur under suitable, rapidly periodic forcing.

\subsection{Blow-up in sup-norm}\label{Sup}

We show how the $H^1$ Schrödinger blow-up asserted in theorem \ref{thmpsi} actually implies $L^\infty$ blow-up.

To align the settings of quadratic heat and Schrödinger equations \eqref{heat}, \eqref{psi}, let $0\leq r \nearrow r^*<\infty$ denote $H^1$ blow-up, analogously to \eqref{defblow-up}, \eqref{blow-up}, for complex time rays $t=r\exp(\mi\theta)$ along any fixed inclined direction $|\theta|\leq\pi/2$.
In other words, consider solutions $w(r)$ of
\begin{equation}
\label{wtheta}
e^{-\mi\theta} w_r = w_{xx}+6 w^2-\lambda\,,
\end{equation}
such that
\begin{equation}
\label{blow-uptheta}
\|w(r)\|_{H^1}\rightarrow\infty \quad\textrm{for}\quad r\nearrow r^*\,.
\end{equation}
A more standard notion requires blow-up in sup-norm, instead, i.e.
\begin{equation}
\label{supblow-up}
\lim\sup \|w(r)\|_{L^\infty} = \infty \quad\textrm{for}\quad r\nearrow r^*\,.
\end{equation}
We claim that the two conditions \eqref{blow-uptheta}, \eqref{supblow-up} are equivalent, in our setting.

Indeed, standard semigroup theory and the Sobolev embedding $H^1 \hookrightarrow L^\infty$ show that \eqref{supblow-up} implies \eqref{blow-uptheta}.

To establish the converse direction, we argue by contraposition: suppose $\|w(r)\|_{L^\infty} \leq C$ violates \eqref{supblow-up}, remaining uniformly bounded, for $0\leq r<r^*$.
For classical solutions, we differentiate \eqref{wtheta} with respect to $x$, and multiply with the conjugate complex $\overline{w}_x$\,. 
Integration over $x$, by parts, then provides the following differential inequality for the $L^2$-norm $\|w_x\|_{L^2}$\,:
\begin{equation}
\label{H1estimate}
    \tfrac{1}{2}\,\partial_r\,\|w_x\|_{L^2}^2 = -\cos\theta\,\|w_{xx}\|_{L^2}^2 + 12C\, \|w_x\|_{L^2}^2 \leq 12C\,  \|w_x\|_{L^2}^2\,.  \\
\end{equation}
%Note invariance of estimate \eqref{H1estimate} under the scaling \eqref{nscale}, by a homogeneous factor $n^8$.
The differential inequality implies at most exponential growth, and hence uniform bounds in $H^1$, for inclinations $|\theta|\leq\pi/2$ and $0\leq r \leq r^*$. 
In weak form, the same $H^1$-bound extends to mild solutions.
This completes the contraposition, by violation of \eqref{blow-uptheta}.
The same argument extends to scalar nonlinearities $f(x,w)$ with bounded derivatives, locally uniformly in $w$.

With the particular choices $\theta=\mp\pi/2,\ r=s$, and $w=\psi$, theorem \ref{thmpsi} then asserts $L^\infty$ Schrödinger blow-up \eqref{supblow-up} of $\psi(s)$, and blow-down, for real $\psi(0) := \Gamma(r_0)$ and at finite times $s\rightarrow \pm s^*$.

\subsection{ODEs revisited}\label{ODEd}

To highlight the peculiarity of the quadratic nonlinearity $f(w)=w^2-1$ studied so far, we now collect some facts concerning solutions $w=w(t)$ of polynomial scalar ODEs
\begin{equation}
\label{ODEw}
\dot{w}=f(w)=(w-e_0)\cdot\ldots\cdot(w-e_{d-1})=w^d+\ldots+f_0
\end{equation}
in complex time $t=r+\mi s$.
For the quadratic case $d=2$, recall section \ref{ODE2}.
Any polynomial ODE of degree $d$ may be brought to this univariate form with normalized top coefficient, of course, by a suitable complex scaling of $w$.
We assume the $d$ complex zeros $e_j$ of $f$ to be simple.
In other words, $W=e_j$ are the $d$ distinct equilibria of \eqref{ODEw}.
We summarize results from \cite{fiedlerShil}.

For elementary integration of \eqref{ODEw} see section \ref{Integration}.
Section \ref{Riemann} takes the Riemann surface viewpoint of classical complex analysis.
From section \ref{Infinity} onwards, we distinguish real time and study Poincaré compactification of \eqref{ODEw} at $w=\infty$.
We also provide a glossary for the remaining ODE discussion.
Section \ref{PerHet} addresses and excludes periodic and saddle-saddle heteroclinic orbits.
This leads to a classification of connection graphs and orbit equivalence of the remaining Morse flows in section \ref{Conn}; see theorem \ref{Classthm}.
We conclude our ODE excursion with the cyclotomic examples $d=3,4$ in section \ref{Cyclotomic}.

\subsubsection{Elementary integration}\label{Integration}

Separation of variables and elementary partial fractions, or the residue theorem, provide explicit solutions
\begin{equation}
\label{sovw}
t \ =\ \int^{w} \omega \ = \ \sum_{j=0}^{d-1}\,\tfrac{1}{f'(e_j)} \,\log(w(t)-e_j)\,.
\end{equation}
The integrand $\omega$ abbreviates the differential form $dw/f(w)$, which is holomorphic for
$w$ in the punctured Riemann sphere
\begin{equation}
\label{Cd}
w\in\widehat{\mathbb{C}}_d:=\widehat{\mathbb{C}}\setminus\{e_0,\ldots,e_{d-1}\}\,.
\end{equation}
In terms of $z=1/w$, in fact, we obtain the locally convergent power series
\begin{equation}
\label{sovz}
\begin{aligned}
   t &\,=\, \int_\infty^{w}\,\omega\,=\, -\tfrac{1}{d-1}z^{d-1}(1+c_1z+\ldots)
\,,
\end{aligned}
\end{equation}
centered at $z=0$, alias $w=\infty$, for some coefficients $c_1,\ldots$ .
This yields the \emph{maximal analytic continuation}
\begin{equation}
\label{alltw}
t\ \in\ \mathrm{range}\,\mathcal{T}+ \int_\infty^w\, \omega\,.
\end{equation}
Here the \emph{period homomorphism} on the fundamental group $\pi_1(\widehat{\mathbb{C}}_d)$,
\begin{equation}
\label{T}
\begin{aligned}
\mathcal{T}: \pi_1(\widehat{\mathbb{C}}_d) &\rightarrow \mathbb{C} \\
							\gamma &\mapsto \int_\gamma\,\omega
\end{aligned}
\end{equation}
captures the ambiguity  in \eqref{sovw}, which arises from the logarithms and the precise integration paths.
Note that the fundamental group $\pi_1=\pi_1(\widehat{\mathbb{C}}_d)$ coincides with the free group $\mathbb{F}_{d-1}$ on $d-1$ generators, viz. small loops around the equilibrium punctures $e_j$ with $0<j<d$.
Because $(\mathbb{C},+)$ is Abelian, $\ker \mathcal{T}$ contains the commutator subgroup $\pi_1'$ of $\pi_1$\,.
The homomorphism theorem implies that $\pi_1'\trianglelefteq\ker \mathcal{T}\trianglelefteq \pi_1$ are normal subgroups, and
\begin{equation}
\label{Thom}
\mathrm{range}\,\mathcal{T}\, = \,2\pi\mi\, \langle \eta_1,\ldots,\eta_{d-1}\rangle_\mathbb{Z}\,
\cong\, \pi_1/\ker \mathcal{T} \,\cong\, (\mathbb{Z}^k,+)\,,
\end{equation}
where $0<j,k<d$, and $\eta_j:=1/f'(e_j)\neq 0$ denote the reciprocal derivatives at $e_j$\,.
In \cite{fiedlerShil}, theorem 1.1, we have observed that $\eta_j$ can be prescribed arbitrarily, for suitable choices of $e_j$, subject to the two constraints
\begin{align}
\label{sumseta}
\sum_{j=0}^{d-1}\,\eta_j\,&=\,0\,,\quad\mathrm{and}\quad \\
\label{sumsetaJ}
 \sum_{j\in J}\,\eta_j\,&\neq\,0
\end{align}
for any nonempty subset $\emptyset\neq J \subsetneq\{0,\ldots,d-1\}$.

\subsubsection{The Riemann surface}\label{Riemann}

For further background on tools from complex analysis see for example \cite{Ilyashenko, Forster, Jost, Lamotke}.
The maximal analytic continuation \eqref{alltw} of the integral $t$ of the differential form $\omega=dw/f(w)$ on $w\in\widehat{\mathbb{C}}_d$ defines a Riemann surface $\mathcal{R}$, which we represent as the set
\begin{equation}
\label{Rdef}
\mathcal{R} := \{(w,t)\in \widehat{\mathbb{C}}_d\times\mathbb{C}\,\vert\, \eqref{alltw}\ \mathrm{holds}\, \}\,.
\end{equation}
In the topology of the Riemann surface $\mathcal{R}$, however, the fibers over any $w_0\in\widehat{\mathbb{C}}_d$ are discrete.
More precisely, let $(w_0,t_0)\in\mathcal{R}$.
Then the difference $T=t_1(w;w_0)-t_2(w;w_0)$ of any two local expansions $t_1,t_2$ of $t$, centered at the same $w_0\in\widehat{\mathbb{C}}_d$\,, is a period $T\in\mathrm{range}\,\mathcal{T}\cong\mathbb{Z}^k$ of the period map $\mathcal{T}$; see \eqref{Thom}.

For polynomials $f$ of degree $d>2$, the two Riemann surfaces $(w,t)\in\mathcal{R}$ and $w\in\widehat{\mathbb{C}}_d$ are hyperbolic, with the complex upper half plane $\mathbb{H}$ as the shared universal cover:
\begin{equation}
\label{HRCd}
\mathbb{H}\xrightarrow{\mathbf{p}}\mathcal{R}\xrightarrow{\mathbf{q}_w}\widehat{\mathbb{C}}_d\,.
\end{equation}
The projection $\mathbf{q}_w(w,t)=:w$ and the universal coverings $\mathbf{p},\,\mathbf{p}\circ \mathbf{q}_w$ are unbranched, unlimited, normal covering maps.
This describes \eqref{HRCd} as principal fiber bundles with discrete fibers given by the corresponding deck transformation groups
\begin{align}
\label{deckHRCd}
    \mathrm{deck}(\mathbf{p}\circ \mathbf{q}_w)\,&\cong\, \pi_1(\widehat{\mathbb{C}}_d) \,\cong\, \mathbb{F}_{d-1}\,,   \\
\label{deckHR}
    \mathrm{deck}(\mathbf{p})\,\,\   &\cong\,\, \pi_1(\mathcal{R}) \;\cong\,\ker\,\mathcal{T}\,,   \\
\label{deckRCd}
    \mathrm{deck}(\mathbf{q}_w)\; &\cong\, \mathbb{F}_{d-1}/\ker\,\mathcal{T} \,\cong\, \mathrm{range} \,\mathcal{T}\,=\,2\pi\mi\,\langle \eta_1,\ldots,\eta_{d-1}\rangle_\mathbb{Z}\,.
\end{align}

Global integration \eqref{alltw} amounts to the covering projection $\mathbf{q}_t(w,t):=t$,
\begin{equation}
\label{RC}
\mathbf{q}_t:\ \mathcal{R}\rightarrow\mathbb{C}\,.
\end{equation}
Branching points are all periods $t\in\mathrm{range}\,\mathcal{T}$, each of branching multiplicity $d-1$, with associated ramification point $(\infty,t)\in\mathcal{R}$.
The topological closure $\mathbb{B}$ of the branching set $\mathrm{range}\,\mathcal{T}\subset\mathbb{C}$ depends on the reciprocal derivatives $\eta_j=1/f'(e_j)\neq 0$.
This makes the closure $\mathbb{B}$ real linear equivalent to one the following five options
\begin{equation}
\label{closT}
\mathbb{Z},\quad \mathbb{Z}\times \mathbb{Z},\quad \mathbb{R},\quad \mathbb{R}\times \mathbb{Z},\quad \mathbb{C}\cong \mathbb{R}^2\,.
\end{equation}
With $\eta_j$ realizable at least as prescribed in \eqref{sumseta}, \eqref{sumsetaJ}, all five cases actually arise for suitable polynomial nonlinearities $f$.
In the topology of \,$\widehat{\mathbb{C}}_d\times\mathbb{C}$, the description \eqref{Rdef} is an embedding of the Riemann surface $\mathcal{R}$ in the first two cases, only.
These are also the only cases where the set of branching points is discrete and, therefore, where the projection $\mathbf{q}_t$ qualifies as an unlimited branched covering. 
In the remaining three cases, branching points are densely accumulating in $\mathbb{B}$, and we only obtain an immersion $(w,t)\in\mathcal{R}\subset\widehat{\mathbb{C}}_d\times\mathbb{C}$ of the densely accumulating local $t$-sheets defined by $t=t(w;w_0)+T$, for periods $T\in\mathrm{range}\,\mathcal{T}$.

\begin{figure}[t]
\centering \includegraphics[width=\textwidth]{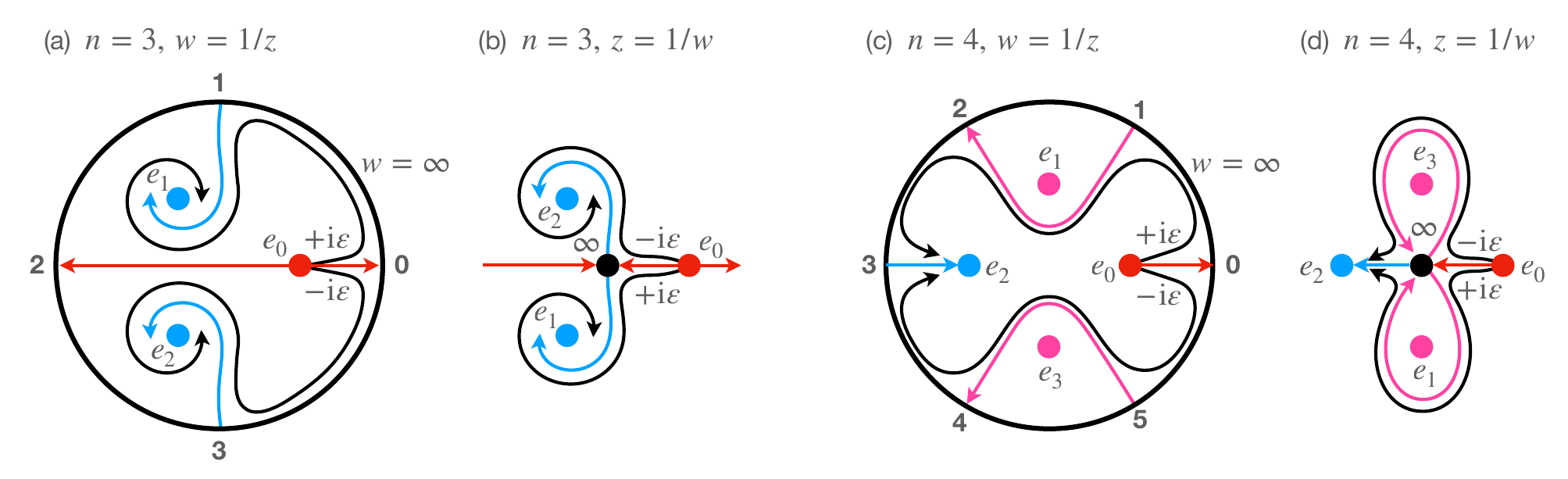}
\caption{\emph{
Schematic phase portraits, in real time, of complex-valued ODEs \eqref{ODEw}, \eqref{rho},  \eqref{alpha} for cyclotomic vector fields $\dot{w}=w^d-1$; see \eqref{cyclotomic}.
For $d=3$ see $w$ in (a), and $z=1/w$ in (b).
Similarly, $w$ in (c) and $z=1/w$ in (d) refer to $d=4$.
The invariant circle $\rho=|z|=0,\ \alpha\in\mathbb{S}^1$ of the Poincaré compactification at $w=\infty$ is marked black in (a), (c).
Interior equilibria $e_j\in\mathbb{D}$ are stable sinks (blue), Hopf centers (purple) surrounded by families of periodic orbits, or unstable sources (red).
Unstable blow-down orbits (blue) emanate from $w=\infty$, alias $z=1/w=0$, at odd-labeled saddles $\mathbf{k}=\mathbf{1},\mathbf{3},\mathbf{5}$.  
Stable blow-up orbits (red) run towards the even-labeled saddles $\mathbf{k}=\mathbf{0},\mathbf{2},\mathbf{4}$.
In (c), two pairs of interior saddle separatrices coincide (purple) in each of the interior saddle-connections $\mathbf{1}\leadsto \mathbf{2}$ and $\mathbf{5}\leadsto \mathbf{4}$.
In (d), the two purple orbits become homoclinic to $w=\infty$.\\
The black trajectories marked $\pm\mi\eps$ are complex perturbations of the red real blow-up orbits $e_0\leadsto \mathbf{0}$, say with initial conditions $w(0)=2\pm\mi\eps$.
For $d=3$ in (a), they closely follow the heteroclinic chains $e_0\leadsto \mathbf{0}\leadsto \mathbf{1}\leadsto e_1$ and $e_0\leadsto \mathbf{0}\leadsto \mathbf{3}\leadsto e_2$, respectively.
The intermediate boundary connections $\mathbf{0}\leadsto \mathbf{1}$ and $\mathbf{0}\leadsto \mathbf{3}$ within the invariant boundary circle $\mathbb{S}^1$ concatenate initial red blow-up to terminal blue blow-down.
In the polar view (b), centered at $w=\infty,\ z=0$, the boundary connections are conflated into $z=0$.
Note the \emph{markedly distinct limits of the two perturbations}, for $\eps\searrow 0$, given by the two distinct remaining blow-up-down concatenations $e_0\leadsto \infty\leadsto e_1$ and $e_0\leadsto \infty\leadsto e_2$\,.\\
In case $d=4$ (d), the initial red blow-up and terminal blue blow-down limits of both perturbations $\pm\eps$ coincide.
Each limit $e_0\leadsto \infty\leadsto\infty\leadsto e_3$ contains an additional interior blow-down-up part $\infty\leadsto\infty$ (purple) which is homoclinic to $z=0$.
The two counterrotating purple homoclinic orbits, however, are markedly distinct: their lobes surround the Hopf centers $e_1$ and $e_3$, respectively, in opposite orientation.
In (c), this is manifest by the two purple interior saddle-connections $\mathbf{1}\leadsto \mathbf{2}$ and $\mathbf{5}\leadsto \mathbf{4}$.
Together with the heteroclinic boundary connections within the black invariant boundary circle $\mathbb{S}^1$, which run between the same adjacent saddles in opposite direction, we obtain two heteroclinic cycles.
Their interiors are properly foliated by families of nested, synchronously iso-periodic orbits of minimal periods $\mp\pi/2$ around the Hopf centers $e_1$ and $e_3$\,; see section \ref{PerHet}.
}}
\label{fig3}
\end{figure}

On $t\in\mathbb{C}$, i.e. after projection by $\mathbf{q}_t$\,, the local flow $\Phi^t$ defined by our original ODE \eqref{ODEw} simply acts by complex time shift.
At the ramification points $(w,t)=(\infty,t_0)\in\mathcal{R}$ of branch points $t_0\in\mathrm{range}\,\mathcal{T}$ of multiplicity $d-1$, however, the regularization \eqref{sovz} shows how $d-1$ different real-time trajectories $r\mapsto z(t_0+r)$ cross each other, at $r=0$. 
See also section \ref{Cyclotomic} and figure \ref{fig3} (b),(d).
In conclusion, real time shift in $t$ does not lift to a real-time flow on the Riemann surface $\mathcal{R}$, for polynomial degrees $d>2$.

Outside $w=\infty$, the covering projections $\mathbf{q}_w$ and $\mathbf{p}\circ \mathbf{q}_w$  lift the local flow $\Phi^t$ on $w\in\widehat{\mathbb{C}}_d\setminus\{\infty\}$ to local flows on the appropriate subsets of $\mathcal{R}$ outside ramification, and of the universal cover $\mathbb{H}$, respectively.
The lifted local flows $\Phi^t$ become equivariant with respect to the groups $\mathrm{deck}(\mathbf{q}_w)\cong \mathrm{range}\, \mathcal{T}$ and $\mathrm{deck}(\mathbf{p}\circ \mathbf{q}_w)\cong\mathbb{F}_{d-1}$; see \eqref{deckRCd}, \eqref{deckHRCd}.

Only in the quadratic case $d=2$ of section \ref{ODE2}, we obtain an unbranched global flow on $\mathcal{R}$.
For the hyperbolic upper half plane $\mathbb{H}$ in \eqref{HRCd}, we just have to substitute the universal cover $\mathbb{C}$ to obtain the parabolic cylinder $w\in\widehat{\mathbb{C}}_2$\,.
The period map features a single purely imaginary period $2\pi\mi/f'(e_1)=-\pi\mi$, in case $f(w)=w^2-1$. 
In particular, $\pi_1(\mathcal{R})\cong\ker \mathcal{T}=\{0\}$ is trivial.
Therefore $\mathcal{R}$ is biholomorphically equivalent to its simply connected universal cover $\mathbb{C}$, and the unbranched projection $\mathbf{q}_t: \mathcal{R}\rightarrow\mathbb{C}$ of \eqref{RC} is biholomorphic as well.
The composition $\mathbf{q}_w\circ\mathbf{q}_t^{-1}$ of \eqref{HRCd} and \eqref{RC} provides the unique solution $\mathbb{C}\ni t\mapsto w(t)\in\widehat{\mathbb{C}}_2$ of ODE \eqref{ODEw}, for initial condition $w=\infty$ at time $t=0$.

For all degrees $d\geq 2$, of course, these results confirm and extend the PDE result of theorem \ref{thmheat} in the simplest homogeneous ODE case \eqref{ODEw}.
Indeed, any real-time heteroclinic orbit $w=\Gamma:W_-\leadsto W_+$ between any two hyperbolic equilibria $e_j$ has to be contained in $\widehat{\mathbb{C}}_d$\,.
In case $d>2$, let $t^*=r_0+\mi s^*$ denote any branch point of the branched covering $\mathbf{q}_t:\mathcal{R}\rightarrow\mathbb{C}$ such that $\lvert s^*\rvert>0$ is minimal, for that choice of $r_0$\,. 
In the unbranched case $d=2$, we simply pick $t^*$ from $\mathbf{q}_t(\mathbf{p}_w^{-1}(\infty))$.
Then $s\mapsto \Gamma(r_0+\mi s)$ starts at $\Gamma(r_0)$ and blows up (or down) to the ramification point $(w,t)=(\infty,t^*)$, at imaginary Schrödinger time $s=s^*$.

Only in the quadratic Masuda case $d=2$, however, the isolated blow-up at $t=t^*$ can be circumnavigated in complex time.
For $d>2$, unique continuation by circumnavigation fails due to branching at $w=\infty$. 
In fact, we will see in the next sections, how two basic options for continuation will lead to two different target equilibria $W_+$\,, after circumnavigation, once real time flow is resumed.

\subsubsection{Poincaré compactification}\label{Infinity}

Motivated by our PDE paradigm, we now emphasize \emph{real time} $r$ for our complex-valued solutions $r\mapsto w(r+\mi s)$.
See section\ref{Cyclotomic} for a detailed discussion of the cyclotomic case $f(w)=w^d-1$.
As a remedy for the $d-1$ trajectories of \eqref{sovz} which cross at $z=1/w=0$, we replace stereographic projection by Poincaré compactification.

Specifically, we introduce polar coordinates $z=\rho\exp(\mi\alpha)$. 
Then \eqref{ODEw} reads
\begin{align}
\label{rho}
\dot{\rho}\ &=\rho\big(-\cos((d-1)\alpha)+\ldots \big)\,;   \\
\label{alpha}
\dot{\alpha}\ &=\phantom{\rho\big(-}\ \sin((d-1)\alpha)+\ldots\,.
\end{align}
Here we have vested \eqref{ODEw} with a local Euler multiplier $\rho^{d-1}$, and we have omitted higher order terms in $\rho$.
This ``blow up'' at $z=0$, in the sense of singularity theory, replaces $w=\infty$ by the circle $\rho=0,\ \alpha\in\mathbb{S}^1$.
In other words, local polar coordinates for $z=1/w$ essentially compactify the plane $w\in\mathbb{C}=\mathbb{R}^2$ to the closed unit disc $\mathbf{D}$, instead of the Riemann sphere $\widehat{\mathbb{C}}$.
Indeed the boundary circle $\mathbb{S}^1$ of $\mathbf{D}$ corresponds to the circle $\alpha\in\mathbb{S}^1$ at $\rho=0$.
See figure \ref{fig3} (a),(c) for this \emph{Poincaré compactification} of $w\in\mathbb{C}$.

On the invariant boundary circle $\rho=0$, a total of $2(d-1)$ hyperbolic saddle equilibria $\alpha_k=\pi k/(d-1)\in\mathbb{S}^1$ appear, for $k\,\textrm{mod}\,2(d-1)$.
In figure \ref{fig3} (a) and (c), we label the saddles at $\alpha_k$ by $\mathbf{k}$\,.
Within the boundary circle, they are alternatingly unstable, at even $\mathbf{k}$, and stable, at odd $\mathbf{k}$.
Their corresponding stable and unstable manifold counterparts, arriving from and emanating into the open disk $\mathbb{D}$ of $|w|=|1/z|<\infty$, are marked red and blue, respectively.
In \eqref{ODEw}, they mark red blow-up and blue blow-down of $w(t)$, in finite real original time $t=r$.
They also delimit the $2(d-1)$ local hyperbolic sectors of the equilibrium $z=0$, according to the planar classification of degenerate saddles by Poincaré. 
See section VII.9 in \cite{Hartman} and figure \ref{fig3} (b),(d).

Let us compare Poincaré compactification \eqref{rho}, \eqref{alpha} with the branching at $t=0$ associated to the ramification point $(w,t)=(\infty,0)$ of the branched covering projection $\mathbf{q}_t$\,; see section \ref{Riemann}.
The $d-1$ branches converging towards the boundary circle (red) simply are the $d-1$ lifts, by the $t$-projection $\mathbf{q}_t$\,, of the negative real axis $t=r<0$ to the sheet of the ramification point in the Riemann surface $\mathcal{R}$.
Indeed, the $d-1$ lifted trajectories $(w(r),r)\in\mathcal{R}$ given by $z(r)=1/w(r)$ in \eqref{sovz} parametrize the red blow-up trajectories in figure \ref{fig3}, which arrive at the boundary circle $z=0$ for $r\nearrow 0$.
Similarly, the $(d-1)$ blue blow-down branches emanating from the boundary circle correspond to real $t=r>0$, for $r\searrow 0$.
The alternating blue and red branches are $d-1$ real analytic curves in the complex $z$-plane which intersect, at $z=0$, under equal asymptotic angles $\pi/(d-1)$; see \eqref{alpha} and figure \ref{fig3}, again.

\textbf{Glossary}. We provide a brief glossary of terms and color codings concerning the real-time dynamics of Poincaré compactifications, for perusal in this paper.
See figure \ref{fig3} (a),(c) for an illustration of all terms.
Dynamics reside on the closed unit disk $\mathbf{D}$ with \emph{interior} $\mathbb{D}$ and invariant \emph{boundary circle} $\mathbb{S}^1$.
The $d$ \emph{interior equilibria} $e_j\in\mathbb{D},\ 0\leq j<d$, possess nonzero complex linearizations $f'(e_j)=1/\eta_j$\,. 
They are called (blue) \emph{sinks}, (purple) \emph{Hopf centers}, and (red) \emph{sources} in case $\Re f'(e_j)$ is strictly negative, zero, or strictly positive, respectively.
Heteroclinic orbits $e_j\leadsto e_k$ from sources $e_j$ to sinks $e_k$ in $\mathbb{D}$ (black) foliate open regions and are called \emph{interior source/sink heteroclinics}.

The $2(d-1)$ \emph{boundary equilibria} $\mathbf{k}\in\mathbb{S}^1$ are equidistantly spaced, at angles $\alpha_k=\pi k/{(d-1)}\in\mathbb{S}^1,\ 0\leq k < 2(d-1)$.
They enumerate all hyperbolic saddles, of alternating stability and instability inside $\mathbb{S}^1$.
For even $\mathbf{k}$, the unique (red) half branches of their stable manifolds arriving from the interior $\mathbb{D}$ are called  \emph{blow-up orbits}.
Their unique (blue) unstable manifold counterparts for odd $\mathbf{k}$, emanating into $\mathbb{D}$, are called \emph{blow-down orbits}.
In the exceptional cases when two such manifolds coincide in $\mathbb{D}$, they form a (purple) saddle-saddle heteroclinic orbit between boundary saddles, which we call an \emph{interior saddle-connection} or \emph{blow-down-up} orbit.
In all other cases, blue blow-down orbits are heteroclinic orbits towards sinks, and red blow-up orbits are heteroclinic orbits emanating from sources.
All non-equilibrium orbits within the invariant boundary circle $\mathbb{S}^1$ are saddle-saddle heteroclinic orbits between adjacent boundary saddles, which we call \emph{boundary saddle-connections}.

\subsubsection{Periodic orbits and saddle connections}\label{PerHet}

Following results in \cite{fiedlerShil}, section 4, we explore some implications of the real degeneracy
\begin{equation}
\label{sumsetabif}
 \sum_{j\in J}\,\Re\,\eta_j\,=\,0,
\end{equation}
with some given nonempty index set $\emptyset\neq J\subsetneq \{0,\ldots,d-1\}$.
Compare with our complex nondegeneracy assumption \eqref{sumsetaJ} for the realization of prescribed inverse derivatives $\eta_j=1/f'(e_j)$.
Due to \eqref{sumseta}, degeneracy \eqref{sumsetabif} holds for the index set $J$, if and only if it holds for the $J$-complement $\emptyset\neq J^c\subsetneq \{0,\ldots,d-1\}$.

Suppose there exists an interior saddle-connection $\gamma:\mathbf{A}\leadsto\mathbf{B}$, in the open disk $\mathbb{D}$ of the Poincaré compactification, between boundary saddle equilibria $\mathbf{A}, \mathbf{B} \in \mathbb{S}^1$.
Let $J$ enumerate the interior source/sink equilibria $e_j$ which are on the opposite side of $e_0$ in $\mathbb{D}\setminus\gamma$. 
Then $J, J^c$ are nonempty, and degeneracy \eqref{sumsetabif} holds.
We formulate a partial converse in section \ref{Conn}.

Suppose $J=\{j\}$ is a singleton. 
Then $f'(e_j)=1/\eta_j\neq 0$ is purely imaginary, and the equilibrium $e_j$ is a Hopf center.
E.g. by Poincaré linearization in imaginary time, a neighborhood of $e_j$ is foliated by a synchronous family of nested iso-periodic orbits of minimal period $p=2\pi\mi\eta_j$\,.
The boundary of that family in the open unit disk $\mathbb{D}$ is the interior saddle-connection associated to real degeneracy \eqref{sumsetabif}.

Conversely, any nonstationary periodic orbit in $\mathbb{D}$ encircles a unique equilibrium $e_j$, which turns out to be a Hopf center.
In imaginary time, section \ref{ODE2} provides the simplest example $d=2$.
For $d=4$ in real time, see figure \ref{fig3} (c),(d).

Absence of real degeneracies, a generic property among polynomials $f$ of degree $d\geq 2$, therefore implies the absence of any Hopf centers, periodic orbits and interior saddle-connections.
In other words, the dynamics of the Poincaré compactification in $\mathbb{D}$ is then of Morse type.

\subsubsection{Phase portraits, connection graphs, and trees}\label{Conn}

In this section we classify the Poincaré compactifications \eqref{rho}, \eqref{alpha} of flows to \eqref{ODEw}.
See \cite{fiedlerShil}, sections 1.9 and 4 for details.
We consider degrees $d\geq2$. 
Contrary to \eqref{sumsetabif}, we now assume the following \emph{strong nondegeneracy condition}
\begin{equation}
\label{sumsetaRe}
 \sum_{j\in J}\,\Re\,\eta_j\,\neq\,0
\end{equation}
to hold, for all nonempty index sets $\emptyset\neq J\subseteq\{0,\ldots,d-1\}$\,.
In particular, all $d$ interior equilibria $e_j$ of \eqref{ODEw} are hyperbolic, i.e. $\Re\,f'(e_j)=\Re(\eta_j/|\eta_j|^2)\neq0$.
Recall that this makes $e_j$ a \emph{sink}, for $\Re\,f'(e_j)<0$, and a \emph{source} for $\Re\,f'(e_j)>0$.
Moreover, the flow in the open disk $\mathbb{D}$ is of Morse type; see section \ref{PerHet}.

We classify and count the compactified flows, in real time $t=r$, up to \emph{orientation preserving} $C^0$ \emph{orbit equivalence}.
In other words, two flows on the closed unit disk $\mathbf{D}$ are considered orbit equivalent, if there exists an orientation preserving disk homeomorphism $H: \mathbf{D}\rightarrow\mathbf{D}$ which maps real-time orbits to real-time orbits.
We do not require $H$ to conjugate the real-time flows, and we do allow $H$ to reverse the real time direction.
So, $H$ conjugates real-time orbits, as sets, but not necessarily their specific time parametrizations.
We call the equivalence classes \emph{phase portraits}.

We associate phase portraits of Poincaré flows to equivalence classes of planar undirected trees $\mathbf{t}$ with $d$ vertices and $d-1$ edges, say contained in the open unit disk $\mathbb{D}$.
Abstractly, and in our setting, the trees are known as the \emph{connection graphs} or \emph{Morse graphs} among the sources and sinks in $\mathbb{D}$.
For many other applications of this concept, see for example \cite{Conley, brfi88, brfi89, Mischaikow, firoSFB, firoFusco, Yorke,Yorke2} and the many references there. 
Vertices of the connection graph $\mathbf{t}$ are the sources and sinks, in our case.
Nonempty families of interior source/sink heteroclinic orbits $e_j\leadsto e_k$ are represented by a single edge, each.
Note absence of undirected cycles in the connection graph $\mathbf{t}$, because any sector between two edges attached to a source contains at least one blow-up orbit.
Since source and sink vertices alternate along $\mathbf{t}$, we may leave edges undirected, up to global time reversal.

An obvious planar embedding of $\mathbf{t}\subset\mathbb{D}$ is given by construction of the source/sink vertices and the choice of representative heteroclinic edges.
We consider two planar trees as \emph{equivalent}, if there exists an orientation preserving disk homeomorphism $H$ of $\mathbb{D}$ which acts as a graph isomorphism on the trees.
In other words, $H$ maps vertices to vertices, and edges to edges, preserving their adjacency relations and the left cyclic orderings of corresponding edges around each vertex.
The direction of edges may be reversed under $H$.
We do not require any distinguished vertices or edges to be mapped to each other, like any vertices marked as ``roots'' or otherwise labeled. 

A third view point are \emph{diagrams of} $d-1$ nonintersecting, undirected \emph{ chords} of the closed unit disk $\mathbf{D}$, up to rotation.
Here the chords are closed straight lines between their $2(d-1)$ end points on the disk boundary $\mathbb{S}^1$, spaced equidistantly at angles $\beta_k:=\pi (k+\tfrac{1}{2})/(d-1),\ 0\leq k<2(d-1)$.
We call two \emph{chord diagrams} equivalent, if they coincide up to a proper rotation of $\mathbf{D}$.
%The rotation is allowed to reverse the direction of chords.
Chord diagrams are also called ``noncrossing (single-armed) handshakes of $2(d-1)$ people on a round table'', or ``partitions of $2(d-1)$ elements into noncrossing blocks of size 2\,''; see \cite{oeis}.

In section \ref{PerHet}, we have observed how interior saddle-connections $\gamma:\mathbf{A}\leadsto \mathbf{B}$ imply a violation \eqref{sumsetabif} of real strong nondegeneracy condition \eqref{sumsetaRe}.
The index sets $j\in J,J^c$ of the violation mark sources/sinks $e_j$ on opposite sides of $\gamma$.
If the interior saddle-connection $\gamma$ is unique, then the two connection trees among sources and sinks in $J,J^c$, separately, turn out to be connected.
To formulate a converse, we now assume $J,\,J^c$ to be a nonempty complementary partition such that the real degeneracy \eqref{sumsetabif} holds for that pair of complementary index sets.
We also assume each of the two connection trees among sources and sinks in $J,J^c$, separately, to be connected.
Then there exists a unique interior saddle-connection $\gamma$, for the Poincaré compactification of \eqref{ODEw}. 

Under real strong nondegeneracy assumption \eqref{sumsetaRe}, the following two statements hold true, in terms of the three orientation preserving equivalence classes just described.
See \cite{fiedlerShil}, theorem 1.3.
\begin{enumerate}[(i)]
\item  The real-time phase portraits of the Poincaré compactifications correspond, one-to-one, to the above planar embeddings of their connection graphs $\mathbf{t}$.
\item Equivalently, they correspond, one-to-one, to certain chord diagrams of $d-1$ unlabeled, undirected, nonintersecting chords of the unit circle.
\end{enumerate}

For equivalence of planar trees and chord diagrams, as well as combinatorial enumerations, see \cite{oeis, Treecount1,Treecount2}. 
Up to orientation preserving equivalence, the number of unlabeled, unrooted, undirected, planar trees with $d$ vertices, alias chord diagrams of $d-1$ unlabeled, undirected, nonintersecting chords of the unit circle, is given by entry $A_{d-1}$ of sequence A002995 in the online encyclopedia of integer sequences \cite{oeis}.
The counts for $2\leq d\leq 16$ are
\begin{equation}
\label{Treecount}
1, 1, 2, 3, 6, 14, 34, 95, 280, 854, 2694, 8714, 28640, 95640, 323396, \ldots\ .
\end{equation}
An explicit expression for the general counts $A_{d-1}$ was derived by \cite{Treecount1}, theorem 2, as
\begin{equation}
\label{Treecount1}
\begin{aligned}
A_{d-1}\,=\,\tfrac{1}{2(d-1)d}\, \tbinom{2(d-1)}{d-1}\,&+\,\tfrac{1}{4(d-1)}\tbinom{d}{d/2}\,+\,\tfrac{1}{d-1}\,\phi(d-1) \,+\\
                  &+\,\tfrac{1}{2(d-1)}\,\sum_{k=2}^{d-2}\,\tbinom{2k}{k}\, \phi(\tfrac{d-1}{k}) \,.           
\end{aligned}
\end{equation}
Here $\phi$ denotes the Euler totient count of coprime elements, and the sum only runs over proper divisors $k$ of $d-1$.
For odd $d$, the second summand on the right is omitted.

Although each strongly nondegenerate orbit equivalence class, in the sense of assumption \eqref{sumsetaRe}, corresponds to one of those trees, it was not clear, at first, whether all trees are actually realized by Poincaré compactifications.
But, yes, they are.

\begin{thm}\label{Classthm} \emph{\cite{fiedlerShil}}
Unlabeled, unrooted, undirected, planar tree with $d$ vertices, alias chord diagrams of $d-1$ nonintersecting, undirected chords, are realized, one-to-one, by Poincaré compactifications of \eqref{ODEw}, for suitable polynomials $f$ of degree $d$.
\end{thm}

\subsubsection{Cyclotomic examples}\label{Cyclotomic}

Geometric analysis of the special cyclotomic cases
\begin{equation}
\label{cyclotomic}
\dot{w}=f(w)=w^d-1
\end{equation} 
for $d\geq 2$ illustrates how Masuda's paradigm is limited to the quadratic case $d=2$.
We refer to \cite{fiedlerShil}, section 1.8 for further details.
The equilibria are $w=e_j=\exp(2\pi\mi j/d)$ with $f'(e_j)=de_{-j}$ and $\eta_j=\tfrac{1}{d}\,e_j$\,, for $j\,\mathrm{mod}\,d$.
The complex nondegeneracy assumption \eqref{sumsetaJ} is violated, if and only if $d$ is nonprime.
This example demonstrates how nondegeneracy \eqref{sumsetaJ} may not be necessary, after all, for realization of reciprocal derivatives $\eta_j=1/f'(e_j)$ by complex polynomials $f$.

From \eqref{rho}, \eqref{alpha}, more dynamically, we recall invariance of the boundary circle $\rho=0,\ \alpha\in\mathbb{S}^1$, for the Poincaré compactification at $w=\infty$, with $2(d-1)$ alternating hyperbolic saddles $\mathbf{k}$\,.
Since the cyclotomic polynomial $f$ possesses real coefficients, complex conjugation provides reflection symmetry between the flows in the upper and lower half plane of $w=u+\mi v$.
Indeed $w(t)$ solves \eqref{ODEw}, if and only if $\overline{w}(t)$ does.
In particular, the horizontal real axis is invariant.
This determines the horizontal blow-up / blow-down orbits of the boundary equilibria $\mathbf{k}=\mathbf{0},\ \mathbf{d-1}$\,, and the real blow-up orbit, or orbits, of $e_0$\,.
For even $d$, the real heteroclinic orbit $e_0\leadsto e_{d/2}$ extends to, both, real blow-up $e_0\leadsto\mathbf{0}$ and to real blow-down $\mathbf{d-1}\leadsto e_{d/2}$ via $w=\infty$, in complex time.
For odd $d$, blow-up $e_0\leadsto\mathbf{0}$ and $e_0\leadsto\mathbf{d-1}$ occurs in both real directions.
By reflection symmetry, it remains to study the upper half plane $v>0$ of $w=u+\mi v$.

Consider the case $d=3$ of figure \ref{fig3} (a), (b), first.
Strong nondegeneracy \eqref{sumsetaRe} then implies the Morse structure of the compactified Poincaré dynamics in the open unit disk $\mathbb{D}$.
The Poincaré-Bendixson theorem \cite{Hartman} therefore identifies the sink $e_1$ as the only possible $\boldsymbol{\omega}$-limit  set of the blue blow-down separatrix of $\mathbf{1}$ in the upper half plane; see also section \ref{PerHet}.
All other trajectories $w$ in the open upper half plane are heteroclinic of type $e_0\leadsto e_1$.

In case $d=4$ of figure \ref{fig3} (c), as for any even $d$, we encounter time reversibility of \eqref{cyclotomic}.
The time reversor is horizontal reflection at the vertical imaginary axis.
For even $d$, indeed, $w(t)$ solves \eqref{cyclotomic}, whenever $-\overline{w}(-t)$ does.
In particular, the Hopf center equilibria $e_1, e_3 =\pm\mi$ are locally surrounded by periodic orbits, only, up to the purple interior saddle-connections $\mathbf{1}\leadsto\mathbf{2}$ and $\mathbf{5}\leadsto\mathbf{4}$.
%In fact, the unstable blow-down separatrix of $\mathbf{1}$ (usually blue, here purple) has to cross the imaginary axis $\Re w=0<\Im w$. 
%Indeed, the first quadrant $\{\Re w>0,\ \Im w>0\}$ does not contain equilibria and, therefore, cannot contain periodic orbits.
%Reflecting on the first crossing, say at time $t=0$, reversibility implies that the unstable blow-down separatrix has to coincide with the stable blow-up separatrix of $\mathbf{2}$, usually colored red.
%We have therefore colored the resulting interior saddle-connection $\mathbf{1}\leadsto\mathbf{2}$ purple.
%Similar arguments show that the resulting heteroclinic loop between $\mathbf{1}$ and $\mathbf{2}$ is filled with periodic orbits around the center $e_1$.
In (d), i.e. upon identification of the boundary circle $\rho=|z|=0,\ \alpha\in\mathbb{S}^1$ with $w=\infty$, the heteroclinic loops become counterrotating homoclinic to $z=1/w=0$.
Symmetrically, the homoclinic lobes remain foliated by iso-periodic families of synchronously counterrotating periodic orbits, as described in section \ref{PerHet}.
Indeed their minimal periods $\mp\pi/2$, in original time and prior to any rescaling, all coincide.
The periodic families are unbounded in $w=1/z$.
All remaining trajectories $w$ in the open upper half plane are again interior source/sink heteroclinic, of type $e_0\leadsto e_2$.

Figure \ref{fig3} also demonstrates the behavior of the global trajectories $w=\Gamma_{\pm\mi\eps}$ under slight perturbations $w(0)=2\pm\mi\eps,\ \eps\searrow 0$, of the real blow-up initial condition $w(0)=2$.
In the case $d=3$ of figure \ref{fig3} (a), the interior source/sink heteroclinic $w=\Gamma_{+\mi\eps}: e_0\leadsto e_1$ in the upper half plane is labeled by $+\mi\eps$ (black).
For small $\eps>0$, it closely follows the concatenated heteroclinic chain $e_0\leadsto \mathbf{0}\leadsto \mathbf{1}\leadsto e_1$. 
The interior heteroclinic $\Gamma_{-\mi\eps}: e_0\leadsto e_2$ in the lower half plane\,, in contrast, labeled $-\mi\eps$, closely follows the \emph{different} concatenation $e_0\leadsto \mathbf{0}\leadsto \mathbf{3}\leadsto e_2$.
In sharp contrast to the quadratic case $d=2$, therefore, the red blow-up cannot be circumnavigated in complex time, consistently. 

The polar view of (b), centered at $w=\infty$ alias $z=1/w=0$, conflates the boundary circle $\rho=|z|=0,\ \alpha\in\mathbb{S}^1$ to a single point.
This shortens the black trajectories $\Gamma_{\pm\mi\eps}$ to perturbations of the two \emph{distinct concatenations} $e_0\leadsto \infty\leadsto e_1$ and $e_0\leadsto \infty\leadsto e_2$\,.
In the limit $\eps\searrow 0$, this amounts to a shared red blow-up  $e_0\leadsto \infty$, followed by two different subsequent blue blow-down orbits $\infty\leadsto e_1$ and $\infty\leadsto e_2$\,.
Again, this is in sharp contrast to the case $d=2$ of figure \ref{fig2}, where both approximations to the cyan blow-up orbit $1\leadsto\infty$ continue along \emph{the same} cyan blow-down $\infty\leadsto-1$.
Since this happens in finite original time $w(t)$, the Masuda paradigm of complex near-recovery for near-homogeneous PDE solutions, after real blow-up, turns out to be a peculiarity of the quadratic nonlinearity $d=2$.

The first nonprime degree $d=4$ of figure \ref{fig3}, (d) leads to heteroclinic perturbations $w=\Gamma_{\pm\mi\eps}: e_0\leadsto e_2$ which, at least formally, seem to limit onto identical concatenations $e_0\leadsto \infty\leadsto\infty\leadsto e_2$, for $\eps\searrow 0$.
Both limits share the initial red blow-up part $e_0\leadsto\infty$ and the terminal blue blow-down part $\infty\leadsto e_2$\,.
The two purple interior saddle-connections, i.e.~the homoclinic blow-down-up orbits $\infty\leadsto\infty$, however, remain distinct.
The limiting homoclinic loop part $+\mi\eps$ of $\Gamma_{+\mi\eps}$ contains the Hopf center $e_1$ in its interior, whereas $-\mi\eps$ of $\Gamma_{-\mi\eps}$ surrounds the Hopf center $e_3$\,, along with their mandatory homoclinic lobes, respectively filled by nested families of counter-rotating, synchronously iso-periodic orbits with constant minimal periods $\mp\pi/2$.
Again, these perturbations run against the quadratic Masuda paradigm.
The observed discrepancies between $\Gamma_{\pm\mi\eps}$ run deeper than ``mere technical'' PDE difficulties like the absence of a heat semiflow in reverse time $r=\Re t<0$.
In fact they already originate from non-unique complex continuation of homogeneous real-time blow-up itself.

For general $d>4$, we may still extend the real blow-up orbit of $e_0=1$, by boundary concatenation to the distinct blow-down orbits emanating from the neighboring boundary saddles $\mathbf{1}$ and $\mathbf{2d-3}$ towards distinct conjugate complex sinks.
Perturbations $\Gamma_{\pm\mi\eps}$ will therefore remain distinct, for arbitrarily small $\eps>0$.
In contrast to the quadratic Masuda paradigm, they do not converge to a single continuation beyond blow-up.

\subsection{Closure and boundaries of fastest unstable manifolds}\label{Boundary}

Let $\partial W^u := \mathrm{clos}\, W^u \setminus W^u$ denote the topological boundary of unstable manifolds $W^u$.
Before we return to complex time $t=r+\mi s$, let us recall the case of dissipative gradient-like flows in real time $t=r$ for a moment.
For some of the intricacies, like Alexander horned spheres which already arise in three-dimensional settings, see the appendix by Laudenbach in \cite{Laudenbach}, and the references there.
Generalizations 
\begin{equation}
\label{heatg}
u_r = u_{xx}+f(x,u,u_x)
\end{equation}
of the complex quadratic heat equation \eqref{heat} have been studied in \cite{WolfrumCorr}, for real $u$ rather than complex $w=u+\mi v$.
Under our usual Neumann boundary conditions, the dynamics remains gradient-like; see \cite{LappicyAttractors, LappicyBlowup} and the references there.
Growth assumptions on $f$ ensure dissipativeness and the existence of a global attractor which consists of equilibria and their heteroclinic orbits, only.
For hyperbolic real equilibria $u=U(x)$, transversality of stable and unstable manifolds holds automatically.
In fact, the decomposition of the global attractor into the disjoint union of unstable manifolds $W^u(U)$ turns out to be a \emph{regular cell complex}.
The cell boundaries $\partial W^u$ themselves possess \emph{signed hemisphere decompositions}.
Based on heteroclinic orbits $U=U_n\leadsto U_0$ in $W^u(U_m)$, that hemisphere decomposition is governed by cascades $U_n\leadsto U_{n-1} \leadsto \ldots U_m$ of intermediate heteroclinic orbits.
Here the Morse indices $i(U_k)=k$ in the cascade drop by 1 at each heteroclinic step.
See \cite{WolfrumCorr} for further details.

Even for the complex quadratic heat equation \eqref{heat}, alas, the boundaries $\partial W^u$ of unstable manifolds are not as easily described.
To simplify matters, we discuss heteroclinic orbits $\Gamma:W_-\leadsto W_+$ emanating from the fastest unstable manifold $W_-^{uu}$ of dimension $d=1$ again; see section \ref{Fastest}. 
The spatially homogeneous source equilibrium $W_-=W_\infty$ leads to spatially homogeneous fastest $\Gamma$, as discussed in sections \ref{ODE2}, \ref{ODEd} above.
For conjectural numerical explorations of blow-up at isolated branching singularities in the technically similar, purely quadratic setting $\lambda=0$ of 1-periodic boundary conditions, see \cite{COS,Fasondini23}.

Locally near any non-homogeneous source $W_-=W_n$\,, and therefore globally, the complex time extension of the fastest heteroclinic orbit $\Gamma:W_n\leadsto W_0$ in $W_-^{uu}$ still parametrizes $W_-^{uu}\setminus\{W_-\}$.
By Poincaré linearization \eqref{P}, in fact, $z:=\exp(\mu_0 t)$ parametrizes $W_-^{uu}\setminus\{W_n\}$, with minimal period $p=2\pi/\mu_0$ in imaginary time.
Here $t=r+\mi s \in\mathbb{C}$, with real part $r$ fixed sufficiently negative.

At the spatially homogeneous, asymptotically stable target $W_+=W_0$ of $\Gamma(r)$, matters look quite different. 
Comparing with \eqref{Gammaperiod}, we first notice how $\Gamma(r+\mi s)$ fails to be $p$-periodic, in the imaginary time direction $s$ and for $r$ fixed large positive.
This failure is due to blow-up at intermediate real parts $r$, and will occur for all but a discrete set of parameters $\lambda$.
In particular, the semiflows $\Phi^r$ and $\Phi^{\mi p}$ fail to commute on $\Gamma$: our previous statement \eqref{shift} concerning commutativity in parallelograms was due to our indirect, and ultimately false, assumption of a complex entire heteroclinic orbits $\Gamma$.

We fill in a few details on the absence of imaginary periods $\mi p$ of $\Gamma$, near $W_+$.
Rescaling \eqref{nscale} allows us to restrict to the case $n=1$ of $W:=W_-=W_1\,,\ \lambda > \lambda_{1,0}=\tfrac{2}{3}\pi^4$; see \eqref{biflambda}.
In real time $t=r\nearrow+\infty$, the asymptotics of $\Gamma\rightarrow W_+=W_0$ is well-understood; see for example \cite{brfi86}.
With purely trigonometric eigenfunctions $c_k=\cos(2\pi k x)$, at eigenvalues
\begin{equation}
\label{spec0}
\mu_{+,k}=-2\sqrt{6\lambda}-(2\pi k)^2,
\end{equation}
and with the homogeneous target $W_+=W_0=-\sqrt{\lambda/6}$ shifted to zero, we obtain
\begin{equation}
\label{Gamma+}
\Gamma(r)=\exp(\mu_{+,0}\,r)+o
\end{equation}
in $H^1$-norm.
We have also shifted real time $r$, appropriately.
For $r\nearrow+\infty$, the error term ``$o$'' decays of higher order than any negative exponential rate between $\mu_{+,0}<0$ and $\mu_{+,1}<\mu_{+,0}$\,.
In imaginary time direction $s$, we may also compare the Schrödinger flow $S^s=\Phi^{\mi s}$  with the linearized flow at the target $W_+$\,. 
Uniformly for bounded $s$, we obtain the same estimate for the complex continuation
\begin{equation}
\label{Gammacont}
\Gamma(r+\mi s)=\exp(\mu_{+,0}\,(r+\mi s))+o\,.
\end{equation}

On the other hand, we may invoke $p$-periodicity of $s\mapsto\Gamma(r+\mi s)$ at the source $W_-$\,, as established in \eqref{Gammaperiod} for large negative $r$.
By analytic continuation in $r$, at fixed $s=p$, periodicity clashes with the $r$-asymptotics \eqref{Gammacont} as follows.
Suppose $\Phi^r$ and $\Phi^{\mi\,p}$ commute on $\Gamma$, for large $r>0$.
Then  $\exp(\mu_{+,0}\,\mi p)=1$.
This implies the global resonance
\begin{equation}
\label{resglobal}
\mu_{+,0}+m\,\mu_0=0,
\end{equation}
for some integer $m>0$.
But at the pitchfork bifurcation $\lambda=\lambda_{1,0}$ of $W=W_n$\,, for $n=1$, we have already obtained the expansions 
\begin{align}
\label{mu0h}
\mu_0   &= \ \ 4\pi^2(1+264\, h^2 +\ldots   ) >0   \\
\label{mu+0h}
\mu_{+,0}&=-4\pi^2(1+120\, h^2+\ldots) <0
\end{align}
along the branch parameter $h$.
For \eqref{mu0h}, see \eqref{muphi0}, \eqref{muphi1}, \eqref{muphi2}, with $n=1,\ k=0$.
For \eqref{mu+0h}, see \eqref{spec0} and \eqref{lambdanh}.
This contradicts global resonance \eqref{resglobal}: locally at pitchfork bifurcation $\lambda=\lambda_{1,0}=\tfrac{2}{3}\pi^4$, and hence everywhere, up to a discrete set of resonant parameters $\lambda$.
Therefore $\Phi^r$ and $\Phi^{\mi\,p}$ actually fail to commute on $\Gamma$ near $W_+$\,.

For rationally independent eigenvalues $\mu_{+,0}\,,\ \mu_{+,1}$, linearization at the target $W_+$ hints at yet another possible complication concerning the closure, and the boundary $\partial W_-^{uu}$, of the fastest unstable manifold $W_-^{uu}$ generated by the heteroclinic orbit $\Gamma$.
When both components are present, the resulting quasi-periodicity (or almost-periodicity) in imaginary time $\mi s$, asymptotically for fixed large $r>0$, might cause dense orbits $\Gamma(r+\mi s)$ in tori.
In the more specific second order ODE setting \eqref{twmu+}, this problem will return below in the guise of traveling wave solutions in imaginary time.

Above, we have considered the fastest unstable manifold at $W_-$\,, to fix a unique heteroclinic orbit $\Gamma:W_-\leadsto W_+$\,.
This has selected a unique complex one-dimensional object $t\mapsto\Gamma(t)$ for maximal analytic continuation.
Generated by other real analytic heteroclinic orbits $r\mapsto\Gamma(r)$, many more candidates are left to explore.
For example, heteroclinic orbits between equilibria of adjacent Morse indices, $i(W_+)=i(W_-)-1$, are known to be unique in the real case \cite{brfi88}.
More generally, Wolfrum has selected unique one-dimensional $\Gamma(r)$, characterized by constant zero number $z(\Gamma(r)-W_\pm)=z(W_--W_+)$, whenever there exists any real heteroclinic orbit $W_-\leadsto W_+$\,; see \cite{WolfrumLemma, WolfrumCorr}.
%Above, we have only addressed the very special case $z=0$.
Even under spectral non-resonance at $W_-$\,, Poincaré linearization \eqref{P} will then only assert $s\mapsto\Gamma(r+\mi s)$ is quasi-periodic, say for $r\leq 0$, rather than periodic.
The global complex geometry of such $\Gamma$, and their topological boundary, therefore remains unexplored open territory.

\subsection{Schrödinger and parabolic singularities}\label{Sing}

We discuss some of the previous literature on parabolic and Schrödinger blow-up, in finite complex time.

The updated standard mathematical PDE monograph on parabolic blow-up in \emph{real time} is \cite{Quittner}.
Topics include multi-dimensional $x$, scalar nonlinearities $f$  which combine $u^p$ and $(\nabla u)^q$, and systems of PDEs.
In addition to the standard PDE arsenal of weak versus strong solutions and Sobolev embeddings, methods are basically a combination of three ingredients: the variational structure of energy and Lyapunov functionals, comparison and maximum principles, and, at least locally and for purely homogeneous polynomial nonlinearities like $\lambda=0$, self-similar scalings \eqref{nscale} by arbitrary, not necessarily integer $0<n\in\mathbb{R}$.
Marek Fila figures prominently in that monograph with 25 of his papers.

Consider the purely quadratic real parabolic system \eqref{heat} with $\lambda=0$, for the real and imaginary parts of $w=u+\mi v$.
Restricted to real time, but for real $m$-dimensional $x$, this case has been addressed in \cite{Yanagida} with complete mathematical rigor.
The variational structure is lost, except in the invariant real subspace $v\equiv 0$.
Benefitting from equal diffusion in the components $u$ and $v$, however, comparison principles remain applicable.
Together with self-similarity, and for certain initial conditions, this provides dichotomies between global existence versus blow-up in finite real time.

A very accessible survey on real-time blow-up phenomena of complex-valued ODEs and PDEs, in a similar spirit, is \cite{Kevrekidis}. 
It is aimed at a broad audience, with a strong emphasis on applications and numerical simulations.

What is usually ``forgotten'' in the above approaches, or at least ignored, are the subtle relations among real-time solutions $r\mapsto w(r+\mi s)$, for various fixed $s$.
These relations are induced and encoded by their contingency \eqref{flow} under commuting local flows, or semiflows, in real and imaginary time $r$ and $\mi s$.
Beyond the kinship between the quadratic heat and Schrödinger equations \eqref{heat} and \eqref{psi} which we have explored above, this is a structure worthwhile of mathematical study, as such, even when a direct interpretation of imaginary time remains elusive ``in nature'', for now.

In the simplest, and not quite typical, case of homogeneous solutions, the orange periodic circles of figure \ref{fig2}, in imaginary time, have illustrated such iso-periodic contingency which globally synchronizes the orthogonal foliation by blue real-time heteroclinic orbits. 
And vice versa. 
See also our discussion in sections \ref{ODEd}, \ref{Boundary} above and the synchronicity of nested iso-periodic orbits in section \ref{PerHet} and figure \ref{fig3} (c),(d), due to the residue theorem.

Forty years ago, Kyûya Masuda has started a systematic study of blow-up in \emph{complex time} with his pioneering work on the purely quadratic heat equation \eqref{heat} at parameter $\lambda=0$, and in any space dimension.
Under Neumann boundary conditions, elliptic maximum principles imply that $u\equiv 0$ is the only real equilibrium.
For almost homogeneous real initial conditions $0\leq w_0\not\equiv 0$, Masuda was able to circumvent blow-up, at real time $0<r=r^*(w_0)<+\infty$, via a sectorial detour in complex time $t=r+\mi s$.
See \cite{Masuda2} for proofs, and the earlier announcement \cite{Masuda1}.
Notably, he also cautioned that the detours via positive and negative imaginary parts $s$ agree in their real-time overlap after blow-up, if \emph{and only if} $u_0$ is spatially homogeneous.

More than three decades after this amazing discovery, \cite{COS} pointedly remarked:
\begin{displayquote}
\emph{
\emph{Masuda’s} pioneering work does not seem to have a follower. \ldots
His paper seems to be a good example of a paper of high originality with few
citations.
 It seems to us that, although there are numerous papers studying blow-up of solutions of nonlinear parabolic equations, few investigate singularities in the complex $t$-plane.}
\end{displayquote}
For some further references in the context of the notoriously difficult Navier-Stokes and Euler PDEs, as well as complex $x$ and vorticity ODEs, see also \cite{COS, LiSinai, Fasondini24}.
The most diligent monograph on parabolic blow-up in real time, \cite{Quittner}, kindly mentions \cite{Masuda2} once, in passing, for his complex circumvention of ``complete'' real-time blow-up; see the short remark 27.8(f).

% written by Hannes
In one space dimension, Stuke has addressed the Masuda setting under Dirichlet boundary conditions \cite{Stukediss, Stukearxiv}.
Even though Stuke's setting is quite similar to Masuda's, the techniques are quite different. 
Masuda exploited the homogenous ODE solution.
In the one-dimensional fastest unstable manifolds, Stuke detected how heteroclinic orbits, blow-up, and center manifolds in real time are related via complex time. 
In the case of the semilinear heat equation, a trapping argument similar to \cite{Yanagida} showed that blow-up must occur on the real time axis. 
This led to the paradigm of the present paper.
For the quadratic heat equation, the results also added to Masuda's attempts at a circumnavigation of real-time blow-up, via the complex time-domain.
Continuation was pursued, separately, in the upper complex half plane $s> 0$, and in the lower half plane
$s< 0$.
Return to the real time axis was achieved, only, after a finite slit interval of real latency times $r$ trailing $r=r^*$. 
The fact that the continuations do not match after blow-up was seen as caused by the ratio of asymptotic eigenvalues at the heteroclinic equilibria $W_\pm$\,.
Attempts to reduce the trailing slit interval to the singleton $\{r^*\}$, i.e. to identify blow-up as an isolated singularity in complex time as conjectured by \cite{COS}, have failed so far, even on the rather formal level of matched asymptotic expansions \cite{Fasondini23}.

Jonathan Jaquette and co-authors mostly addressed the purely quadratic case $\lambda=0$ of the Schrödinger equation \eqref{psi}.
Spatial periodicity $x\in\mathbb{T}^m$ is imposed.
By reflection through the boundary, Neumann boundary conditions can be subsumed.

The quadratic Schrödinger equation on $x\in\mathbb{T}^m$ is integrated in \cite{Jaquetteqp}, by explicit recursive ODEs for the spatial Fourier-coefficients.
The one-dimensional case \eqref{psi} of $m=1$ and spatial period 1  is studied in detail, for spatially monochromatic initial conditions $\psi_0(x)=a \exp(2\pi\mi x)$.
Solutions $\psi=\psi(s,x)$ of \eqref{psi} are shown to be explicitly time-periodic, of fixed minimal period $p=1/(2\pi)$, alias frequency $(2\pi)^2$, for $|a|\leq 2\pi^2$.
For initial amplitudes $|a|\geq 4\pi^2$, in contrast, blow-up occurs in finite ``real'' Schrödinger time $s$.
In higher dimensions $m>1$, small solutions are shown to be quasi-periodic in time $s$.
To avoid small divisors, initial conditions $\psi_0$ are assumed to be supported on strictly positive spatial Fourier modes.
The fundamental frequencies $\Omega_j=\omega_j^2$ in time $s$ originate from the spatial frequencies $\omega_j$\,, alias spatial periods $q_j=2\pi/\omega_j$\ of $\mathbb{T}^m$, imposed in the spatial directions $x_j$\,, for $j=1,\ldots,m$.
From our abstract point of view, these results fit well with an elusive partial Poincaré linearization of \eqref{heat}, within the infinite-dimensional strong stable manifold $W^{ss}$ of the trivial equilibrium $W\equiv 0$ at $\lambda=0$.

As an alternative to the above Fourier analysis, an approach based on rigorous numerics is pursued in the trilogy \cite{JaquetteHet,JaquetteStuke,JaquetteMasuda}.
We comment in chronological order.

Rigorous numerics, in the spirit of previous analysis by \cite{Masuda1,Masuda2,Stukediss,Stukearxiv}, has been provided in \cite{JaquetteMasuda}.
The setting is the purely quadratic heat equation \eqref{heat}, $\lambda=0$, under 1-periodic boundary conditions in $m=1$ spatial dimensions.
Specifically, they follow \cite{COS} for the large real initial condition $w_0=50(1-\cos(2\pi x))$.
Corroborating \cite{COS}, they rigorously determine a box, in complex time, where a complex singularity of Masuda type has to occur.
They also consider the inclined Schrödinger-like equation \eqref{wtheta}, i.e. \eqref{heat} along rays $t=r\exp(\mi\theta),\ r\geq 0$.
Under the same initial condition $w_0$\,, in contrast, they prove global existence, and convergence $w(r)\rightarrow 0$ for $r\nearrow\infty$, at fixed angles $\theta \in\{1,2,3,4\}\cdot\pi/12$.
A second conjecture in \cite{COS} aims at such a result for $\theta=\pi/2$, i.e. for the Schrödinger case \eqref{psi}, and any initial condition $\psi_0\in L^2$.
The counterexamples of our theorem \ref{thmpsi}, of course, require $\lambda>0$ rather than $\lambda=0$.

The purely quadratic Schrödinger case \eqref{psi}, $\lambda=0$, is also studied in \cite{JaquetteHet}.
We comment, specifically, on the nontrivial complex equilibria $u_1^i$ and $u_1^{ii}$ of \cite{JaquetteHet}, Theorem 1.7, figure 1, in line with their open question 6.4.
Their result is also reproduced in \cite{JaquetteStuke}.
Both equilibria can in fact be obtained analytically, along the completely classical lines of our section \ref{wp}.
Consider $g_2=0$ in \eqref{lambdan}, with discriminant $\Delta=-27 g_3^2$ and Klein-invariant $J=0$.
The corresponding hexagonal period lattice $\Lambda$ in \eqref{defwp} is characterized uniquely by the modular parameter $\tau:=\exp(\pi\mi/3)$, alias purely imaginary $h=\exp(\pi\mi\tau)=\mi\exp(-\tfrac{\pi}{2}\sqrt{3})$ in \eqref{Wnh}, \eqref{eta}.
Under Neumann boundary conditions at minimal spatial half-period $1/2$, we obtain the equilibrium 
\begin{equation}
\label{ui}
u_1^i(x)= 6\,W_1(x):=-6\,\wp(x+\tfrac{1}{2}\tau).
\end{equation}
Families $W_{\pm n}$ arise, because the parameter $\lambda=0$ is not affected by the self-similar rescalings \eqref{nscale}.

Seeking complex equilibrium solutions for $\lambda=0$, we can utilize the more general invariance of the equilibrium equation \eqref{ODEW}, i.e. of $W_{xx}+6W^2=0$, under any complex scaling factor $\sigma\in\mathbb{C}$. 
This generates equilibrium families
\begin{equation}
\label{sigma}
W^\sigma(x):=\sigma^2\, W(\sigma x + z_0),
\end{equation}
for suitably chosen $\sigma$ and $z_0$\,.
The equilibrium $u_1^{ii}(x)$ of \cite{JaquetteHet} is an example.
Explicitly,
$u_1^{ii}(x)=-6\, \wp^\sigma(x)$ along the line segment $z=\sigma x + z_0,\ 0\leq x\leq\tfrac{1}{2}$, with rescaling factor $\sigma=\sqrt{3}\exp(\pi\mi/6)$ and offset $z_0=\tfrac{1}{2}\tau$.
In other words,
\begin{equation}
\label{uii}
u_1^{ii}(x)=-6\sigma^2\, \wp(\sigma x+\tfrac{1}{2}\tau).
\end{equation}
More generally, proper complex scaling produces complex Neumann equilibria $W^\sigma$ from $\wp$, evaluated on line segments between points $z_0$ and $z_0+\tfrac{1}{2}\sigma$ of the half-period lattice $\tfrac{1}{2}\Lambda$. 
The segments have to avoid the singularities of $\wp$ on $\Lambda$ itself, of course.
The solution $u_1^{ii}(x)$, for example, is based on the inclined segment from $z_0=\tfrac{1}{2}\tau$ to $\tau+\tfrac{1}{2}$\,.
The solution $u_1^{i}(x)$ evaluates $x$ horizontally, from $z_0=\tfrac{1}{2}\tau$ to $\tfrac{1}{2}\tau+\tfrac{1}{2}$\,.
Explicit Fourier expansions \eqref{Wnh}, \eqref{eta} of $\pm\wp$ apply.
It is an elementary, but worthwhile, task to explicitly enumerate all possibilities on the hexagonal lattice.

Under periodic boundary conditions, and even before invoking any rescalings $\sigma$, we obtain a 2-torus of primary complex equilibria $W(x)=-\wp(x+z_0)$.
Only shifts by real $z_0\in\mathbb{T}^2=\mathbb{C}/\Lambda$ correspond to spatial shifts of $x\in\mathbb{T}^1=\mathbb{R}/\mathbb{Z}$.
The exceptional solution for real $z_0$ is spatially quadratically singular at $x=0$, of course.
Subsequent complex scalings $\sigma$ also apply.

Remarkably, computer-assisted proofs establish the existence of heteroclinic PDE orbits between $W\equiv 0$ and the two nontrivial equilibria $W=u_1^i\,,\,u_1^{ii}$\,.
See \cite{JaquetteHet}, theorem 1.9 and, for $u_1^i$, also figure 2 there.
The heteroclinic orbits then run in either direction, due to a reversibility \eqref{defR} akin to \eqref{reversibility}. 
In addition, their theorem 1.5 presents an open set of small homoclinic PDE solutions to $W_\infty:= 0$, with exponentially decaying Fourier modes.

Since the equilibria $u_1^i\,,\,u_1^{ii}$ are not real, reversibility \eqref{reversibility} is adapted here as follows. We note that $\psi(s,x)$ solves the quadratic Schrödinger equation \eqref{psi}, if and only if \begin{equation}
\label{defR}
(R\psi)(s,x):= \overline\psi(-s,\tfrac{1}{2}-x)
\end{equation}
does.
As always, the parameter $\lambda=\bar\lambda$ is assumed to be real, here.
For the purely quadratic case $\lambda=0$ of the Schrödinger equation \eqref{psi}, we note $R$-invariance $R\,W=W$, for any equilibrium $W\in\{0,u_1^i\,,\,u_1^{ii}\}$.

Now assume $\Gamma:W_-\leadsto W_+$ is heteroclinic, for any $W_\pm\in \{0,u_1^i\,,\,u_1^{ii}\}$.
Then $R$-invariance of $W_\pm$\,, and time reversibility of $\Gamma$ under $R$, imply that the heteroclinic orbit $R\,\Gamma:W_+\leadsto W_-$ runs in the opposite direction.
Similar observations, with or without spatial reflection, apply to complex equilibria for general real nonlinearities $\overline f(w) = f(\overline w)$ replacing purely quadratic $f(w)=w^2$.

For complex entire nonlinearities, intriguing questions concern blow-up of complex time extensions $\Gamma(s-\mi r)$, alias solutions $\Gamma(r+\mi s)$ of nonlinear heat equations like \eqref{heat}.
For the purely quadratic heat equation \eqref{wtheta} with $\lambda=0$ and some inclined time rays $t=r\exp(\mi\theta), \ r\geq0$, the spectrum and complex two-dimensional unstable manifold of $W_1$ are beautifully explored in \cite{JaquetteStuke}.
Inclinations $\theta$ are chosen in $\{0,1,2\}\cdot\pi/4$.
The complex eigenspace of one of the two complex conjugate unstable eigenvalues (at $\theta=0$) is probed for computationally rigorous heteroclinic orbits $\Gamma:W_1\leadsto 0$, in angular rays of $1^\circ$ steps.
Results document existence, as well as absence, of $\Gamma$, depending on time inclination $\theta$ and angular degree of departure from $W_1$\,.
For $\theta=0$ or $\pi/4$, blow-up in $r$ is observed in the unstable manifold of $W_1$\,.
Such advanced insights and inspiring results concerning non-real complex equilibria, and the cases $\lambda\leq 0$, certainly lead beyond our present focus on real equilibria.
Several of our undue present constraints have yet to be overcome, like hyperbolicity, fast unstable manifolds, non-resonance, and asymptotically stable targets.
It may require and deserve substantial further analysis to progress towards such more complex directions, with all due diligence and generality.

\subsection{Periodic traveling waves in complex time}\label{Trav}

In this section, we start from the parabolic PDE under spatially $p$-periodic boundary conditions:
\begin{equation}
\label{heatS1}
w_r=w_{xx}+w^2-1\,,\qquad x\in \mathbb{S}^1=\mathbb{R}/p\mathbb{Z}\,.
\end{equation}
Up to positive scaling, this addresses the quadratic case \eqref{heat}; generalizations to complex entire ``pendulum'' nonlinearities $f=f(w)$ are evident.

\subsubsection{Real time}\label{Travreal}
On the whole real line $x\in\mathbb{R}$, \emph{real traveling wave solutions} $w(r,x)=W(\xi)\in\mathbb{R}$, in the real traveling wave coordinate
\begin{equation}
\label{xic}
\xi=x-cr,
\end{equation}
have been studied, ever since pioneering work by Fisher, Kolmogorov and coworkers \cite{Fisher, KPP}.
For a survey of some further developments in real settings, much beyond this first, standard variant and including PDE stability, see for example \cite{SandstedeTravel} and the references there.
Denoting $'=\tfrac{d}{d\xi}$, traveling waves $w=W(\xi)$ arise from eternal solutions of the ODE
\begin{equation}
\label{tw}
W''+cW'+W^2-1 = 0\,.
\end{equation}
For positive wave speeds $c>0$ and real $W,\ \xi$, this is a damped pendulum equation.
All eternal real non-equilibrium solutions are then heteroclinic ``backs'' $W=\Gamma: -1\leadsto +1$.

A family of real periodic traveling waves $W$ of \eqref{tw}, in particular, can only occur in their stationary, standing wave guise $c=0$; see also \eqref{heat0}.
The homoclinic orbit $W=\Gamma(\xi)$ of \eqref{tw}, alias a standing soliton pulse, arises in the limit of real period $p\nearrow\infty$; explicitly
\begin{equation}
\label{pulse}
\Gamma(\xi)= -1+3/\cosh^2(\xi/\sqrt{2})\,.
\end{equation}

\subsubsection{Imaginary time, periodic}\label{Travper}
By now, of course, we routinely extend the pendulum ODE \eqref{tw}  to complex ``time'' $\zeta=\xi+\mi\eta,\ '=\tfrac{d}{d\zeta}$\,, and to complex $W=W(\zeta)$.
As our second variant, we therefore consider \eqref{tw} in  imaginary ``time'' $\zeta=\mi\eta$. 
Define $\psi(s,x):=-W(\mi\eta)$ in the imaginary part
\begin{equation}
\label{etac}
\eta=x-cs
\end{equation}
of the traveling wave coordinate $\zeta$.
Imaginary periods $\mi p$ of $W(\zeta)$ then provide complex solutions of the scaled Schrödinger variant of \eqref{heatS1}, \eqref{psi}:
\begin{equation}
\label{psiS1}
\mi\psi_s=\psi_{xx}+\psi^2-1\,,\qquad x\in \mathbb{S}^1=\mathbb{R}/p\mathbb{Z}\,.
\end{equation}
For example, consider the explicit standing wave \eqref{pulse} for $c=0$, with complex $\zeta=\xi+\mi \eta$ replacing $\xi$. 
Then $\eta\mapsto\Gamma(\xi+\mi \eta)$ generates a spatially periodic standing wave of the Schrödinger equation \eqref{psiS1}, for any fixed real $\xi\neq 0$.
Indeed, \eqref{etac} provides a minimal spatial Schrödinger period of $p=\pi\sqrt{2}$\,.

\subsubsection{Imaginary time, quasi-periodic}\label{Travqper}
For arbitrary wave speeds $c>0$, i.e. for positive damping in \eqref{tw}, the real homoclinic orbit \eqref{pulse} becomes real heteroclinic, $\Gamma:-1\leadsto+1$, from the homogeneous saddle $W\equiv-1$ to the homogeneous sink $W\equiv+1$. But locally, $\Gamma$ is still given by the increasing part of the one-dimensional unstable manifold $W^u$ of $W\equiv-1$.
Let 
\begin{equation}
\label{twmu-}
\mu_-= \tfrac{1}{2}(-c+\sqrt{c^2+8}) >0
\end{equation} 
denote the unstable real eigenvalue at $W\equiv-1$.
Poincaré linearization at $W\equiv-1$, analogously to \eqref{Gammaperiod}, determines the spatial period $p$ in \eqref{etac}.
Explicitly,
\begin{equation}
\label{imp}
p=2\pi/\mu_-=4\pi/(-c+\sqrt{c^2+8}),
\end{equation}
for all $c\geq0$.
This holds near $W\equiv-1$, i.e. for real parts $\xi$ of the analytic extension $\Gamma=\Gamma(\xi+\mi\eta)$ which are fixed sufficiently negative, say $\xi\leq0$.

Fix any wave speed $c>0$, and consider $\xi$ as a parameter.
Then the resulting Schrödinger family $\psi(s,x)=-\Gamma(\xi+\mi(x-cs))$ of rotating waves possesses fixed spatial period $p$ and time period $p/c$. 
With respect to $\xi$, the family $\psi(s,x)$ is parametrized  via the conjugating flow of the traveling wave equation \eqref{tw} in real ``time'' $\xi\leq0$.
The parametrization of the rotating wave family reaches from small oscillations around $W\equiv-1$, at $\xi\searrow-\infty$, up to the first blow-up point, at $\xi\nearrow\min\xi^*(\eta)$.
These arguments are analogous to sections \ref{Fastest} and \ref{HeaSchro}. 
% and \ref{ODEd}.

Specifically, the two eigenvalues $\mu_+$ at the heteroclinic target $W\equiv+1$ of \eqref{tw} are
\begin{equation}
\label{twmu+}
\mu_+= \tfrac{1}{2}(-c\pm\sqrt{c^2-8})\,.
\end{equation} 
For $c>2\sqrt{2}$, both eigenvalues are simple, real, and stable.
They are non-resonant, except for the discrete sequence of 1:$m$ resonant wave speeds $c_m=\sqrt{2}(\sqrt{m}+\sqrt{1/m})\nearrow\infty$.
Non-resonant local Poincaré linearization \eqref{P} at the $\Gamma$-heteroclinic target $W\equiv+1$  then further aggravates the closure problem of the unstable manifold $\Gamma=W^u\setminus\{-1\}$ at the stable target, as we already mentioned in section \ref{Boundary}.

In fact, suppose the two stable eigenvalues $\mu_+$ are non-resonant, and the heteroclinic orbit $\Gamma: -1\leadsto+1$ does not linearize to only one of the  two exact eigenmodes of the two stable eigenvalues $\mu_+<0$.
This only discards a discrete set of exceptional wave speeds $c$.
For any other wave speed, $W(\zeta)$ near $W_+$ turns out to be 2-frequency quasi-periodic in imaginary time $\zeta=\mi\eta$, rather than periodic.
Quasi-periodicity of $\Gamma(\xi+\mi\eta)$ at the target $W\equiv+1$, in imaginary time $\mi\eta$, implies density in a family of 2-tori, parametrized over real times $\xi$, and shrinking to the target equilibrium $W\equiv+1$ for $\xi\nearrow+\infty$.

\subsubsection{Complex-valued blow-up}\label{Travblow}
For a fourth variant, we return to the parabolic PDE \eqref{heatS1} in real time $r$, but for complex-valued $w$.
Recall \cite{Yanagida} for blow-up under Neumann boundary conditions in $x$.
Periodic boundary conditions, instead, allow us to mix real and imaginary ``time'' $\zeta=\xi+\mi\eta$ in \eqref{tw} and consider
\begin{equation}
\label{zetac}
\zeta=cr+\mi x.
\end{equation}
Again we fix $c>0$.
The complex extension of the real heteroclinic orbit $\Gamma:-1\leadsto+1$ of \eqref{tw} then provides a solution $w(r,x)=-\Gamma(cr+\mi x)$ of the parabolic PDE \eqref{heatS1}.
By construction and Poincaré linearization, $\Gamma$ is in the ODE local unstable manifold $W^u$ of $W_-\equiv-1$, for large negative $r$.
By $p$-periodicity in $x=\Im \zeta$ near $W_-$\,, with $p$ from \eqref{imp}, section \ref{HeaSchro} then asserts blow-up of the complex solution $w$ in finite real time $r\nearrow \tfrac{1}{c}\min\xi^*(\eta)$.
In the present setting, this blow-up is governed by the damped pendulum ODE \eqref{tw}.
It may be interesting to pursue this construction further, and to compare the results with the approach of  \cite{Yanagida}.

\subsection{Ultra-invisible chaos and ultra-sharp Arnold tongues: a 1000 \euro\ question}\label{1000}

Even in pure ODE settings, complex time extensions of real analytic heteroclinic and homoclinic orbits $\Gamma$ hold interesting promise and applied interest.
We briefly sketch some recent progress.
For details and many further references, we refer to \cite{FiedlerScheurle, fiedlerClaudia}.
To be specific, consider autonomous ODE systems
\begin{equation}
\label{odef}
\dot{w}=f(\lambda,w).
\end{equation}
The nonlinearity $f$ is assumed to be analytic for small real $\lambda$, all $w\in\mathbb{C}^N$, and real for real $w$. 
Throughout, we assume a hyperbolic trivial equilibrium $0=f(\lambda,0)$.

\begin{figure}[t]
\centering \includegraphics[width=0.6\textwidth]{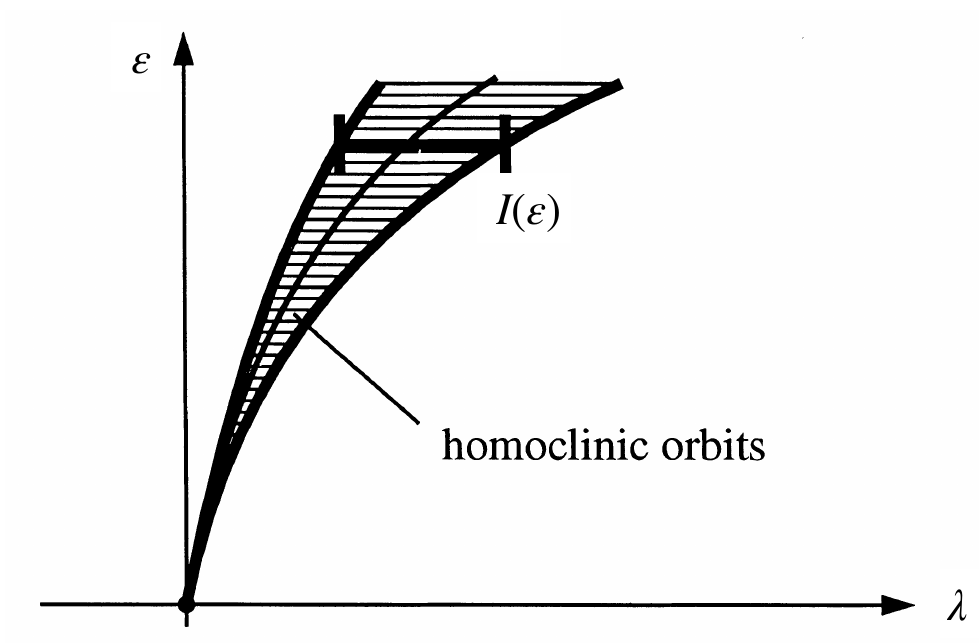}
\caption{\emph{
Schematic splitting region for a real homoclinic orbit $\Gamma$ of the ODE-flow \eqref{odef} (hashed), under discretization with step size $\eps$ or, equivalently, under $\eps$-periodic non-autonomous forcing; see \cite{FiedlerScheurle}.
At fixed levels of $\eps$, the horizontal splitting intervals $\lambda\in I(\eps)$ mark parameters $\lambda$ for which single-round homoclinic orbits occur near $\Gamma$.
For analyticity of $\Gamma(r+\mi s)$ in a strip $|s|\leq a$,  the width $\ell(\eps)$ of the splitting interval $I(\eps)$ is exponentially small in $\eps$ with exponent $a$; see \eqref{split}.
For better visibility, the horizontal width $\ell(\eps)$ of the exponentially flat splitting region has therefore been much exaggerated, in our schematic illustration.
We call the dynamics in the resulting chaotic region ``invisible chaos''.
For complex entire homoclinic orbits $\Gamma(t)$, the exponent $a$ could be chosen arbitrarily large.
Such ultra-exponentially small splittings would lead to ultra-invisible chaos.\\
}}
\label{fig4}
\end{figure}

Consider the fate of a real homoclinic orbit $\Gamma:0\leadsto 0$ of \eqref{odef} under analytic discretizations with small step size $\eps>0$ or, equivalently, under the stroboscopic time-$\eps$ map of nonautonomous forcings with rapid period $\eps$.
Under certain nondegeneracy conditions on $\Gamma$ and the discretization or forcing, this leads to \emph{invisible chaos}: the chaotic region accompanying transverse homoclinics is exponentially thin, both, in parameter space $(\lambda,\eps)$ and along the homoclinic loop.
More precisely, let 
\begin{equation}
\label{astrip}
r\in\mathbb{R},\ |s|\leq a
\end{equation} 
denote a horizontal complex strip where the complex time extension $\Gamma(r+\mi s)$ of $\Gamma(r)$ remains complex analytic.
Then the width $\ell(\eps)$ of the splitting interval $I(\eps)$, where transverse homoclinicity and accompanying shift dynamics prevail, satisfies an exponential splitting estimate
\begin{equation}
\label{split}
\ell(\eps) \leq C(a) \exp(- a/\eps),
\end{equation}
for some constant $C(a)>0$ and all small $\eps \searrow 0$.
Remarkably, the splitting region is small of infinite order $\exp(-a/\eps)$ in $\eps$, even under forcings or discretizations of only finite order $\eps^p$.
See figure \ref{fig4} for illustration.
For thorough ODE surveys on this topic, dating back as far as Poincaré \cite{Poincare3body} and starting technically with \cite{Neishtadt}, see \cite{Gelfreich01,Gelfreich02}.
For PDE extensions see \cite{MatthiesDiss,Matthieshom,MatthiesScheel} and the references there.
In numerical analysis, this topic runs under the name of \emph{backward error analysis}; see for example \cite{MoserPhysics, Reich, WulffOliver} and the more detailed references in \cite{fiedlerClaudia}.

The significance of \emph{complex entire homoclinic orbits} $\Gamma$ is now obvious.
In the complex entire case, an upper estimate \eqref{split} would hold for horizontal complex strips \eqref{astrip} of any width $\pm a$. 
In fact, this would strengthen \eqref{split} to an ultra-exponential estimate
\begin{equation}
\label{ultrasplit}
-\log \ell(\eps) \geq c(1/\eps) >0.
\end{equation}
Here $c(\cdot)$ is a suitable convex function of unbounded positive slope.
We would call this phenomenon \emph{ultra-exponentially small splitting of separatrices}.
Under a further mild nondegeneracy assumption, such splitting would be accompanied by \emph{ultra-invisible chaos}.

Alas, can it actually happen?
In \cite{fiedlerClaudia}, the author has personally offered a
\begin{center}
\textbf{1,000\,\euro\ reward}
\end{center}
for settling this question.
A \emph{negative answer} would prove that such complex entire homoclinic orbits $\Gamma(t)$ cannot exist.
In \cite{fiedlerClaudia}, we have already collected some very partial and incomplete results in that direction. 
However, cases involving complex conjugate, or even just algebraically degenerate real, eigenvalues remain wide open.
This includes the celebrated Shilnikov homoclinic orbits \cite{Shilnikov}.

A \emph{positive answer} would require any explicit example of a real nonstationary homoclinic orbit $\Gamma(t)$ to a hyperbolic equilibrium, for a complex entire ODE \eqref{odef} on $X=\mathbb{C}^N$ with $f$ real for real arguments, such that  $\Gamma(t)$ is complex entire for $t\in\mathbb{C}$. 
Such an example would initiate many challenges, not least for numerical exploration.
Eventually, it might lead towards a whole new theory to address ultra-exponentially small splitting behavior under discretization.

In order to protect well-established colleagues against their own, potentially intensely distracting, financial interests, the prize will be awarded to the first solution by anyone \emph{without} a permanent position in academia.
Priority is defined by submission time stamp at arxiv.org or equivalent repositories.
Subsequent confirmation by regular refereed publication is required.

We conclude with the time reversible, second-order, scalar ODE
\begin{equation}
\label{oderev}
\ddot w +\dot w^2 +w^2-3w=0.
\end{equation}
from \cite{fiedlerClaudia}.
By integrability, the homoclinic orbit $w(t)=\Gamma(t)$ to the hyperbolic trivial equilibrium $w=0$ satisfies
\begin{equation}
\label{revint}
\tfrac{1}{2}\dot{w}^2=\exp(-2w)-1+2w\,-\,\tfrac{1}{2}w^2\,.
\end{equation}
In particular, $\Gamma(t)$ is neither entire nor meromorphic in complex time.
Indeed, theorem 1.1 in \cite{fiedlerClaudia} implies that $\Gamma$ cannot be complex entire.
The term $\exp(-2w)$ in ODE \eqref{revint}, on the other hand, prevents meromorphic singularities.
If the singularity is isolated, it will therefore have to be essential.

Notably, however, the nonlinear ODE \eqref{oderev} does possess a \emph{complex entire $2\pi$-periodic orbit} $w$ given by the harmonic function
\begin{equation}
\label{revper}
w(t) := 2+\sqrt{2}\cos t\,.
\end{equation}
Under discretizations, alias rapid forcings, such complex entire periodic orbits provide devil's staircases with ultra-exponentially thin resonance plateaus, and correspondingly ultra-sharp Arnold tongues.
See \cite{ArnoldODE,Devaney,GuckenheimerHolmes} for a general background.
Further details and references are discussed in section 7 of \cite{fiedlerClaudia}.

Alas, we have only been able to derive an ultra-exponential \emph{upper} estimate \eqref{ultrasplit} for Arnold tongues associated to complex entire periodic orbits. 
This raises the question: how ultra-sharp are they, really?
Complementary \emph{lower} estimates of ultra-sharp Arnold tongues are not known.
For homoclinic, rather than periodic, splittings we even lack any 1,000\ \euro\ ODE example with an ultra-exponential upper estimate.

%%%%%%%%%%%%%%%%%%%%%%%%%%%%%%%%%%%%%%%%%%%%%%%%%%%%%%%%%%%
%\newpage

\end{document}